\documentclass{amsart}

\copyrightinfo{2000}{American Mathematical Society}

\newtheorem{theorem}{Theorem}[section]
\newtheorem{lemma}[theorem]{Lemma}

\newtheorem{proposition}[theorem]{Proposition}

\theoremstyle{definition}
\newtheorem{definition}[theorem]{Definition}

\theoremstyle{remark}
\newtheorem{remark}[theorem]{Remark}

\numberwithin{equation}{section}

\usepackage{amssymb}
\usepackage{longtable}

\usepackage[dvips]{graphicx}

\DeclareMathOperator{\rank}{rank}
\DeclareMathOperator{\supp}{supp}

\DeclareMathOperator{\Sym}{Sym}

\newcommand{\C}{\mathbb C}

\newcounter{primitivenumber}
\newcommand{\atcpntpn}{\addtocounter{primitivenumber}{1}\theprimitivenumber. }

\newsavebox{\aone}
\newsavebox{\aprime}
\newsavebox{\GreyCircle}
\newsavebox{\segm}
\newsavebox{\susp}
\newsavebox{\shortsusp}
\newsavebox{\bifurc}
\newsavebox{\longbifurc}
\newsavebox{\dthree}
\newsavebox{\dn}
\newsavebox{\shortdn}
\newsavebox{\atwo}
\newsavebox{\tosw}
\newsavebox{\tose}
\newsavebox{\tonw}
\newsavebox{\tone}
\newsavebox{\toe}
\newsavebox{\tow}
\newsavebox{\longam}
\newsavebox{\mediumam}
\newsavebox{\shortam}
\newsavebox{\plusaoneaone}
\newsavebox{\plusdm}

\savebox{\aone}{\begin{picture}(600,1800)\put(300,900){\circle*{150}}\put(300,1500){\circle{600}}\put(300,300){\circle{600}}\end{picture}}
\savebox{\aprime}{\begin{picture}(600,900)\put(300,900){\circle*{150}}\put(300,300){\circle{600}}\end{picture}}
\savebox{\GreyCircle}(600,600){\begin{picture}(600,600)(300,300)\put(600,600){\circle{600}}\put(30,25){\multiput(500,350)(150,0){2}{\circle*{70}}\multiput(425,425)(150,0){3}{\circle*{70}}\multiput(350,500)(150,0){4}{\circle*{70}}\multiput(425,575)(150,0){3}{\circle*{70}}\multiput(350,650)(150,0){4}{\circle*{70}}\multiput(425,725)(150,0){3}{\circle*{70}}\multiput(500,800)(150,0){2}{\circle*{70}}}\end{picture}}
\savebox{\segm}{\begin{picture}(1800,0)\multiput(0,0)(1800,0){2}{\circle*{150}}\thicklines\put(0,0){\line(1,0){1800}}\end{picture}}
\savebox{\susp}{\begin{picture}(3600,0)\multiput(0,0)(3600,0){2}{\circle*{150}}\thicklines\multiput(0,0)(2500,0){2}{\line(1,0){1100}}\multiput(1300,0)(400,0){3}{\line(1,0){200}}\end{picture}}
\savebox{\shortsusp}{\begin{picture}(1800,0)\multiput(0,0)(1800,0){2}{\circle*{150}}\thicklines\multiput(0,0)(400,0){5}{\line(1,0){200}}\end{picture}}
\savebox{\bifurc}{\begin{picture}(1200,2400)\multiput(1200,0)(0,2400){2}{\circle*{150}}\thicklines\put(0,1200){\line(1,1){1200}}\put(0,1200){\line(1,-1){1200}}\end{picture}}
\savebox{\longbifurc}{\begin{picture}(1800,3600)\multiput(1800,0)(0,3600){2}{\circle*{150}}\thicklines\put(0,1800){\line(1,1){1800}}\put(0,1800){\line(1,-1){1800}}\end{picture}}
\savebox{\dthree}{\begin{picture}(3600,600)\put(1500,0){\usebox{\GreyCircle}}\multiput(0,300)(1800,0){2}{\usebox{\segm}}\end{picture}}
\savebox{\dn}{\begin{picture}(8700,2400)\put(0,900){\usebox{\GreyCircle}}\multiput(300,1200)(5400,0){2}{\usebox{\segm}}\put(2100,1200){\usebox{\susp}}\put(7500,0){\usebox{\bifurc}}\end{picture}}
\savebox{\shortdn}{\begin{picture}(6900,2400)\put(0,900){\usebox{\GreyCircle}}\multiput(300,1200)(3600,0){2}{\usebox{\segm}}\put(2100,1200){\usebox{\shortsusp}}\put(5700,0){\usebox{\bifurc}}\end{picture}}
\savebox{\atwo}{\begin{picture}(2400,600)\put(300,300){\usebox{\segm}}\multiput(300,300)(1800,0){2}{\circle{600}}\multiput(300,300)(25,25){13}{\circle*{70}}\multiput(600,600)(300,0){4}{\multiput(0,0)(25,-25){7}{\circle*{70}}}\multiput(750,450)(300,0){4}{\multiput(0,0)(25,25){7}{\circle*{70}}}\multiput(2100,300)(-25,25){13}{\circle*{70}}\end{picture}}
\savebox{\tosw}{\begin{picture}(300,300)\multiput(-550,-550)(10,30){16}{\line(0,1){30}}\multiput(-550,-550)(30,10){16}{\line(1,0){30}}\end{picture}}
\savebox{\tose}{\begin{picture}(300,300)\multiput(550,-550)(-10,30){16}{\line(0,1){30}}\multiput(550,-550)(-30,10){16}{\line(-1,0){30}}\end{picture}}
\savebox{\tonw}{\begin{picture}(300,300)\multiput(-550,550)(10,-30){16}{\line(0,-1){30}}\multiput(-550,550)(30,-10){16}{\line(1,0){30}}\end{picture}}
\savebox{\tone}{\begin{picture}(300,300)\multiput(550,550)(-10,-30){16}{\line(0,-1){30}}\multiput(550,550)(-30,-10){16}{\line(-1,0){30}}\end{picture}}
\savebox{\toe}{\begin{picture}(400,400)\multiput(0,0)(20,10){21}{\line(1,0){20}}\multiput(0,400)(20,-10){21}{\line(1,0){20}}\end{picture}}
\savebox{\tow}{\begin{picture}(400,400)\multiput(400,0)(-20,10){21}{\line(-1,0){20}}\multiput(400,400)(-20,-10){21}{\line(-1,0){20}}\end{picture}}
\savebox{\longam}{\begin{picture}(7800,600)\put(300,300){\circle*{150}}\put(2100,300){\circle*{150}}\put(5700,300){\circle*{150}}\put(7500,300){\circle*{150}}\put(300,300){\circle{600}}\put(7500,300){\circle{600}}\multiput(300,300)(25,25){13}{\circle*{70}}\multiput(600,600)(3900,0){2}{\multiput(0,0)(300,0){9}{\multiput(0,0)(25,-25){7}{\circle*{70}}}\multiput(150,-150)(300,0){9}{\multiput(0,0)(25,25){7}{\circle*{70}}}}\multiput(7500,300)(-25,25){13}{\circle*{70}}\thicklines\put(300,300){\line( 1, 0){2250}}\put(7500,300){\line(-1, 0){2250}}\multiput(2850,300)(600,0){4}{\line( 1, 0){300}}\end{picture}}
\savebox{\mediumam}{\begin{picture}(6000,600)\multiput(300,300)(3600,0){2}{\usebox{\segm}}\put(2100,300){\usebox{\shortsusp}}\multiput(300,300)(1800,0){4}{\circle*{150}}\multiput(300,300)(5400,0){2}{\circle{600}}\multiput(300,300)(25,25){13}{\circle*{70}}\multiput(600,600)(3000,0){2}{\multiput(0,0)(300,0){6}{\multiput(0,0)(25,-25){7}{\circle*{70}}}\multiput(150,-150)(300,0){6}{\multiput(0,0)(25,25){7}{\circle*{70}}}}\multiput(5700,300)(-25,25){13}{\circle*{70}}\multiput(2700,600)(300,0){3}{\circle*{70}}\end{picture}}
\savebox{\shortam}{\begin{picture}(4200,600)\put(300,300){\usebox{\susp}}\multiput(300,300)(3600,0){2}{\circle{600}}\multiput(300,300)(25,25){13}{\circle*{70}}\multiput(600,600)(2100,0){2}{\multiput(0,0)(300,0){3}{\multiput(0,0)(25,-25){7}{\circle*{70}}}\multiput(150,-150)(300,0){3}{\multiput(0,0)(25,25){7}{\circle*{70}}}}\multiput(3900,300)(-25,25){13}{\circle*{70}}\multiput(1800,600)(300,0){3}{\circle*{70}}\end{picture}}
\savebox{\plusaoneaone}{\begin{picture}(1950,3000)\put(0,300){\usebox{\bifurc}}\multiput(1200,300)(0,2400){2}{\circle{600}}\multiput(1500,300)(0,2400){2}{\line(1,0){450}}\put(1950,300){\line(0,1){2400}}\end{picture}}
\savebox{\plusdm}{\begin{picture}(6600,2400)\put(0,1200){\usebox{\segm}}\put(1800,1200){\usebox{\susp}}\put(5400,0){\usebox{\bifurc}}\put(1500,900){\usebox{\GreyCircle}}\end{picture}}

\begin{document}

\title{Wonderful varieties of type $D$}

\author{Paolo Bravi}
\address{Dipartimento di Matematica\\ Universit\`a La Sapienza\\ P.le Aldo Moro 2\\ 00185 Roma\\ Italy}
\curraddr{Fachbereich Mathematik\\ Universit\"at Wuppertal\\ Gauss-Str. 20\\ 42097 Wuppertal\\ Germany}
\email{bravi@math.uni-wuppertal.de}

\author{Guido Pezzini}
\address{Dipartimento di Matematica\\ Universit\`a La Sapienza\\ P.le Aldo Moro 2\\ 00185 Roma\\ Italy}
\curraddr{Institut Fourier\\ Universit\'e Joseph Fourier\\ B.P. 74\\ 38402 Saint-Martin d'H\`eres\\ France}
\email{pezzini@mozart.ujf-grenoble.fr}

\date{}



\begin{abstract}
Let $G$ be a connected semisimple group over $\C$, whose simple components have type $\mathsf A$ or $\mathsf D$.
We prove that wonderful $G$-varieties are classified by means of combinatorial objects
called {\em spherical systems}. This is a generalization of a known result of Luna for groups of type $\mathsf A$;
thanks to another result of Luna, this implies also the classification
of all spherical $G$-varieties for the groups $G$ we are considering. For these $G$ we also
prove the smoothness of the embedding of Demazure.
\end{abstract}

\maketitle


\section{Introduction}
Let $G$ be a connected semisimple group over $\C$. In \cite{L01} it is proven that if $G$ is adjoint
of type $\mathsf A$ then wonderful $G$-varieties are classified by means of discrete invariants coming
from the theory of spherical varieties. These invariants are treated as combinatorial objects
called {\em spherical systems}.

In this work we extend this result to the case where $G$ has simple components of type $\mathsf A$ or $\mathsf D$.

Wonderful varieties were studied at first in the theory of linear algebraic groups as compactifications of
certain homogeneous spaces, precisely the symmetric ones (for more details, see \cite{DP83}). They also appear in the theory of
spherical embeddings, as ``canonical'' compactifications of particular spherical homogeneous spaces.
Actually all wonderful varieties can be obtained in this way (\cite{L96}).

Wonderful varieties of low rank have been classified for all $G$, by Ahiezer (rank $1$, in the framework
of compact groups, also obtained with algebraic methods by Brion) and Wasserman (rank $2$). Their results
have been crucial in establishing the axioms defining the spherical systems.

The role of wonderful varieties within the theory of spherical varieties has been made clear in \cite{Kn96} and
\cite{L01}; their importance has grown up to the point that the classification obtained in \cite{L01} and in the present work
implies the classification of all spherical $G$-varieties for the groups $G$ under consideration.

We also prove for these groups a known conjecture of Brion (appeared in \cite{B90}), already proven by Luna in \cite{L02} for $G$ of
type $\mathsf A$, which states that the embedding of Demazure is smooth.

We are grateful to D.~Luna and C.~Procesi for their precious help. 


\section{The classification of wonderful varieties}

\subsection{Wonderful varieties}

We recall some basic definitions and known facts about wonderful and spherical varieties. Let $G$ be a complex semisimple linear algebraic group, $B$ and $B^-$ two mutually opposite Borel subgroups, and $T=B\cap B^-$ a maximal torus.

\begin{definition}[\cite{DP83}] An algebraic $G$-variety $X$ is \textit{wonderful} of rank $r$ if: \begin{itemize} \item[-] $X$ is smooth and complete, \item[-] $G$ has a dense orbit in $X$ whose complement is the union of $r$ smooth prime divisors $D_i$, $i=1,\ldots,r$, with normal crossings, \item[-] the intersection of the divisors $D_i$ is nonempty and, for all $I\subseteq\{1,\ldots,r\}$, \[(\bigcap_{i\in I}D_i)\setminus(\bigcup_{i\notin I}D_i)\] is a $G$-orbit. \end{itemize} \end{definition}

A wonderful $G$-variety is always projective and \textit{spherical} (see \cite{L96}), where in general an irreducible $G$-variety $X$ is said to be spherical if it is normal and it contains an open dense $B$-orbit. This implies that any $B$-proper rational function on $X$ is uniquely determined (up to a scalar factor) by its $B$-weight; the lattice of such weights is called $\Xi_X$. If $X$ is wonderful of rank $r$ then $\Xi_X$ has rank $r$ (\cite{L96}).

Every $\mathbb Q$-valued discrete valuation $\nu$ of $\mathbb C(X)$ defines a functional $\rho_X(\nu)$ on $\Xi_X$ in the following way: $\langle\rho_X(\nu),\chi\rangle=\nu(f_\chi)$ where $\chi\in\Xi_X$ and $f_\chi\in\mathbb C(X)$ is $B$-proper with weight $\chi$. The application $\nu\mapsto\rho_X(\nu)$ is injective if we take its restriction to the set $\mathcal V_{X}$ of all $G$-invariant valuations of $X$. In this way $\mathcal V_{X}$ is identified with a subset of $\mathrm{Hom}_{\mathbb Z} (\Xi_X,\mathbb Q)$, and it turns out to be a polyhedral convex cone.

This cone can be described as the set $\left\{v\in\mathrm{Hom}_{\mathbb Z} (\Xi_X,\mathbb Q)\colon v(\chi)\leq 0\;\;\forall \chi\in\Sigma\right\}$ for a uniquely determined set $\Sigma$ of indecomposable elements of $\Xi_X$ called the \textit{spherical roots} of $X$.

If a prime divisor of $X$ is $B$-stable but not $G$-stable, then it is called a \textit{colour}; the set of all colours is denoted by $\Delta_X$. The discrete valuation $\nu_D$ associated to a colour $D$ defines a functional $\rho_X(\nu_D)$ which is also denoted by  $\rho_X(D)$ for brevity. In general the application $D\mapsto\rho_X(D)$ (where $D$ ranges in $\Delta_X$) is not injective. The stabilizer in $G$ of the union of all colours is a parabolic subgroup containing $B$ denoted by $P_X$.

A subgroup $H$ of $G$ is \textit{spherical} if $G/H$ is a spherical $G$-variety under the action of left multiplication, and it is \textit{wonderful} if $G/H$ admits an embedding which is a wonderful variety. If $H$ is wonderful then $H$ necessarily has finite index in $N_G(H)$ and the wonderful embedding of $G/H$ is unique. Moreover, if a spherical subgroup $H$ is equal to $N_G(H)$ then it is wonderful (\cite{Kn96}).

\subsection{Relations between spherical roots and colours}

If $X$ is a wonderful variety, the set of its spherical roots $\Sigma_X$ is a basis of $\Xi_X$. In this setting one can also retrieve 
$\Sigma_X$ in the following 
way: let $z\in X$ be the unique fixed point of $B^-$ and consider the orbit $Z=G.z$ which is the unique closed orbit in $X$. Then the spherical roots are exactly the $T$-weights appearing in the quotient $T_z X/T_z Z$.

One can associate to each spherical root $\gamma$ a $G$-stable prime divisor $D^\gamma$ such that $\gamma$ is the $T$-weight of $T_zX/T_zD^\gamma$. Consider the intersection of all $G$-invariant prime divisors of $X$ different from $D^\gamma$: this intersection is a wonderful variety of rank $1$, having $\gamma$ as its spherical root.

The set of spherical roots of all wonderful $G$-varieties of rank $1$ is denoted by $\Sigma(G)$; our last considerations amount to say that $\Sigma_X\subset\Sigma(G)$ for any wonderful $G$-variety $X$. For $G$ of adjoint type the elements of $\Sigma(G)$ are always linear combinations of simple roots with nonnegative integer coefficients.

For any given $G$ there exist only finitely many wonderful varieties of rank $1$, and they are classified (\cite{A83}, \cite{B89}). With this regard it is convenient to recall the following.

\begin{definition} A wonderful variety $X$ with open $G$-orbit $X^\circ_G\cong G/H$ is said to be \textit{prime} if $H$ satisfies the following conditions: \begin{itemize} \item[-] if $P$ is a parabolic subgroup such that $P^r \subseteq H \subseteq P$ then $P=G$; \item[-] if $G=G_1 \times G_2$ and $H=H_1 \times H_2$ with $H_i \subseteq G_i$ then $G_1 =G$ or $G_2 =G$. \end{itemize} \end{definition}

See Section~\ref{SS31} for a motivation of this definition. We report in Table~\ref{T1} the list of rank one prime wonderful $G$-varieties for $G$ an adjoint group of type $\mathsf A$~$\mathsf D$, up to isomorphism. Two $G$-varieties, $X$ and $X^\prime$, are here considered to be isomorphic if there exists an outer automorphism $\sigma$ of $G$ such that $X$ is $G$-isomorphic to $X^\prime$ endowed with the action of $G$ through $\sigma$.

\textit{How to read Table~\ref{T1}}. For any rank one prime wonderful $G$-variety $X$ we report, in Column~4, the homogeneous space $G/H\cong X^\circ_G$ and, moreover, the following informations: (Column~1) the spherical root as linear combination of simple roots denoted as usually, (Column~2) the type of the spherical root (the type will be used throughout the text to refer to the corresponding spherical root according with this table), (Column~3) the Luna diagram of the spherical system of $X$.

\begin{table}
\begin{small}
\begin{longtable}{|p{2.7cm}|p{1cm}|c|p{4.5cm}|}
\caption{Rank one prime wonderful $G$-varieties for $G$ adjoint of type $\mathsf A$~$\mathsf D$, up to isomorphism.}
\label{T1}\\
\hline
Spherical root & Type & Diagram & $G$ / $H$ \\
\hline
\hline
\endfirsthead
\caption{Rank one prime wonderful $G$-varieties (continuation).}\\
\hline
Spherical root & Type & Diagram & $G$ / $H$ \\
\hline
\hline
\endhead
$\alpha_1$ & $a_1$ & \parbox[c]{600\unitlength}{\begin{picture}(600,2400)\put(0,300){\usebox{\aone}}\end{picture}} & $SL(2)$ / $GL(1)$ \\
\hline
$\alpha_1+\ldots+\alpha_n$ & $a_n$ & \parbox[c]{6000\unitlength}{\begin{picture}(6000,1200)\put(0,300){\usebox{\mediumam}}\end{picture}} & $SL(n+1)$ / $GL(n)$, $n\geq2$ \\
\hline
$2\alpha_1$ & $a_1^\prime$ & \parbox[c]{600\unitlength}{\begin{picture}(600,2400)\put(0,300){\usebox{\aprime}}\end{picture}} & $SL(2)$ / $N_G(GL(1))$ \\
\hline
$\alpha_1+\alpha_1^\prime$ & $a_1\times a_1$ (or $d_2$) & \parbox[c]{3000\unitlength}{\begin{picture}(3000,2100)(0,-300)\put(300,750){\circle{600}}\put(300,750){\circle*{150}}\put(2700,750){\circle{600}}\put(2700,750){\circle*{150}}\put(300,450){\line( 0,-1){450}}\put(300,0){\line( 1, 0){2400}}\put(2700,0){\line( 0, 1){450}}\end{picture}} & $(SL(2)\times SL(2))$ / $(SL(2)\cdot C_G)$, $SL(2)$ embedded diagonally \\
\hline
$\alpha_1+2\alpha_2+\alpha_3$ & $d_3$ & \parbox[c]{3600\unitlength}{\begin{picture}(3600,1200)\put(0,300){\usebox{\dthree}}\end{picture}} & $SL(4)$ / $(Sp(4)\cdot C_G)$ \\
\hline
\hline
$2\alpha_1+\ldots+2\alpha_{n-2}+\alpha_{n-1}+\alpha_n$ & $d_n$ & \parbox[c]{6900\unitlength}{\begin{picture}(6900,3000)\put(0,300){\usebox{\shortdn}}\end{picture}} & \mbox{$SO(2n)$ / $(SO(2n-1)\cdot C_G)$,} \mbox{$n\geq4$} \\
\hline
\end{longtable}
\end{small}
\end{table}

Rank two prime wonderful varieties have been classified by Wasserman (\cite{Wa96}); in Table~\ref{T2} we report the ones for adjoint groups of type $\mathsf A$~$\mathsf D$ (up to isomorphism, as in Table~\ref{T1}).

\textit{How to read Table~\ref{T2}}. For any rank two prime wonderful $G$-variety $X$ we report, in Column~3, the homogeneous space $G/H\cong X^\circ_G$ and, moreover: (Column~1) the label according with \cite{Wa96} and the Luna diagram of the spherical system, (Column~2) the set of spherical roots and the set of colours with the table of values taken by the functional associated to each colour on the spherical roots. When necessary, in Column~3, we describe the subgroup $H$ by providing its connected center $C^\circ$, the commutator subgroup $(L,L)$ of a Levi subgroup $L$ and the Lie algebra of its unipotent radical $H^u$ in a rather concise way: in a matrix we report the generic element $(c_1,\ldots,c_{\rank \Xi(H)})$ of $C^\circ$ times the generic element of $(L,L)$, plus the generic element of $Lie(H^u)$ delimited by a line. We denote with $\pi$ the usual isogeny $\pi\colon SL(2)\times SL(2)\to SO(4)$ and with $j$ the usual injection $j\colon SO(2n-1)\to SO(2n)$.

\begin{small}
\begin{longtable}{|p{4cm}|p{3.5cm}|p{4cm}|}
\caption{Rank two prime wonderful $G$-varieties for $G$ adjoint of type $\mathsf A$~$\mathsf D$, up to isomorphism.}
\label{T2} \\
\hline
Diagram & $\Sigma$, $\delta_{\alpha_i}(\gamma_j)$ & $G$ / $H$ ($C^\circ$, $(L,L)$, $Lie(H^u)$) \\
\hline
\endfirsthead
\caption{Rank two prime wonderful $G$-varieties (continuation).} \\
\hline
Diagram & $\Sigma$, $\delta_{\alpha_i}(\gamma_j)$ & $G$ / $H$ ($C^\circ$, $(L,L)$, $Lie(H^u)$) \\
\hline
\endhead
\hline
\endfoot
\hline
$A1$ \begin{center}
\begin{picture}(2400,900)\multiput(0,0)(1800,0){2}{\usebox{\aprime}}\thicklines\put(300,900){\line(1,0){1800}}\end{picture}
\end{center} &
\begin{scriptsize}
$\{2\alpha_1,\ 2\alpha_2\}$,
\vspace{1ex}\begin{center}$\begin{array}{r|cc}
\hline
\delta^\prime_{\alpha_1} & 2 & -1 \\
\delta^\prime_{\alpha_2} & -1 & 2
\end{array}$\end{center}
\end{scriptsize} &
$SL(3)$ / $(SO(3)\cdot C_G)$ \\
\hline
$A2$ \begin{center}
\begin{picture}(6600,1500)\multiput(300,750)(4200,0){2}{\multiput(0,0)(1800,0){2}{\circle*{150}}}\multiput(300,750)(4200,0){2}{\multiput(0,0)(1800,0){2}{\circle{600}}}\put(300,450){\line(0,-1){450}}\put(300,0){\line(1,0){4200}}\put(4500,0){\line(0,1){450}}\put(2100,1050){\line(0,1){450}}\put(2100,1500){\line(1,0){4200}}\put(6300,1500){\line(0,-1){450}}\thicklines\multiput(300,750)(4200,0){2}{\line(1,0){1800}}\end{picture}
\end{center} &
\begin{scriptsize}
$\{\alpha_1+\alpha_1^\prime,\ \alpha_2+\alpha_2^\prime\}$,
\vspace{1ex}\begin{center}$\begin{array}{r|cc}
\hline
\delta_{\alpha_1}=\delta_{\alpha_1^\prime} & 2 & -1 \\
\delta_{\alpha_2}=\delta_{\alpha_2^\prime} & -1 & 2
\end{array}$\end{center}
\end{scriptsize} &
\mbox{$(SL(3)\times SL(3))$ /} \mbox{$(SL(3)\cdot C_G)$}, \mbox{\begin{scriptsize}$SL(3)$ embedded diagonally\end{scriptsize}} \\
\hline
$A3$ \begin{center}
\begin{picture}(7200,600)\multiput(0,300)(1800,0){5}{\circle*{150}}\multiput(1500,0)(3600,0){2}{\usebox{\GreyCircle}}\thicklines\put(0,300){\line(1,0){7200}}\end{picture}
\end{center} &
\begin{scriptsize}
$\{\alpha_1+2\alpha_2+\alpha_3,\ \alpha_3+2\alpha_4+\alpha_5\}$
\vspace{1ex}\begin{center}$\begin{array}{r|cc}
\hline
\delta_{\alpha_2} & 2 & -1 \\
\delta_{\alpha_4} & -1 & 2
\end{array}$\end{center}
\end{scriptsize} &
$SL(6)$ / $(Sp(6)\cdot C_G)$ \\
\hline
$A4(i)$ \begin{center}
\begin{picture}(4200,2400)\put(1800,300){\usebox{\aone}}\put(1200,1600){\put(150,0){\usebox{\tow}}}
\multiput(300,1200)(3600,0){2}{\circle*{150}}\multiput(300,1200)(3600,0){2}{\circle{600}}\put(300,900){\line(0,-1){900}}\put(300,0){\line(1,0){3600}}\put(3900,0){\line(0,1){900}}\thicklines\put(300,1200){\line(1,0){3600}}\end{picture}
\end{center} &
\begin{scriptsize}
$\{\alpha_1+\alpha_3,\ \alpha_2\}$
\vspace{1ex}\begin{center}$\begin{array}{r|cc}
\hline
\delta_{\alpha_1}=\delta_{\alpha_3} & 2 & -1 \\
\delta_{\alpha_2}^+ & -1 & 1 \\
\delta_{\alpha_2}^- & -1 & 1
\end{array}$\end{center}
\end{scriptsize} &
$SL(4)$ / $S(GL(2)\times GL(2))$ \\
\hline
$A4(ii)$ \begin{center}
\begin{picture}(4200,2400)\put(1800,300){\usebox{\aprime}}\multiput(300,1200)(3600,0){2}{\circle*{150}}\multiput(300,1200)(3600,0){2}{\circle{600}}\put(300,900){\line(0,-1){900}}\put(300,0){\line(1,0){3600}}\put(3900,0){\line(0,1){900}}\thicklines\put(300,1200){\line(1,0){3600}}\end{picture}
\end{center} &
\begin{scriptsize}
$\{\alpha_1+\alpha_3,\ 2\alpha_2\}$
\vspace{1ex}\begin{center}$\begin{array}{r|cc}
\hline
\delta_{\alpha_1}=\delta_{\alpha_3} & 2 & -2 \\
\delta_{\alpha_2}^\prime & -1 & 2
\end{array}$\end{center}
\end{scriptsize} &
$SL(4)$ / $N_G(SL(2)\times SL(2))$ \\
\hline
$A4(iii)$ \begin{center}
\begin{picture}(9600,1800)\put(1800,600){\usebox{\mediumam}}\put(0,150){\begin{picture}(9600,1500)\multiput(300,750)(9000,0){2}{\circle*{150}}\multiput(300,750)(9000,0){2}{\circle{600}}\put(300,450){\line(0,-1){450}}\put(300,0){\line(1,0){9000}}\put(9300,0){\line(0,1){450}}\thicklines\multiput(300,750)(7200,0){2}{\line(1,0){1800}}\end{picture}}\end{picture}
\end{center} &
\begin{scriptsize}
$\{\alpha_1+\alpha_n,\ \alpha_2+\ldots+\alpha_{n-1}\}$
\vspace{1ex}\begin{center}$\begin{array}{r|cc}
\hline
\delta_{\alpha_1}=\delta_{\alpha_n} & 2 & -1 \\
\delta_{\alpha_2} & -1 & 1 \\
\delta_{\alpha_{n-1}} & -1 & 1
\end{array}$\end{center}
\end{scriptsize} &
\mbox{$SL(n+1)$ /} \mbox{$S(GL(2)\times GL(n-1))$} \\
\hline
$A5$ \begin{center}
\begin{picture}(4200,2400)\put(1800,300){\usebox{\aone}}
\multiput(300,1200)(3600,0){2}{\circle*{150}}\multiput(300,1200)(3600,0){2}{\circle{600}}\put(300,900){\line(0,-1){900}}\put(300,0){\line(1,0){3600}}\put(3900,0){\line(0,1){900}}\thicklines\put(300,1200){\line(1,0){3600}}\end{picture}
\end{center} &
\begin{scriptsize}
$\{\alpha_1+\alpha_3,\ \alpha_2\}$
\vspace{1ex}\begin{center}$\begin{array}{r|cc}
\hline
\delta_{\alpha_1}=\delta_{\alpha_3} & 2 & -1 \\
\delta_{\alpha_2}^+ & 0 & 1 \\
\delta_{\alpha_2}^- & -2 & 1
\end{array}$\end{center}
\end{scriptsize} &
$SL(4)$ / $H$,
\begin{scriptsize}
\vspace{1ex}\begin{center}$\left(\begin{array}{cc}
c\,A & 0 \\
\cline{1-1}
\multicolumn{1}{|c|}{M} & c^{-1}\,A \\
\cline{1-1}
\end{array}\right)$\end{center}
for $A\in SL(2)$, $c\in GL(1)$, $M=M^\prime+M^{\prime\prime}\in \mathfrak{sl}(2)\oplus\C$, $M^{\prime\prime}=0$
\end{scriptsize} \\
\hline
$A6(i)$ \begin{center}
\begin{picture}(2400,2100)\multiput(0,0)(1800,0){2}{\usebox{\aone}}
\thicklines\put(300,900){\line(1,0){1800}}\end{picture}
\end{center} &
\begin{scriptsize}
$\{\alpha_1,\ \alpha_2\}$
\vspace{1ex}\begin{center}$\begin{array}{r|cc}
\hline
\delta_{\alpha_1}^+ & 1 & 0 \\
\delta_{\alpha_1}^- & 1 & -1 \\
\delta_{\alpha_2}^+ & 0 & 1 \\
\delta_{\alpha_2}^- & -1 & 1
\end{array}$\end{center}
\end{scriptsize} &
$SL(3)$ / $H$,
\begin{scriptsize}
\begin{center}$\left(\begin{array}{ccc}
c_1 & 0 & 0 \\
0 & c_2 & 0 \\
\cline{1-2}
\multicolumn{1}{|c}{\ast} & \multicolumn{1}{c|}{0} & \ast \\
\cline{1-2}
\end{array}\right)$\end{center}
for $c_1,c_2\in GL(1)$
\end{scriptsize} \\
\hline
$A6(ii)$ \begin{center}
\begin{picture}(7800,2100)\put(300,900){\usebox{\segm}}\put(0,0){\usebox{\aone}}\put(1800,600){\usebox{\mediumam}}
\end{picture}
\end{center} &
\begin{scriptsize}
$\{\alpha_1,\ \alpha_2+\ldots+\alpha_n\}$
\vspace{1ex}\begin{center}$\begin{array}{r|cc}
\hline
\delta_{\alpha_1}^+ & 1 & 0 \\
\delta_{\alpha_1}^- & 1 & -1 \\
\delta_{\alpha_2} & -1 & 1 \\
\delta_{\alpha_n} & 0 & 1
\end{array}$\end{center}
\end{scriptsize} &
$SL(n+1)$ / $H$,
\begin{scriptsize}
\begin{center}$\left(\!\begin{array}{ccc}
c_1 & 0 & 0 \\
0 & \ast & 0 \\
\cline{1-2}
\multicolumn{1}{|c}{\ast} & \multicolumn{1}{c|}{0} & c_2\,A \\
\cline{1-2}
\end{array}\!\right)$\end{center}
for $A\in SL(n-1)$, $c_1,c_2\in GL(1)$
\end{scriptsize} \\
\hline
$A6(iii)$ \begin{center}
\begin{picture}(9600,1800)\multiput(0,600)(5400,0){2}{\usebox{\shortam}}\put(3900,900){\usebox{\segm}}\end{picture}
\end{center} &
\begin{scriptsize}
$\{\alpha_1+\ldots+\alpha_m,\ \alpha_{m+1}+\ldots+\alpha_n\}$
\vspace{1ex}\begin{center}$\begin{array}{r|cc}
\hline
\delta_{\alpha_1} & 1 & 0 \\
\delta_{\alpha_m} & 1 & -1 \\
\delta_{\alpha_{m+1}}^+ & -1 & 1 \\
\delta_{\alpha_n}^- & 0 & 1
\end{array}$\end{center}
\end{scriptsize} &
$SL(n+1)$ / $H$,
\begin{scriptsize}
\begin{center}$\left(\begin{array}{ccc}
c_1\,A_1 & 0 & 0 \\
0 & \ast & 0 \\
\cline{1-2}
\multicolumn{1}{|c}{\ast} & \multicolumn{1}{c|}{0} & c_2\,A_2 \\
\cline{1-2}
\end{array}\right)$\end{center}
for $A_1\in SL(m)$, $A_2\in SL(n-m-1)$, $c_1,c_2\in GL(1)$
\end{scriptsize} \\
\hline
$A7(i)$ \begin{center}
\begin{picture}(2400,2550)\put(300,900){\usebox{\segm}}\multiput(0,0)(1800,0){2}{\usebox{\aone}}\multiput(300,1800)(1800,0){2}{\line(0,1){450}}\put(300,2250){\line(1,0){1800}}
\end{picture}
\end{center} &
\begin{scriptsize}
$\{\alpha_1,\ \alpha_2\}$
\vspace{1ex}\begin{center}$\begin{array}{r|cc}
\hline
\delta_{\alpha_1}^+=\delta_{\alpha_2}^+ & 1 & 1 \\
\delta_{\alpha_1}^- & 1 & -2 \\
\delta_{\alpha_2}^- & -2 & 1
\end{array}$\end{center}
\end{scriptsize} &
$SL(3)$ / $H$,
\begin{scriptsize}
\begin{center}$\left(\begin{array}{ccc}
c & 0 & 0 \\
\cline{1-1}
\multicolumn{1}{|c|}{m} & 1 & 0 \\
\cline{2-2}
\multicolumn{1}{|c}{\ast} & \multicolumn{1}{c|}{-m} & c^{-1} \\
\cline{1-2}
\end{array}\right)$\end{center}
for $c \in GL(1)$, $m\in\C$
\end{scriptsize} \\
\hline
$A7(ii)$ \begin{center}
\begin{picture}(4200,1800)\multiput(0,600)(1800,0){2}{\usebox{\atwo}}\end{picture}
\end{center} &
\begin{scriptsize}
$\{\alpha_1+\alpha_2,\ \alpha_2+\alpha_3\}$
\vspace{1ex}\begin{center}$\begin{array}{r|cc}
\hline
\delta_{\alpha_1} & 1 & -1 \\
\delta_{\alpha_2} & 1 & 1 \\
\delta_{\alpha_3} & -1 & 1
\end{array}$\end{center}
\end{scriptsize} &
$SL(4)$ / $H$,
\begin{scriptsize}
\begin{center}$\left(\begin{array}{ccc}
c & 0 & 0 \\
\cline{1-1}
\multicolumn{1}{|c|}{M_1} & A & 0 \\
\cline{2-2}
\multicolumn{1}{|c}{\ast} & \multicolumn{1}{c|}{M_2} & c^{-1} \\
\cline{1-2}
\end{array}\right)$\end{center}
for $A\in SL(2)$, $c \in GL(1)$, $M_1,M_2\in\C^2\cong(\C^2)^\ast$, $M_1+M_2=0$
\end{scriptsize} \\
\hline
$A8$ \begin{center}
\begin{picture}(2400,2550)\multiput(0,0)(1800,0){2}{\usebox{\aone}}\multiput(300,1800)(1800,0){2}{\line(0,1){450}}\put(300,2250){\line(1,0){1800}}\end{picture}
\end{center} &
\begin{scriptsize}
$\{\alpha_1,\ \alpha_2\}$
\vspace{1ex}\begin{center}$\begin{array}{r|cc}
\hline
\delta_{\alpha_1}^+=\delta_{\alpha_2}^+ & 1 & 1 \\
\delta_{\alpha_1}^- & 1 & -1 \\
\delta_{\alpha_2}^- & -1 & 1
\end{array}$\end{center}
\end{scriptsize} &
\mbox{$(SL(2)\times SL(2))$ / $B^-$}, \begin{scriptsize}\mbox{$B^-$ Borel} subgroup of $SL(2)$ embedded diagonally\end{scriptsize} \\
\hline
\hline
$D1(i)$ \begin{center}
\begin{picture}(8700,2700)\put(0,300){\usebox{\aone}}\put(300,1200){\usebox{\segm}}\put(1800,0){\usebox{\shortdn}}
\put(600,1600){\put(50,0){\usebox{\toe}}}\end{picture}
\end{center} &
\begin{scriptsize}
$\{\alpha_1,\ 2\alpha_2+\ldots+2\alpha_{n-2}+\alpha_{n-1}+\alpha_n\}$
\vspace{1ex}\begin{center}$\begin{array}{r|cc}
\hline
\delta_{\alpha_1}^+ & 1 & -1 \\
\delta_{\alpha_1}^- & 1 & -1 \\
\delta_{\alpha_2} & -1 & 2
\end{array}$\end{center}
\end{scriptsize} &
\mbox{$SO(2n)$ /} \mbox{$(GL(1)\times SO(2n-2))$} \\
\hline
$D1(ii)$ \begin{center}
\begin{picture}(8700,2700)\put(0,300){\usebox{\aprime}}\put(300,1200){\usebox{\segm}}\put(1800,0){\usebox{\shortdn}}\end{picture}
\end{center} &
\begin{scriptsize}
$\{2\alpha_1,\ 2\alpha_2+\ldots+2\alpha_{n-2}+\alpha_{n-1}+\alpha_n\}$
\vspace{1ex}\begin{center}$\begin{array}{r|cc}
\hline
\delta_{\alpha_1}^\prime & 2 & -1 \\
\delta_{\alpha_2} & -2 & 2
\end{array}$\end{center}
\end{scriptsize} &
$SO(2n)$ / $(N_{SL(2)}(GL(1))\times SO(2n-2))$ \\
\hline
$D3$ \begin{center}
\begin{picture}(5400,3300)\put(0,1200){\usebox{\dthree}}\put(3600,300){\usebox{\bifurc}}\multiput(4800,300)(0,2400){2}{\circle{600}}\thicklines\multiput(4800,300)(0,2000){2}{\line(0,1){400}}\multiput(4800,700)(-200,200){4}{\line(-1,0){200}}\multiput(4600,700)(-200,200){4}{\line(0,1){200}}\multiput(4800,2300)(-200,-200){4}{\line(-1,0){200}}\multiput(4600,2300)(-200,-200){4}{\line(0,-1){200}}\end{picture}
\end{center} &
\begin{scriptsize}
$\{\alpha_1+2\alpha_2+\alpha_3,\ \alpha_3+\alpha_4+\alpha_5\}$
\vspace{1ex}\begin{center}$\begin{array}{r|cc}
\hline
\delta_{\alpha_2} & 2 & -1 \\
\delta_{\alpha_4} & -1 & 1 \\
\delta_{\alpha_5} & -1 & 1
\end{array}$\end{center}
\end{scriptsize} &
$SO(10)$ / $GL(5)$ \\
\hline
$D4(i)$ \begin{center}
\begin{picture}(7650,3300)\put(0,1200){\usebox{\mediumam}}\put(5700,0){\usebox{\plusaoneaone}}\end{picture}
\end{center} &
\begin{scriptsize}
$\{\alpha_1+\ldots+\alpha_{n-2},\ \alpha_{n-1}+\alpha_n\}$
\vspace{1ex}\begin{center}$\begin{array}{r|cc}
\hline
\delta_{\alpha_1} & 1 & 0 \\
\delta_{\alpha_{n-2}} & 1 & -2 \\
\delta_{\alpha_{n-1}}=\delta_{\alpha_n} & -1 & 2
\end{array}$\end{center}
\end{scriptsize} &
$SO(2n)$ / $H$,
\begin{scriptsize}
\begin{center}$\left(\begin{array}{ccc}
\!c\,A_1\! & 0 & 0 \\
\cline{1-1}
\multicolumn{1}{|c|}{M} & \!\pi(A_2,A_2)\! & 0 \\
\cline{2-2}
\multicolumn{1}{|c}{\!\ast\!} & \multicolumn{1}{c|}{\!\ast\!} & \!\ast\! \\
\cline{1-2}
\end{array}\right)$\end{center}
for $A_1\in SL(n-2)$, $A_2\in SL(2)$, $c \in GL(1)$, $M=M^\prime+M^{\prime\prime}\in(\C^{n-2})^\ast\otimes\C^4\cong(\C^{n-2})^\ast\otimes(\mathfrak{sl}(2)\oplus\C)\cong((\C^{n-2})^\ast\otimes(\mathfrak{sl}(2))\oplus(\C^{n-2})^\ast)$, $M^{\prime\prime}=0$ 
\end{scriptsize} \\
\hline
$D4(ii)$ \begin{center}
\begin{picture}(8700,2700)\put(0,300){\usebox{\aone}}\put(300,1200){\usebox{\segm}}\put(1800,0){\usebox{\shortdn}}
\end{picture}
\end{center} &
\begin{scriptsize}
$\{\alpha_1,\ 2\alpha_2+\ldots+2\alpha_{n-2}+\alpha_{n-1}+\alpha_n\}$
\vspace{1ex}\begin{center}$\begin{array}{r|cc}
\hline
\delta_{\alpha_1}^+ & 1 & 0 \\
\delta_{\alpha_1}^- & 1 & -2 \\
\delta_{\alpha_2} & -1 & 2
\end{array}$\end{center}
\end{scriptsize} &
$SO(2n)$ / $H$,
\begin{scriptsize}
\begin{center}$\left(\begin{array}{ccc}
c & 0 & 0 \\
\cline{1-1}
\multicolumn{1}{|c|}{M} & j(A) & 0 \\
\cline{2-2}
\multicolumn{1}{|c}{\ast} & \multicolumn{1}{c|}{\ast} & \ast \\
\cline{1-2}
\end{array}\right)$\end{center}
for $A\in SO(2n-3)$, $c \in GL(1)$, $M=M^\prime+M^{\prime\prime}\in\C^{2n-2}\cong\C^{2n-3}\oplus\C$, $M^{\prime\prime}=0$ 
\end{scriptsize} \\
\hline
$D4(iii)$ \begin{center}
\begin{picture}(10500,2700)\put(0,900){\usebox{\shortam}}\put(3900,0){\usebox{\plusdm}}\end{picture}
\end{center} &
\begin{scriptsize}
$\{\alpha_1+\ldots+\alpha_m,\ 2\alpha_{m+1}+\ldots+2\alpha_{n-2}+\alpha_{n-1}+\alpha_n\}$
\vspace{1ex}\begin{center}$\begin{array}{r|cc}
\hline
\delta_{\alpha_1} & 1 & 0 \\
\delta_{\alpha_m} & 1 & -2 \\
\delta_{\alpha_{m+1}} & -1 & 2
\end{array}$\end{center}
\end{scriptsize} &
$SO(2n)$ / $H$,
\begin{scriptsize}
\begin{center}$\left(\begin{array}{ccc}
c\,A_1 & 0 & 0 \\
\cline{1-1}
\multicolumn{1}{|c|}{M} & j(A_2) & 0 \\
\cline{2-2}
\multicolumn{1}{|c}{\ast} & \multicolumn{1}{c|}{\ast} & \ast \\
\cline{1-2}
\end{array}\right)$\end{center}
for $A_1\in SL(m)$, $A_2\in SO(2n-2m-1)$, $c \in GL(1)$, $M=M^\prime+M^{\prime\prime}\in(\C^m)^\ast\otimes\C^{2n-2}\cong(\C^m)^\ast\otimes(\C^{2n-3}\oplus\C)\cong((\C^m)^\ast\otimes\C^{2n-3})\oplus(\C^m)^\ast$, $M^{\prime\prime}=0$ 
\end{scriptsize} \\
\hline
$D7$ \begin{center}
\begin{picture}(7200,3300)\multiput(300,1500)(3600,0){2}{\usebox{\segm}}\put(2100,1500){\usebox{\shortsusp}}\put(5700,300){\usebox{\bifurc}}\put(300,1500){\circle{600}}\multiput(6900,300)(0,2400){2}{\circle{600}}\multiput(300,1500)(25,25){13}{\circle*{70}}\put(600,1800){\multiput(0,0)(300,0){6}{\multiput(0,0)(25,-25){7}{\circle*{70}}}\multiput(150,-150)(300,0){6}{\multiput(0,0)(25,25){7}{\circle*{70}}}}\put(3600,1800){\multiput(0,0)(300,0){7}{\multiput(0,0)(25,-25){7}{\circle*{70}}}\multiput(150,-150)(300,0){6}{\multiput(0,0)(25,25){7}{\circle*{70}}}}\multiput(300,1500)(25,-25){13}{\circle*{70}}\put(600,1200){\multiput(0,0)(300,0){6}{\multiput(0,0)(25,25){7}{\circle*{70}}}\multiput(150,150)(300,0){6}{\multiput(0,0)(25,-25){7}{\circle*{70}}}}\put(3600,1200){\multiput(0,0)(300,0){7}{\multiput(0,0)(25,25){7}{\circle*{70}}}\multiput(150,150)(300,0){6}{\multiput(0,0)(25,-25){7}{\circle*{70}}}}\multiput(2700,1200)(0,600){2}{\multiput(0,0)(300,0){3}{\circle*{70}}}\thicklines\put(6900,2700){\line(-1,0){400}}\multiput(6500,2700)(-200,-200){5}{\line(0,-1){200}}\multiput(6500,2500)(-200,-200){4}{\line(-1,0){200}}\multiput(5700,1700)(-30,-10){5}{\line(-1,0){30}}\put(6900,300){\line(-1,0){400}}\multiput(6500,300)(-200,200){5}{\line(0,1){200}}\multiput(6500,500)(-200,200){4}{\line(-1,0){200}}\multiput(5700,1300)(-30,10){5}{\line(-1,0){30}}\end{picture}
\end{center} &
\begin{scriptsize}
$\{\alpha_1+\ldots+\alpha_{n-2}+\alpha_{n-1},\ \alpha_1+\ldots+\alpha_{n-2}+\alpha_n\}$
\vspace{1ex}\begin{center}$\begin{array}{r|cc}
\hline
\delta_{\alpha_1} & 1 & 1 \\
\delta_{\alpha_{n-1}}^+ & 1 & -1 \\
\delta_{\alpha_n}^- & -1 & 1
\end{array}$\end{center}
\end{scriptsize} &
$SO(2n)$ / $H$,
\begin{scriptsize}
\begin{center}$\left(\begin{array}{cccc}
c\,A & 0 & 0 & 0 \\
\cline{1-1}
\multicolumn{1}{|c|}{M} & 1 & 0 & 0 \\
\cline{2-2}
\multicolumn{1}{|c}{-M} & \multicolumn{1}{c|}{0} & 1 & 0 \\
\cline{3-3}
\multicolumn{1}{|c}{\ast} & \ast & \multicolumn{1}{c|}{\ast} & \ast \\
\cline{1-3}
\end{array}\right)$\end{center}
for $A\in SL(n-1)$, $c \in GL(1)$, $M\in (\C^{n-1})^\ast$
\end{scriptsize} \\
\hline
\end{longtable}
\end{small}

Now let $S$ be the set of simple roots associated to $B$ and let $X$ be a wonderful $G$-variety. For every $\alpha\in S$, let $P_{\{\alpha\}}$ be the standard parabolic subgroup associated to $\alpha $; let $\Delta_X(\alpha)$ denote the set of non-$P_{\{\alpha\}}$-stable colours. For brevity we say that $\alpha$ ``moves'' the colours in $\Delta_X(\alpha)$, and we should point out that a colour is always moved by some simple root.

\begin{lemma}[\cite{L97}]\label{lemma:grcel}For all $\alpha\in S$, $\Delta_X(\alpha)$ has at most two elements. The following four distinct cases can occur. \begin{enumerate} \item $\Delta_X(\alpha) = \emptyset$, this happens exactly when the open $B$ orbit is stable under $P_{\{\alpha\}}$. The set of all such $\alpha$ is denoted by $S_X^p$; we have obviously $P_{S_X^p}= P_X$. \item $\Delta_X(\alpha)$ has two elements, this happens exactly when $\alpha \in \Sigma_X$. The two colours in $\Delta_X(\alpha)$ are denoted by $D_\alpha^+$, $D_\alpha^-$, and we have: \[\langle \rho(D_\alpha^+), \gamma\rangle + \langle \rho(D_\alpha^-), \gamma\rangle = \langle \alpha^\vee,\gamma\rangle\] for every $\gamma \in \Sigma_X$. We denote by ${\bf A}_X$ the union of $\Delta_X(\alpha)$ for $\alpha\in S\cap\Sigma_X$. \item $\Delta_X(\alpha)$ has one element and $2\alpha\in\Sigma_X$. The colour in $\Delta_X(\alpha)$ is denoted by $D^\prime_\alpha$, and: \[\langle \rho(D^\prime_\alpha), \gamma\rangle = \frac{1}{2} \langle \alpha^\vee, \gamma\rangle\] for every $\gamma \in \Sigma_X$. \item The remaining case, i.e.\ $\Delta_X(\alpha)$ has one element, but $2\alpha\notin\Sigma_X$. In this case, the colour in $\Delta_X(\alpha)$ is denoted by $D_\alpha$, and: \[\langle \rho(D_\alpha), \gamma\rangle = \langle \alpha^\vee, \gamma\rangle\] for every $\gamma \in \Sigma_X$. \end{enumerate} \end{lemma}


\begin{lemma}[\cite{L01}]\label{lemma:colours} For all $\alpha,\beta\in S$, the condition $\Delta_X(\alpha)\cap\Delta_X(\beta)\neq\emptyset$ occurs only in the following two cases: \begin{itemize} \item[-] if $\alpha,\beta\in S\cap\Sigma_X$ then it can happen that the cardinality of $\Delta_X(\alpha)\cup\Delta_X(\beta)$ is equal to $3$, \item[-] if $\alpha$ and $\beta$ are orthogonal and $\alpha+\beta$ (or $\frac{1}{2}(\alpha+\beta)$) belongs to $\Sigma_X$, then $D_\alpha=D_\beta$.\end{itemize} \end{lemma}

The relations appearing in these two lemmas have been found using some analysis of the cases in rank $1$ and $2$, and they will appear in the next section as axioms for an abstract combinatorial structure called \textit{spherical system}. These objects take into account the set $S^p_X$ of simple roots moving no colour, the set $\Sigma_X$ of spherical roots and the subset of colours ${\bf A}_X$.

\subsection{Spherical systems}

\begin{definition} Let $(S^p,\Sigma,\mathbf A)$ be a triple such that $S^p \subset S$, $\Sigma \subset \Sigma(G)$ and $\mathbf A$ a finite set endowed with an application $\rho\colon \mathbf A \to \Xi ^\ast$, where $\Xi = \langle \Sigma \rangle \subset \Xi (T)$. For every $\alpha \in \Sigma \cap S$, let $\mathbf A (\alpha)$ denote the set $\{\delta \in \mathbf A \colon \delta(\alpha)=1 \}$. Such a triple is called a \textit{spherical system} for $G$ if: \begin{itemize} \item[(A1)] for every $\delta \in \mathbf A$ and $\gamma \in \Sigma$ we have $\langle \rho (\delta),\gamma \rangle \leq 1$, and if $\delta (\gamma)=1$ then $\gamma \in S\cap\Sigma$; \item[(A2)] for every $\alpha \in \Sigma \cap S$, $\mathbf A(\alpha)$ contains two elements and by denoting with $\delta_\alpha^+$ and $\delta_\alpha^-$ these elements, it holds $\langle \rho (\delta_\alpha^+),\gamma \rangle + \langle \rho (\delta_\alpha^-),\gamma \rangle = \langle \alpha^\vee , \gamma \rangle$, for every $\gamma \in \Sigma$; \item[(A3)] the set $\mathbf A$ is the union of $\mathbf A(\alpha)$ for all $\alpha\in\Sigma \cap S$; \item[($\Sigma 1$)] if $2\alpha \in \Sigma \cap 2S$ then $\frac{1}{2}\langle\alpha^\vee, \gamma \rangle$ is a nonpositive integer for every $\gamma \in \Sigma \setminus \{ 2\alpha \}$; \item[($\Sigma 2$)] if $\alpha, \beta \in S$ are orthogonal and $\alpha + \beta \in \Sigma$ (or $\frac{1}{2}(\alpha + \beta) \in \Sigma$) then $\langle \alpha ^\vee , \gamma \rangle = \langle \beta ^\vee , \gamma \rangle$ for every $\gamma \in \Sigma$; \item[(S)] for every $\gamma \in \Sigma$, there exists a wonderful $G$-variety $X$ of rank $1$ with $\gamma$ as spherical
root and $S^p =S^p_X$. \end{itemize} \end{definition}

The cardinality (rank) of $\Sigma$ is the \textit{rank} of the spherical system.

Using the finite list of rank one prime wonderful varieties, in case the group $G$ is of type $\mathsf A$~$\mathsf D$ the last axiom can be rewritten as follows: \begin{itemize} \item[(S)] (type $\mathsf A$~$\mathsf D$) for every $\gamma \in \Sigma$, \[\{\alpha \in S \colon \langle \alpha ^\vee , \gamma \rangle=0 \} \cap \supp(\gamma) \subset S^p \subset \{\alpha \in S \colon \langle \alpha ^\vee , \gamma \rangle=0 \}, \] \end{itemize} where the set $\supp(\gamma)$\index{$\supp(\gamma)$, support of $\gamma$} is the \textit{support} of $\gamma$, namely $\{\alpha \in S \colon \langle \alpha ^\ast , \gamma \rangle\neq 0 \}$ (here $\{\alpha^\ast\}$ denotes the dual basis of $S$).

\begin{remark}For every spherical system the spherical roots are linearly independent and $\Sigma$ is a basis of $\Xi$.\end{remark}

On a wonderful variety $X$ the information contained in $S_X^p$, $\Sigma_X$ and ${\mathbf A}_X$ is enough to ``recover'' the full set of colours $\Delta_X$ at least as an abstract set endowed with an application $\rho_X\colon\Delta_X \to \Xi ^\ast$. This is easily done using lemmas \ref{lemma:grcel} and \ref{lemma:colours}; let us see this procedure applied to an abstract spherical system $(S^p,\Sigma,\mathbf A)$. Its \textit{set of colours} $\Delta$ is defined as follows:
\[
\Delta = \Delta^a \sqcup \Delta ^{a^\prime} \sqcup \Delta ^b
\]
where \begin{itemize} \item[] $\Delta^a = \mathbf A$, \item[] $\Delta ^{a^\prime} = S\cap\frac 12\Sigma$, \item[] $\Delta ^b = S^b/_\sim$, where $S^b=\left(S\setminus(\Sigma\cup\frac12\Sigma\cup S^p)\right)$ and $\alpha \sim \beta$ if and only if $\alpha \perp \beta$ and $\alpha + \beta \in \Sigma$ (or $\frac{1}{2}(\alpha+\beta)\in \Sigma$). \end{itemize} The set of elements of $\Delta$ ``corresponding'' to $\alpha \in S$ is denoted by $\Delta (\alpha)$, and we say that these are the colours {\em moved} by $\alpha$. The element corresponding to $\alpha \in S\cap\frac12\Sigma$ in $\Delta^{a^\prime}$ is denoted by $\delta_\alpha^\prime$, and the element corresponding to $\alpha \in S^b$ in $\Delta^b$ is denoted by $\delta_\alpha$. The map $\rho\colon\mathbf A \to \Xi ^\ast$ can be extended to $\Delta$ as follows: \begin{itemize} \item[] $\langle \rho (\delta_\alpha^\prime), \gamma \rangle = \frac{1}{2} \langle \alpha ^\vee , \gamma \rangle \quad \forall \alpha \in S\cap\frac12\Sigma, \forall \gamma \in \Xi$; \item[] $\langle \rho (\delta_\alpha), \gamma \rangle = \langle \alpha ^\vee , \gamma \rangle \quad \forall \alpha \in S^b, \forall \gamma \in \Xi$. \end{itemize}

The definition of spherical system is such that the following lemmas hold (see \cite{L01}):

\begin{lemma}\label{L1} For every wonderful $G$-variety $X$ the triple $(S^p_X,\Sigma_X, \mathbf A _X)$ is a spherical system. \end{lemma}

\begin{lemma}\label{L2} The map $X\mapsto (S^p_X, \Sigma _X, \mathbf A _X)$ is a bijection between rank one (resp.\ rank two) wonderful varieties (up to $G$-isomorphism) and rank one (resp. rank two) spherical systems. \end{lemma}

In \cite{L01} it is proven that spherical systems classify wonderful $G$-varieties for $G$ adjoint of type $\mathsf A$. Our main theorem is the generalization of this result to type $\mathsf A$~$\mathsf D$.

\begin{theorem}\label{thmAD} Let $G$ be a semisimple adjoint algebraic group of mixed type $\mathsf A$ and $\mathsf D$. Then the map associating to a wonderful $G$-variety $X$ the triple $(S^p_X, \Sigma _X, \mathbf A _X)$ is a bijection between wonderful $G$-varieties (up to $G$-isomorphism) and spherical systems for $G$. \end{theorem}

\begin{remark} Obviously, two $G$-isomorphic wonderful $G$-varieties have the same spherical system. If two wonderful $G$-varieties are isomorphic, namely $G$-isomorphic up to an outer automorphism of $G$, their spherical systems are equal up to a permutation of the set $S$ of simple roots.\end{remark}

A useful way to represent a spherical system for $G$ is provided by the Luna diagram constructed on the Dynkin diagram of the root system of $G$. We focus on adjoint groups of type $\mathsf A$~$\mathsf D$, following \cite{L01} in type $\mathsf A$.

For every $\alpha \in S$, we draw one or two circles near the corresponding vertex to represent the colours moved by $\alpha$. Different circles corresponding to the same colour are joined by a line.

In the Luna diagram each spherical root in $\Sigma$ is denoted as in Table~\ref{T1}\footnote{For the type $a_m$, $m\geq 2$, we modify the convention of D.~Luna and we partially follow that of R.~Camus (\cite{Cam01}). Notice that the zigzag line is drawn only to denote this type of spherical root and not to join the two circles, since they correspond to two distinct colours.} on the vertices corresponding to $\supp(\gamma)$.

In order to represent the set $\mathbf A$ and the map $\rho\colon \mathbf A \to \Xi^\ast$, we arrange the diagram in such a way that the circle above a simple root $\alpha$ (corresponding to a colour denoted by $\delta^+_\alpha$) satisfies $\langle\rho(\delta^+_\alpha),\gamma\rangle\geq-1$ for all $\gamma\in\Sigma$. Moreover we draw one arrow, $>$ or $<$, near the circle above the simple root $\alpha$ pointing towards the spherical root $\gamma$ when $\gamma$ is nonorthogonal to $\alpha$ and $\langle\rho(\delta^+_\alpha),\gamma\rangle=-1$. These data are sufficient to recover the map $\rho$, using axioms A1 and A2 of spherical systems. In some cases, we report the arrows analogously for the circles below the simple roots in order to avoid confusion.


\section{Reduction to the primitive cases}

\subsection{Luna's lemmas}\label{SS31}

\subsubsection*{Localization}

(see also \cite{L97}) Let $X$ be a wonderful $G$-variety, let $\Sigma^\prime$ be a subset of the set $\Sigma_X$ of the spherical roots. The variety \[ X^{\Sigma^\prime}=\bigcap_{\gamma\in\,\Sigma\setminus\Sigma^\prime} D^\gamma\] is called the \textit{localization} of $X$ in $\Sigma^\prime$. It is a wonderful $G$-variety with \begin{itemize} \item[] $S^p_{X^{\Sigma^\prime}}=S^p_X$, \item[] $\Sigma_{X^{\Sigma^\prime}}=\Sigma^\prime$, \item[] $\mathbf A _{X^{\Sigma^\prime}}=\mathbf A_X(\Sigma^\prime\cap S)$, namely the union of $\mathbf A_X(\alpha)$ for all $\alpha\in\,\Sigma\cap S$, and $\rho_{X^{\Sigma^\prime}}$ is the restriction of $\rho_X$ to $\mathbf A _{X^{\Sigma^\prime}}$ followed by the projection $\Xi^\ast\to\langle\Sigma^\prime\rangle^\ast$.\end{itemize}

Let $S^\prime$ be a subset of the set $S$ of simple roots. Let $P_{S^\prime}$ be the standard parabolic subgroup containing $B$ and associated to $S^\prime$ and $L^\prime$ a Levi subgroup. Let $C^\prime$ denote the radical of $L^\prime$ and let $X^{C^\prime}\subset X$ denote the fixed-point set of $C^\prime$. The connected component of $X^{C^\prime}$ containing the unique point $z$ fixed by $B^-$ is the \textit{localization} $X^{S^\prime}$ of $X$ in $S^\prime$. The variety $X^{S^\prime}$ is wonderful under the action of $L^\prime$, with \begin{itemize} \item[] $S^p_{X^{S^\prime}}=S^p_X\cap S^\prime$, \item[] $\Sigma_{X^{S^\prime}}=\{\gamma\in\Sigma_X\colon \supp(\gamma)\subset S^\prime\}$, \item[] $\mathbf A_{X^{S^\prime}}=\mathbf A(\Sigma\cap S^\prime)$ and $\rho_{X^{S^\prime}}$ is the restriction of $\rho_X$ followed by the projection $\Xi^\ast\to\langle\Sigma_{X^{S^\prime}}\rangle^\ast$. \end{itemize}

In general, let $(S^p,\Sigma,\mathbf A)$ be a spherical system for $G$. Let $S^\prime$ be a subset of $S$ and $\Sigma^\prime$ a subset of $\Sigma$ such that $\supp(\Sigma^\prime)\subset S^\prime$. The \textit{localization} of $(S^p,\Sigma,\mathbf A)$ in $(S^\prime,\Sigma^\prime)$ is the spherical system $(S^{\prime p},\Sigma^\prime,\mathbf A^\prime)$ where $S^{\prime p}=S^p\cap S^\prime$, $\mathbf A^\prime=\mathbf A(S\cap\Sigma^\prime)$ and $\rho^\prime$ is the restriction of $\rho$ followed by the projection $\Xi^\ast\to\langle\Sigma^\prime\rangle^\ast$.

The variety $(X^{\Sigma^\prime})^{S^\prime}=(X^{S^\prime})^{\Sigma^\prime}$ is a wonderful $L^\prime$-variety whose spherical system is the localization of $(S^p_X,\Sigma_X,\mathbf A_X)$ in $(S^\prime,\Sigma^\prime)$.

\subsubsection*{Quotient}

Let $(S^p,\Sigma,\mathbf A)$ be a spherical system, let $\Delta$ be its set of colours. Let $V$ be the $\mathbb Q$-vector space $\mathbb Q \otimes \Xi^\ast$, and let $\mathcal V$ be the cone $\{v\in V\colon \langle v,\gamma\rangle\leq0\ \forall \gamma\in\Sigma\}$. A subset $\Delta^\prime\subset\Delta$ is said to be \textit{distinguished} if the intersection of the interior of the convex cone generated by the set $\rho(\Delta^\prime)$ with the cone $-\mathcal V$ is nonempty, namely if there exists a linear combination of the elements of $\rho(\Delta^\prime)$ with positive coefficients: \[ \phi=\sum _{\delta\in\Delta^\prime} c_{\delta}\, \rho(\delta),\quad c_{\delta}>0\ \forall \delta\in\Delta^\prime,\] such that $\langle\phi,\gamma\rangle\geq 0$ for all $\gamma\in\Sigma$. Let $\Sigma(\Delta^\prime)$ denote the maximal subset of spherical roots such that there exists a linear combination $\phi$ as above such that $\langle\phi,\gamma\rangle>0$ for all $\gamma\in\Sigma(\Delta^\prime)$.
The subset $\Delta^\prime\subset\Delta$ is distinguished if and only if there exists a vector subspace $V^\prime$ of $V$ such that the convex cone generated by $\rho(\Delta^\prime)$ and $V^\prime\cap\mathcal V$ is equal to $V^\prime$, and $V^\prime\cap\mathcal V$ is a face of $\mathcal V$. If such a subspace $V^\prime$ exists, it is unique and it is denoted by $V(\Delta^\prime)$. We have \[ V(\Delta^\prime)=\langle\rho(\Delta^\prime)\cup\{\gamma^\ast\colon
\gamma\in\Sigma(\Delta^\prime)\}\rangle.\] where $\{\gamma^\ast\}$ is the dual
basis of $\Sigma$.

Let $(S^p,\Sigma,\mathbf A)$ be a spherical system and $\Delta^\prime$ a distinguished subset of $\Delta$. Let $\Xi/\Delta^\prime$ be the subgroup of $\Xi$ defined by: \begin{eqnarray*} \frac{\Xi}{\Delta^\prime} & = & \{\xi\colon \langle v,\xi\rangle=0\ \forall v\in V(\Delta^\prime)\}\\ & = & \{\xi\colon \langle\rho(\delta),\xi\rangle=0\ \forall\delta\in\Delta^\prime\ \mathrm{and}\ \langle\gamma^\ast,\xi\rangle=0\ \forall\gamma\in\Sigma(\Delta^\prime)\}.\end{eqnarray*} The triple $(S^p/\Delta^\prime,\Sigma/\Delta^\prime, \mathbf A/\Delta^\prime)$ defined as follows is called the \textit{quotient} triple of $(S^p,\Sigma,\mathbf A)$ by $\Delta^\prime$ and it is also denoted by $(S^p,\Sigma,\mathbf A)/\Delta^\prime$. \begin{itemize} \item[] The set $S^p/\Delta^\prime$ is $\{\alpha\in S^p\colon \Delta(\alpha)\subset\Delta^\prime\}\subset S^p.$ \item[] The set $\Sigma/\Delta^\prime$\index{$\Sigma/\Delta^\prime$} is the set of the indecomposable elements of the semigroup \[\{\sum_{\gamma\in\Sigma}c_\gamma \gamma\,\in\Xi/\Delta^\prime\colon c_\gamma\geq0\ \forall\gamma\in\Sigma\}.\] \item[] The set $\mathbf A/\Delta^\prime\subset\mathbf A$ is the union of those $\mathbf A(\alpha)$ which satisfy $ \mathbf A(\alpha) \cap \Delta^\prime=\emptyset$, and the map $\rho/\Delta^\prime\colon\mathbf A/\Delta^\prime\to\left(\Xi/\Delta\prime\right)^\ast$ is the restriction of $\rho$ followed by the projection $\Xi^\ast\to\left(\Xi/\Delta^\prime\right)^\ast$.\end{itemize} In general, the set $\Sigma/\Delta^\prime$ is not always a subset of $\Sigma$. If it is a basis of $\Xi/\Delta^\prime$ then the distinguished subset $\Delta^\prime$ is said to have the \textit{$\ast$-property} of Luna. In this case the quotient triple is a spherical system.

\begin{remark}It is not clear whether there exist examples of distinguished subsets without $\ast$-property. All distinguished subsets considered in the following have the $\ast$-property. \end{remark}

Let $\Phi\colon X\to X^\prime$ be a dominant $G$-morphism between wonderful $G$-varieties. Set $\Delta_\Phi=\{D\in \Delta_X\colon\Phi(D)=X^\prime\}$.

\begin{proposition}[\cite{L01}]\label{morphisms} Let $X$ be a wonderful $G$-variety: the application $\Phi\mapsto\Delta_\Phi$ is a bijection between the set of dominant $G$-morphisms with connected fibers of $X$ onto another wonderful $G$-variety, and the set of distinguished subsets of  $\Delta_X$ having the $\ast$-property.

Moreover, for any such $\Phi\colon X\to X^\prime$ the spherical system of $X^\prime$ is the quotient triple $(S^p_X,\Sigma_X,\mathbf A_X)/\Delta_\Phi$.\end{proposition}

The distinguished subset $\Delta^\prime$ of $\Delta$ is \textit{smooth} if the dimension of the space $V(\Delta^\prime)$ is equal to the dimension of the face $V(\Delta^\prime)\cap\mathcal V$ of $\mathcal V$. Equivalently, if $V(\Delta^\prime)=\langle\gamma^\ast\colon \gamma\in\Sigma(\Delta^\prime)\rangle$ and, hence, $\Sigma/\Delta^\prime=\Sigma\setminus\Sigma(\Delta^\prime)\subset\Sigma$. In particular, if $\Delta^\prime$ is smooth then it has the $\ast$-property.

\begin{proposition}[\cite{L01}] Let $\Phi\colon X\to X^\prime$ be a dominant $G$-morphism with connected fibers between wonderful $G$-varieties. Then $\Phi$ is smooth if and only if $\Delta_\Phi$ is smooth.\end{proposition}

The distinguished subset $\Delta^\prime$ of $\Delta$ is \textit{parabolic} if $V(\Delta^\prime)=V$, or equivalently if $\Sigma(\Delta^\prime)=\Sigma$, namely $\Sigma/\Delta^\prime=\emptyset$. In particular, in this case $\Delta^\prime$ is smooth.

\begin{proposition}[\cite{L01}] Let $X$ be a wonderful $G$-variety and let $S^\prime$ be a subset of $S$. The map $\Phi\mapsto\Delta_\Phi$ is a bijection between the set of the $G$-morphisms $\Phi\colon X\to G/P_{S^\prime}$ and the set of the parabolic distinguished subsets of $\Delta_X$; the bijection is such that $S^\prime=S^p/\Delta^\prime$.\end{proposition}

\subsubsection*{Parabolic induction}

Let $P$ be a parabolic subgroup of $G$ and let $L\subset P$ be a Levi subgroup. Let $Y$ be a wonderful $L$-variety. The \textit{parabolic induction} obtained from $Y$ by $P$ is the wonderful $G$-variety $G\times_P Y$, where $Y$ is considered as a $P$-variety where the radical $P^r$ of $P$ acts trivially and $G\times_P Y$ is the algebraic quotient of $G\times Y$ by the following action of $P$: $p.(g,y)=(g\,p^{-1},p.y)$, for all $p\in P$, $g\in G$, $y\in Y$.

On the other hand, a wonderful variety $X$ is obtained by parabolic induction by $P$ if and only if $P^r\subset H\subset P$, where $G/H$ is the open $G$-orbit of $X$.

\begin{proposition}[\cite{L01}]\label{P4} Let $X$ be a wonderful $G$-variety and let $S^\prime$ be a subset of $S$. Let $P^-_{S^\prime}$ denote the parabolic subset containing $B^-$ associated to $S^\prime$. Then $X$ is obtained by parabolic induction by $P^-_{S^\prime}$ if and only if $(\supp(\Sigma)\cup S^p)\subset S^\prime$.

In this case the $G$-morphism $\Phi\colon X\to G/P^-_{S^\prime}$ is unique and $\Delta_\Phi$ is equal to the parabolic distinguished subset $\Delta(S^\prime)$. \end{proposition}

A spherical system is \textit{cuspidal} if $\supp(\Sigma)=S$. In particular, a cuspidal spherical system cannot be obtained by parabolic induction by any proper parabolic subgroup.

\subsubsection*{Direct and fiber product}

A spherical system $(S^p,\Sigma,\mathbf A)$ is said to be \textit{reducible} if there exists a partition of $S$ into two subsets $S_1,S_2$ such that: \begin{itemize} \item[] $S_1 \perp S_2$, \item[] $\forall \gamma\in\Sigma$ $\supp(\gamma)\subset S_1$ or $\supp(\gamma)\subset S_2$, \item[] $\forall \delta\in \mathbf A(\Sigma\cap S_1)$, $\rho(\delta)$ is null on $S_2$, and vice versa. \end{itemize} In this case $(S^p,\Sigma,\mathbf A)$ is the \textit{direct product} of the spherical systems ($S^p\cap S_1$, $\Sigma_1$, $\mathbf A(\Sigma \cap S_1)$) and $(S^p\cap S_2, \Sigma_2, \mathbf A(\Sigma \cap S_2))$, where $\Sigma_i=\{\gamma\in\Sigma\colon \supp(\gamma)\subset S_i\}$. For this notion of reducibility, the obvious statement of ``unique factorization'' is true.

Let $(S^p,\Sigma, \mathbf A)$ be a spherical system. Let $\Delta_1$ and $\Delta_2$ be two distinguished subsets of $\Delta$, obviously $\Delta_3=\Delta_1\cup \Delta_2$ is distinguished. The subsets $\Delta_1$ and $\Delta_2$ \textit{decompose into fiber product} the spherical system $(S^p,\Sigma,\mathbf A)$ if: \begin{enumerate} \renewcommand{\theenumi}{\roman{enumi}} \item $\Delta_1\neq\emptyset$, $\Delta_2\neq\emptyset$ and $\Delta_1 \cap \Delta_2 = \emptyset$, \item $\Delta_1$, $\Delta_2$ and $\Delta_3$ have the $\ast$-property, \item $(\Sigma\setminus(\Sigma/\Delta_1)) \cap (\Sigma\setminus(\Sigma/\Delta_2))=\emptyset$, \item $((S^p/\Delta_1)\setminus S^p)\perp ((S^p/\Delta_2)\setminus S^p)$, \item $\Delta_1$ or $\Delta_2$ is smooth. \end{enumerate} Notice that $\Sigma\setminus(\Sigma/\Delta^\prime)=\Sigma(\Delta^\prime)$ when $\Delta^\prime$ is smooth.

\begin{proposition} [\cite{L01}]\label{P5} Let $(S^p,\Sigma,\mathbf A)$ be a spherical system, let $\Delta_1$ and $\Delta_2$ be two distinguished subsets that decompose $(S^p,\Sigma,\mathbf A)$. Let us suppose that, for $i=1,2,3$, there exists a wonderful $G$-variety $X_i$ unique up to $G$-isomorphism with spherical system $(S^p,\Sigma,\mathbf A)/\Delta_i$. Then there exists two $G$-morphisms $\Phi_1\colon X_1\to X_3$, $\Phi\colon X_2\to X_3$ such that the fiber product $X_1\times_{X_3}X_2$ is a wonderful $G$-variety with spherical system $(S^p,\Sigma,\mathbf A)$ and every wonderful $G$-variety with this spherical system is $G$-isomorphic to $X_1\times_{X_3}X_2$.\end{proposition}

Notice that the above operation of factorization is a particular case of this decomposition into fiber product. A reducible spherical system, $S=S_1\sqcup S_2$, is decomposed by $\Delta(S_1)$ and $\Delta(S_2)$, and it corresponds to a direct product of two wonderful varieties.

\subsubsection*{Projective fibration}

Let $X$ and $Y$ be two wonderful $G$-varieties. A $G$-morphism $\Phi\colon X\to Y$ is called a \textit{projective fibration} if it is smooth, all its fibers are isomorphic to a projective space of fixed dimension $n$ and $rank(X)=rank(Y)+n$. We can look at the spherical systems in order to characterize these morphisms $\Phi$ in terms of conditions on $\Delta_\Phi$.

Let $(S^p,\Sigma,\mathbf A)$ be a spherical system. An element $\delta$ of $\mathbf A$ is \textit{projective} if $\langle\rho(\delta),\gamma\rangle\geq0$ for all $\gamma\in\Sigma$. Notice that this condition is equivalent to $\langle\rho(\delta),\gamma\rangle\in\{0,1\}$, thanks to axiom A1.

Let $\delta\in\mathbf A$ be a projective element, let $S_\delta$ denote the set $\{\alpha\in\Sigma\cap S\colon \langle\rho(\delta),\alpha\rangle=1\}$. The subset $\{\delta\}$ is distinguished and smooth, the quotient triple $(S^p,\Sigma,\mathbf A)/\{\delta\}$ is a spherical system where: \begin{itemize} \item[] $S^p/\{\delta\}=S^p$, \item[] $\Sigma/\{\delta\}=\Sigma\setminus S_\delta$, \item[] $\mathbf A/\{\delta\}=\mathbf A(S\setminus S_\delta)$, and $\rho/\{\delta\}$ is the restriction of $\rho$ followed by the projection on $(\Xi/\{\delta\})^\ast$. \end{itemize}

\begin{proposition} \label{P6}Let $G$ be of type $\mathsf A$~$\mathsf D$. Let $(S^p,\Sigma,\mathbf A)$ be a spherical system with a projective element $\delta\in\mathbf A$. Let us suppose that there exists a wonderful $G$-variety $Y$ unique up to $G$-isomorphism with spherical system $(S^p,\Sigma,\mathbf A)/\{\delta\}$. Then there exists a wonderful $G$-variety $X$ unique up to $G$-isomorphism with spherical system $(S^p,\Sigma,\mathbf A)$ and if $\Phi_\delta\colon X\to X_\delta$ denotes the $G$-morphism such that $\Delta_{\Phi_\delta}=\{\delta\}$ we have that $X_\delta$ is $G$-isomorphic to $Y$.\end{proposition}

In \cite{L01} the above proposition is stated only for the type $\mathsf A$, but its proof holds also in the case of type $\mathsf A$~$\mathsf D$. Indeed, it is enough to notice that the key lemma (3.6.2, \textit{loc.\ cit.}) is true whenever no spherical root $\alpha$ being also a simple root lies in the support of any other spherical root $\gamma\in\Sigma$. This occurs in type $\mathsf A$ as well as it occurs in type $\mathsf A$~$\mathsf D$, as can be seen from the list of rank two spherical systems.

\subsection{Proof of Theorem~\ref{thmAD}: Reduction}
\label{Reduction}

In order to prove that for all spherical systems $(S^p,\Sigma,\mathbf A)$ for $G$ adjoint of type $\mathsf A$~$\mathsf D$ there exists a wonderful $G$-variety unique up to a $G$-isomorphism, it is sufficient to prove that this holds for all spherical systems called 
\textit{primitive}, that are: \begin{itemize} \item[-] cuspidal, \item[-] without projective elements, \item[-] indecomposable, and in particular irreducible. \end{itemize} Indeed, we can proceed by induction on the rank of $S$ and on the rank of $\Sigma$ and we can use propositions \ref{P4}, \ref{P5} and \ref{P6} about parabolic inductions, fiber products and projective fibrations.

The list of primitive spherical systems is contained in Table~\ref{T3}. Such a list can be obtained by a combinatorial reduction, that is explained in the next two paragraphs: we report some technical definitions and outline the strategy, similar to that of \cite{L01}.

\textit{How to read Table~\ref{T3}}. For any primitive wonderful $G$-variety $X$ we report the label and the Luna diagram of the spherical system and the homogeneous space $G/H\cong X^\circ_G$. When necessary we describe the subgroup $H$, as in Table~\ref{T2}, by providing its connected center $C^\circ$, the commutator subgroup $(L,L)$ of a Levi subgroup $L$ and the Lie algebra of its unipotent radical $H^u$. In some cases, where the component $d(n_i)$ occurs, we use further special notations: if $n_i=2$ (with $d(2)$ we denote $aa(1,1)$) $\phi_i(A)=\pi(A,a)$ for $A\in SL(2)$, $\pi\colon SL(2)\times SL(2) \to SO(4)$ isogeny and $V_i\cong \mathfrak{sl}(2)$; if $n_i>2$ $\phi_i(A)=j(A)$ for $A\in SO(2n_i-1)$, $j\colon SO(2n_i-1)\to SO(2n_i)$ injection and $V_i\cong \C^{2n_i-1}$.


\begin{small}
\begin{longtable}{|p{12cm}|}
\caption{Primitive wonderful $G$-varieties with rank at least three for $G$ adjoint of type $\mathsf A$~$\mathsf D$, up to isomorphism.}
\label{T3} \\
\hline
Diagram, $H$, ($C^\circ$, $(L,L)$, $Lie(H^u)$) \makebox[200pt][r]{$G$}\\
\hline
\hline
\endfirsthead
\caption{Primitive wonderful $G$-varieties (continuation).} \\
\hline
Diagram, $H$, ($C^\circ$, $(L,L)$, $Lie(H^u)$) \makebox[200pt][r]{$G$}\\
\hline
\hline
\endhead
\hline
\endfoot

\atcpntpn $\mathbf{ao(n)}$.\makebox[308pt][r]{$SL(n+1)$}\begin{center}
\begin{picture}(7800,900)\multiput(300,900)(5400,0){2}{\usebox{\segm}}\put(2100,900){\usebox{\susp}}\multiput(0,0)(5400,0){2}{\multiput(0,0)(1800,0){2}{\usebox{\aprime}}}\end{picture}
\end{center}
$SO(n+1)\cdot C_G$ \\
\hline

\atcpntpn $\mathbf{ac(n)}$, $n\geq3$ odd. \makebox[260pt][r]{$SL(n+1)$} \begin{center}
\begin{picture}(10800,600)\put(0,0){\usebox{\dthree}}\put(3600,300){\usebox{\susp}}\put(7200,0){\usebox{\dthree}}\end{picture}
\end{center} $Sp(n+1)\cdot C_G$\\
\hline

\atcpntpn $\mathbf{aa(p+q+p)}$. \makebox[270pt][r]{$SL(2p+q+1)$}\begin{center}
\begin{picture}(15000,1500)\multiput(300,1200)(10800,0){2}{\usebox{\susp}}\multiput(3900,1200)(5400,0){2}{\usebox{\segm}}\put(5400,900){\usebox{\shortam}}\multiput(300,1200)(14400,0){2}{\circle{600}}\multiput(3900,1200)(7200,0){2}{\circle{600}}
\multiput(3900,900)(7200,0){2}{\line(0,-1){450}}\put(3900,450){\line(1,0){7200}}\multiput(300,900)(14400,0){2}{\line(0,-1){900}}\put(300,0){\line(1,0){14400}}
\end{picture}
\end{center} $S(GL(p+q)\times GL(p+1))$\\
\hline
\end{longtable}
\end{small}

\begin{remark}Some primitive spherical systems reported in Table~\ref{T3} are associated to the wonderful embedding of symmetric varieties (see \ref{reductivecases}). In Table~\ref{T3.2} we report the corresponding Satake diagrams.\end{remark}

\begin{small}
\begin{longtable}{|p{4.7cm}|p{6.3cm}|}
\caption{Symmetric cases.}
\label{T3.2}\\
\hline
Spherical system, $G$, $H$ & \multicolumn{1}{|c|}{Satake diagram}\\
\hline
\hline
\endfirsthead
\mbox{$\mathbf{ao(n)}$, $G=SL(n+1)$,} \mbox{$H=SO(n+1)\cdot C_G$} & \begin{center}\includegraphics[scale=0.3]{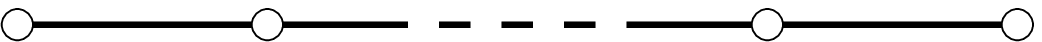}\end{center}\\
\hline
\mbox{$\mathbf{ac(n)}$, $G=SL(n+1)$,} \mbox{$H=Sp(n+1)\cdot C_G$} & \begin{center}\includegraphics[scale=0.3]{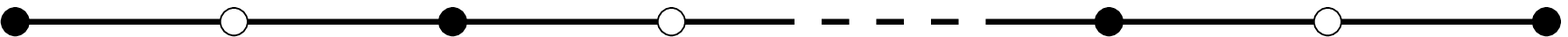}\end{center}\\
\hline
\mbox{$\mathbf{aa(p+q+p)}$, $q\geq2$,} \mbox{$G=SL(2p+q+1)$,} \mbox{$H=S(GL(p+q)\times GL(p+1))$} & \begin{center}\includegraphics[scale=0.3]{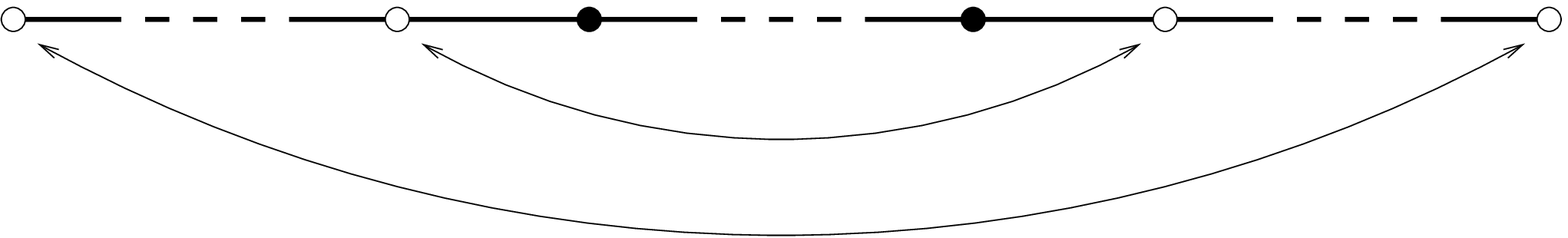}\end{center}\\
\hline
\mbox{$\mathbf{aa(p,p)}$,} \mbox{$G=SL(p+1)\times SL(p+1)$,} \mbox{$H=SL(p+1)\cdot C_G$} & \begin{center}\includegraphics[scale=0.3]{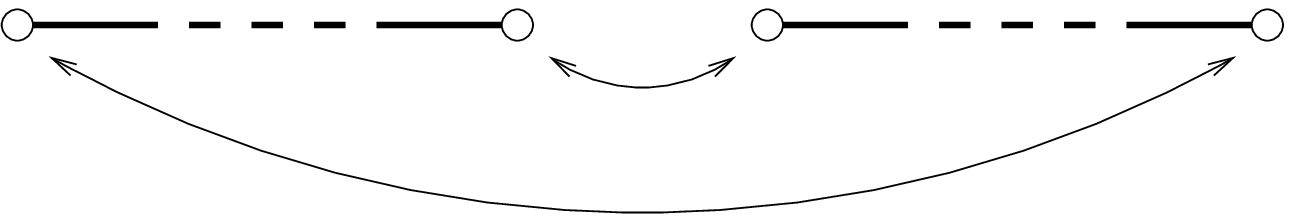}\end{center}\\
\hline
\mbox{$\mathbf{aa(p+1+p)}$, $G=SL(2p+2)$,} \mbox{$H=N_G(SL(p+1)\times SL(p+1))$} & \begin{center}\includegraphics[scale=0.3]{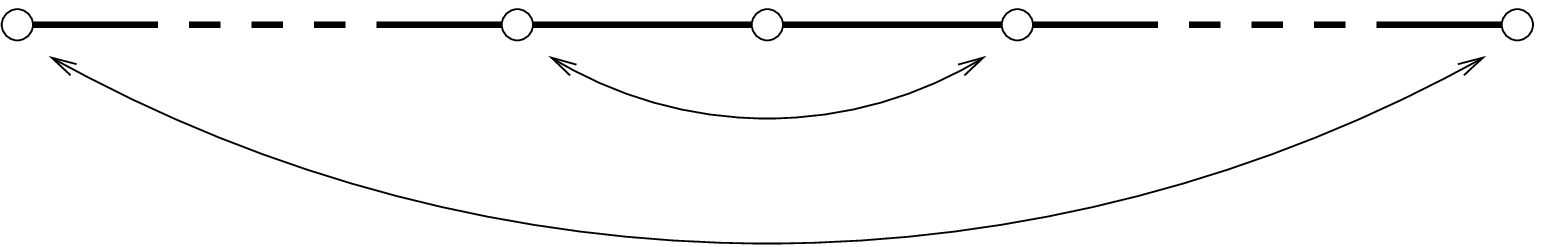}\end{center}\\
\hline
\mbox{$\mathbf{do(p+q)}$, $G=SO(2p+2q)$,} \mbox{$H=S(O(p+1)\times O(2q+p-1))$} & \begin{center}\includegraphics[scale=0.3]{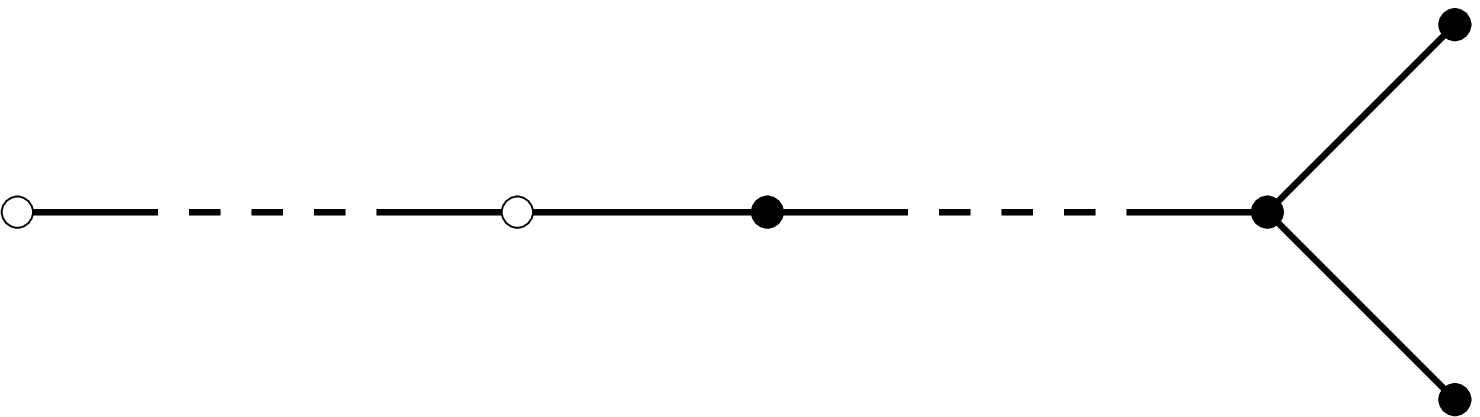}\end{center}\\
\hline
\mbox{$\mathbf{do(p+2)}$, $G=SO(2n)$,} \mbox{$H=S(O(n-1)\times O(n+1))$} & \begin{center}\includegraphics[scale=0.3]{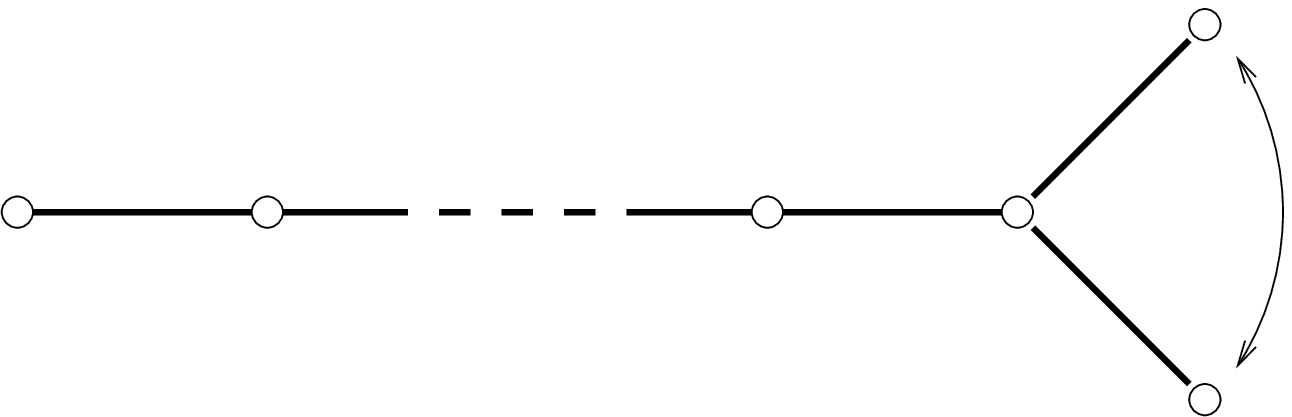}\end{center}\\
\hline
\mbox{$\mathbf{do(n)}$, $G=SO(2n)$,} \mbox{$H=N_G(SO(n)\times SO(n))$} & \begin{center}\includegraphics[scale=0.3]{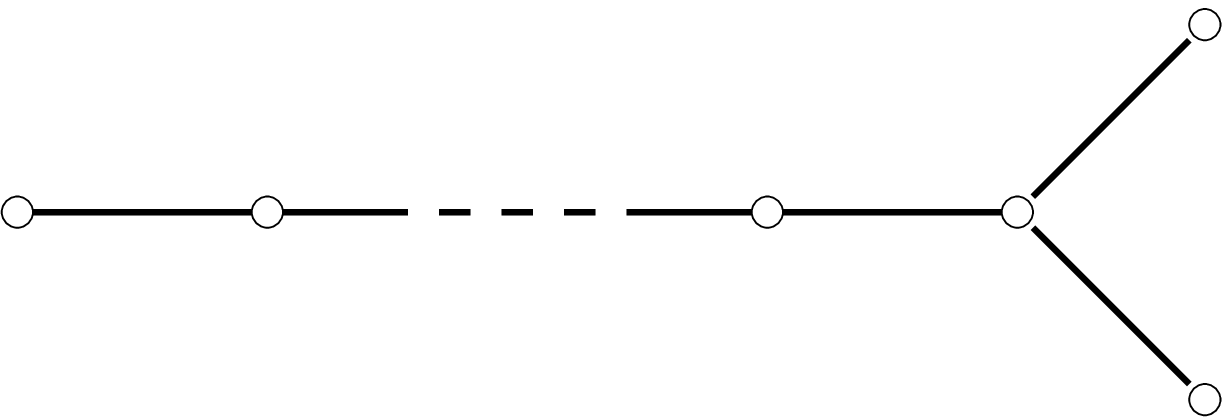}\end{center}\\
\hline
\mbox{$\mathbf{dc(n)}$, $n$ odd, $G=SO(2n)$,} \mbox{$H=GL(n)$} & \begin{center}\includegraphics[scale=0.3]{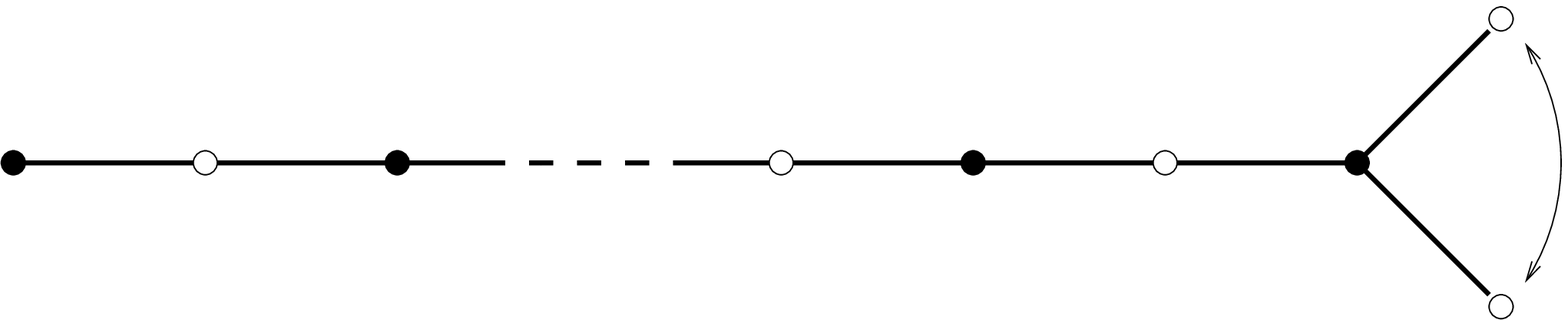}\end{center}\\
\hline
\mbox{$\mathbf{dc^\prime(n)}$, $n$ even, $G=SO(2n)$,} \mbox{$H=N_G(GL(n)$} & \begin{center}\includegraphics[scale=0.3]{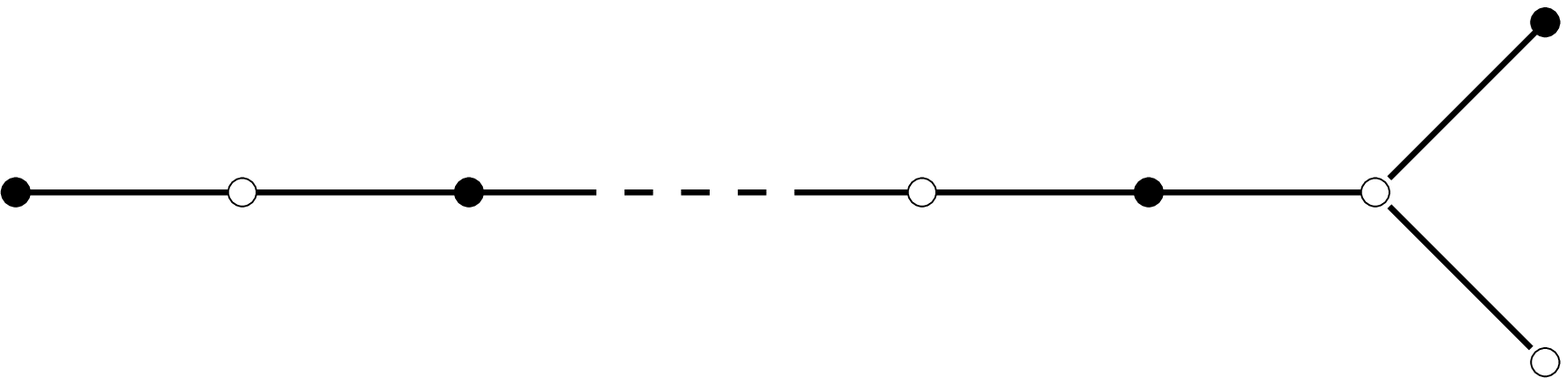}\end{center}\\
\hline
\mbox{$\mathbf{dd(p,p)}$, $G=SO(2p)\times SO(2p)$,} \mbox{$H=SO(2p)\cdot C_G$} & \begin{center}\includegraphics[scale=0.3]{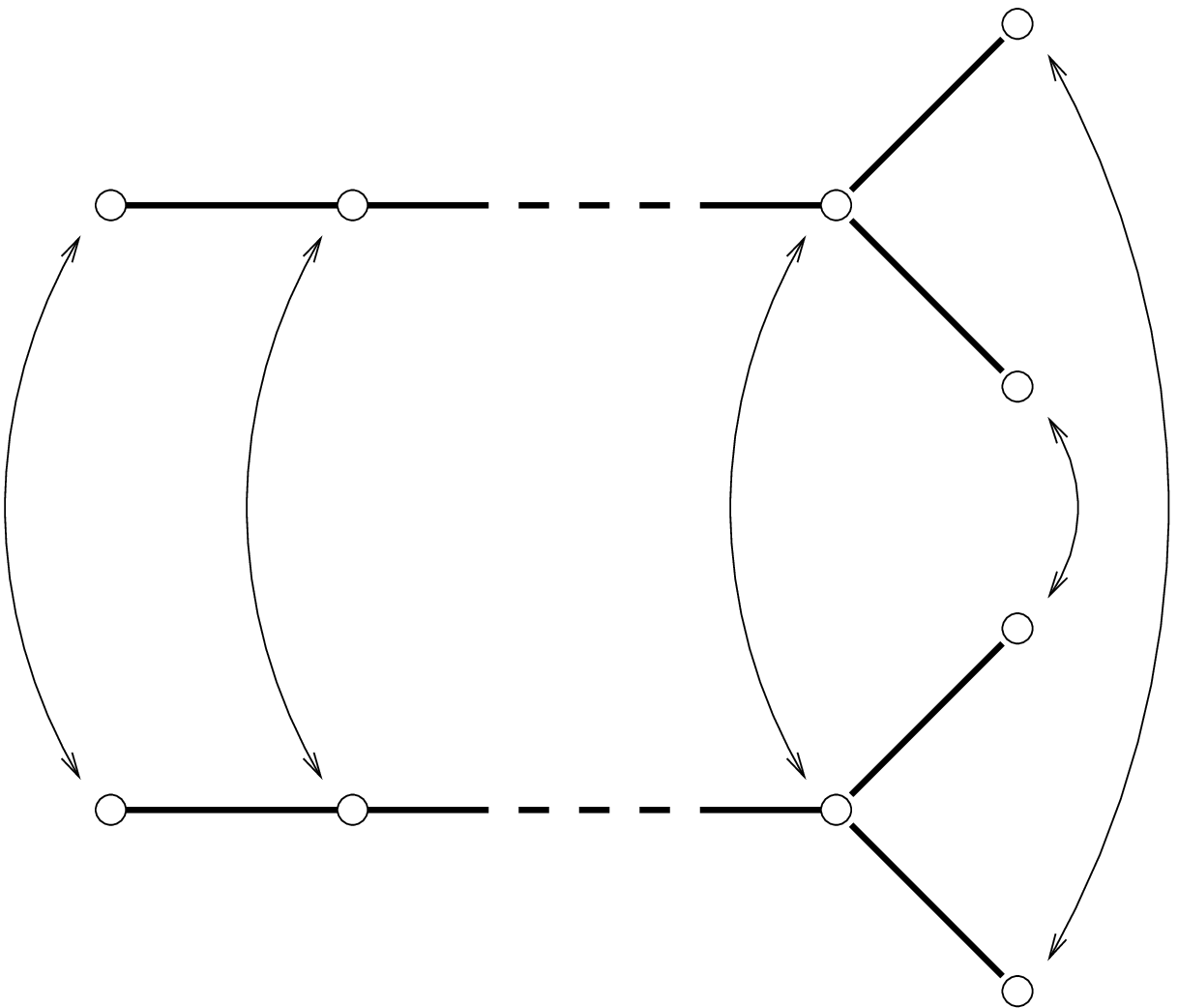}\end{center}\\
\hline
\end{longtable}
\end{small}

\subsubsection*{Strong $\Delta$-connectedness}

Let $(S^p,\Sigma,\mathbf A)$ be a spherical system and let $\Delta$ be its set of colours. Two spherical roots $\gamma_1,\gamma_2\in\Sigma$ are said to be \textit{strongly $\Delta$-adjacent} (\cite{L01}) if for all $\delta\in\Delta(\{\gamma_1\})$ we have $\langle\rho(\delta),\gamma_2\rangle\neq0$ and, vice versa, for all $\delta\in\Delta(\{\gamma_2\})$ we have $\langle\rho(\delta),\gamma_1\rangle\neq0$.

See Column~2 of Table~\ref{T2} and notice that, in type $\mathsf A$~$\mathsf D$, two spherical roots with overlapping supports are strongly $\Delta$-adjacent.

Let $\Sigma^\prime$ be a subset of $\Sigma$. The localization $(S^p\cap\supp(\Sigma^\prime),\Sigma^\prime,\mathbf A(S\cap \Sigma^\prime))$ of $(S^p,\Sigma,\mathbf A)$ in $\Sigma^\prime$ is said to be \textit{strongly $\Delta$-connected} if for every couple of spherical roots $\gamma_1,\gamma_2\in\Sigma^\prime$ there exists a finite sequence of spherical roots in $\Sigma'$, one strongly $\Delta$-adjacent to the next one, starting with $\gamma_1$ and ending with $\gamma_2$. If $\Sigma^\prime\subset\Sigma$ is maximal with this property we say that $(S^p\cap\supp(\Sigma^\prime),\Sigma^\prime,\mathbf A(S\cap \Sigma^\prime))$ is a \textit{strongly $\Delta$-connected component}\footnote{There exists a weaker property called \textit{$\Delta$-connectedness} also introduced by Luna in \cite{L01}, but its generalization to type $\mathsf A$~$\mathsf D$ doesn't seem to be as useful as in type $\mathsf A$.} of $(S^p,\Sigma,\mathbf A)$.

Let $(S^p,\Sigma,\mathbf A)$ be a spherical system. Let $\Sigma^\prime$ be a subset of spherical roots. Let $\Delta(\Sigma^\prime)$ denote the subset of colours $\delta\in\Delta(\supp(\Sigma^\prime))$ such that $\langle\rho(\delta),\gamma\rangle=0$ for all $\gamma\in\Sigma\setminus\Sigma^\prime$. We say that the localization $(S^p\cap\supp(\Sigma^\prime),\Sigma^\prime,\mathbf A(S\cap \Sigma^\prime))$ of $(S^p,\Sigma,\mathbf A)$ is \textit{erasable}\footnote{This requirement is slightly weaker than that found in \cite{L01}.} if there exists a nonempty smooth distinguished subset of colours $\Delta^\prime$ included in $\Delta(\Sigma^\prime)$. Notice that such a subset is smooth and distinguished with respect to the spherical system $(S^p,\Sigma,\mathbf A)$ if and only if it is smooth and distinguished with respect to the localization $(S^p\cap\supp(\Sigma^\prime),\Sigma^\prime,\mathbf A(S\cap \Sigma^\prime))$.

Weakening this requirement, we say that the localization ($S^p\cap\supp(\Sigma^\prime)$, $\Sigma^\prime$, $\mathbf A(S\cap \Sigma^\prime)$) is \textit{quasi-erasable} if there exists a nonempty distinguished subset of colours $\Delta^\prime$ included in $\Delta(\Sigma^\prime)$ and satisfying the $\ast$-property.

\begin{lemma}\label{L2a} Let $(S^p,\Sigma,\mathbf A)$ be a spherical system. Let $\Sigma_1$ and $\Sigma_2$ be two disjoint subsets of $\Sigma$ giving two quasi-erasable localizations such that at least one of them is erasable. Then the corresponding subsets $\Delta^\prime_1\subset\Delta(\Sigma_1)$ and $\Delta^\prime_2\subset\Delta(\Sigma_2)$ decompose the spherical system $(S^p,\Sigma,\mathbf A)$.\end{lemma}

\begin{proof} By definition, the subsets $\Delta^\prime_1$ and $\Delta^\prime_2$ are distinguished. \begin{enumerate} \renewcommand{\theenumi}{\roman{enumi}} \item They are different from the empty set as required and no colour can be in the intersection of $\Delta(\Sigma_1)$ and $\Delta(\Sigma_2)$ thanks to lemma \ref{lemma:colours}.\item The subsets $\Delta^\prime_1$ and $\Delta^\prime_2$ satisfy the $\ast$-property. Their union also satisfies the $\ast$-property since \[\frac{\Xi}{\Delta^\prime_1\cup\Delta^\prime_2}=\frac{\Xi}{\Delta^\prime_1}\cap\frac{\Xi}{\Delta^\prime_2}\] and \[\frac{\Sigma}{\Delta^\prime_1\cup\Delta^\prime_2}=\bigg(\Sigma\setminus(\Sigma_1\cup\Sigma_2)\bigg)\cup\frac{\Sigma_1}{\Delta^\prime_1}\cup\frac{\Sigma_2}{\Delta^\prime_2}.\] \item The subsets $\Sigma_1$ and $\Sigma_2$ are disjoint. \item Let us suppose that $\alpha\in \supp(\Sigma_1)$ and $\beta\in \supp(\Sigma_2)$ with $\alpha$ nonorthogonal to $\beta$: analyzing all the cases for $\alpha,\beta$ of type $a,a^\prime,b,p$ we get that neither $\Delta(\alpha)$ can be included in $\Delta(\Sigma_1)$ nor $\Delta(\beta)$ can be included in $\Delta(\Sigma_2)$. \item By the hypothesis, $\Delta^\prime_1$ or $\Delta^\prime_2$ is smooth. \end{enumerate} \end{proof}

A localization $(S^p\cap\supp(\Sigma^\prime),\Sigma^\prime,\mathbf A(S\cap \Sigma^\prime))$ of a strongly $\Delta$-connected component is said to be \textit{isolated} if the partition $\supp(\Sigma^\prime)\sqcup (\supp(\Sigma)\setminus \supp(\Sigma^\prime))$ gives a factorization of the localization in $\supp(\Sigma)$ of the spherical system $(S^p,\Sigma,\mathbf A)$. An isolated localization is a fortiori erasable.

\subsubsection*{Reduction}

The aim of this paragraph is to list all primitive spherical systems up to an automorphism of the Dynkin diagram of G.

The list of the primitive spherical systems of type $\mathsf A$ can be taken from \cite{L01} and it is included in Table~\ref{T3}; we have added two further cases which were missing in \cite{L01}.

The first step to obtain the complete list in type $\mathsf A$~$\mathsf D$ consists in finding all the cuspidal strongly $\Delta$-connected spherical systems. One can proceed by taking a spherical root, attaching to it any other spherical root to obtain a strongly $\Delta$-connected pair and going on recursively by attaching a new spherical root anywhere it is possible.

We obtain the following, where the labels are those of Table~\ref{T3}. Starting from a spherical root of type $d_m$: $ac(n)$ $n\geq3$ odd, $d(n)$ $n\geq4$ with only one spherical root of type $d_n$ ($ac(3)$ will be denoted also by $d(3)$), $do(p+q)$ $q\geq3$, $dc(n)$ $n$ even, $dc^\prime(n)$ $n$ even, $dc(n)$ $n$ odd, $dc^\sharp(n)$ (the case $D1(i)$ of Table~\ref{T2}). Starting from a root of type $a_1\times a_1$ without roots of type $d_m$: $aa(p+q+p)$, $aa(p,p)$, $aa^\ast(p+1+p)$, $do(p+2)$, $dd(p,p)$ $p\geq4$. Starting from a root of type $a^\prime_1$ without roots of type $d_m$ and $a_1\times a_1$: $ao(n)$, $do(n)$. Starting from a root of type $a_m$, $m\geq2$, with only roots of type $a_m$: $aa(n)$ $n\geq2$ with only one spherical root of type $a_q$, $ac^\ast(n)$ $n\geq3$, $dc^\ast(n)$, $ds(n)$, $ds^\ast(4)$.

Now we must find all cuspidal strongly $\Delta$-connected spherical systems with $\Sigma$ containing only spherical roots of type $a_1$. In this case we have $S^p=\emptyset$ and $\Sigma=S$.

\begin{lemma}\label{L3} If such a spherical system has a projective element $\delta$, then $\langle\rho(\delta),\alpha\rangle=1$ for all $\alpha\in\Sigma$, namely: \[\delta\in\bigcap_{\alpha\in\Sigma}\mathbf A(\alpha).\]\end{lemma}

\begin{proof} Two spherical roots $\alpha$ and $\beta$ of type $a_1$ are strongly $\Delta$-adjacent if and only if there is a colour $\delta$ in $\mathbf A(\alpha)\cap\mathbf A(\beta)$, namely we have the cases $A7(i)$ and $A8$ of Table~\ref{T2}. If $\alpha$ and $\beta$ are spherical roots of a spherical system with the above properties, this colour $\delta$ is the unique element of $\mathbf A$ that can be projective.
\end{proof}

We will call \textit{$n$-comb} a spherical system of rank $n$ satisfying the hypotheses of the above lemma.

In order to obtain all cuspidal strongly $\Delta$-connected spherical systems with only spherical roots of type $a_1$ without projective elements, one can proceed analogously as above. The list is the following: $aa(1)$ with only one spherical root of type $a_1$, $ax(1,1,1)$, $ay(p,p)$ $p\geq2$, $ay(p,(p-1))$ $p\geq2$, $ay(n)$ $n\geq4$ even, $az(3,3)$, $az(3,2)$, $az(3,1)$, $ae_6(6)$, $ae_6(5)$, $ae_7(7)$, $ae_7(6)$, $ae_7(5)$, $dy(p,p)$ $p\geq4$, $dy(p,(p-1))$ $p\geq4$, $dy(p,(p-2)$ $p\geq4$, $dz(4,4)$, $dz(4,3)$, $dz(4,2)$, $dz(4,1)$, $de_6(7)$, $de_6(6)$, $dy(7)$, $de_7(8)$, $de_7(7)$, $de_8(8)$, $de_8(7)$.

The next step will consist in gluing two or more strongly $\Delta$-connected components together. The following consideration about the properties of being erasable or isolated will be useful to simplify this procedure.

Let $(S^p,\Sigma,\mathbf A)$ be a cuspidal spherical system and let $\Sigma^\prime$ give a strongly $\Delta$-connected component $(S^p\cap\supp(\Sigma^\prime),\Sigma^\prime,\mathbf A(S\cap\Sigma^\prime))$. There are only the following two situations where an element $\delta\in\mathbf A$ may correspond to a colour not belonging to $\Delta(\Sigma^\prime)$: \begin{itemize}\item[-] The colour $\delta$ is said to be \textit{external} if there is a unique simple root $\alpha$ with $\delta\in\mathbf A(\alpha)$ and it is not in $\Delta(\Sigma^\prime)$ since it takes a negative value on a spherical root $\gamma$ outside $\Sigma^\prime$. Then the root $\alpha$ is nonorthogonal to $\gamma$. \item[-] The colour $\delta$ is said to be a \textit{bridge} if there are more than one simple root, let us say $\alpha_i$ for $i=1,\ldots,r$, such that $\delta$ is in $\mathbf A(\alpha_i)$, and it takes a negative value on a spherical root $\gamma$ outside $\Sigma^\prime$. Then all the roots $\alpha_i$ are nonorthogonal to $\gamma$, which is necessarily of type $a_m$, and $r$ is equal to $2$ or $3$.\end{itemize}

In most cases a look at the set of colours of a strongly $\Delta$-connected component easily shows that it is necessarily isolated, erasable or quasi-erasable.

The strongly $\Delta$-connected component $(S^p\cap\supp(\Sigma^\prime),\Sigma^\prime,\mathbf A(S\cap\Sigma^\prime))$ is isolated if it is isomorphic to one of these: $do(p+q)$, $dc(n)$ for $n\geq6$ even, $dc^\prime(n)$, $dc(n)$ for $n$ odd, $aa(p+q+p)$ for $p\geq2$, $aa(p,p)$ for $p\geq2$, $aa^\ast(p+1+p)$, $do(p+2)$, $dd(p,p)$, $ao(n)$, $do(n)$, $ay(n)$ for $n\geq6$, $dy(p,p)$ for $p\geq5$, $de_6(6)$, $dy(7)$, $de_7(7)$, $de_8(7)$.

If the strongly $\Delta$-connected component $(S^p\cap\supp(\Sigma^\prime),\Sigma^\prime,\mathbf A(S\cap\Sigma^\prime))$ is isomorphic to $ac(n)$ for $n\geq5$, $aa(1+q+1)$ for $q\geq2$, $ay(4)$ or $ae_7(5)$ then either it is isolated or $\supp(\Sigma^\prime)$ is nonorthogonal to the support of an $n$-comb and $(S^p,\Sigma,\mathbf A)$ has a projective element.

The strongly $\Delta$-connected component $(S^p\cap\supp(\Sigma^\prime),\Sigma^\prime,\mathbf A(S\cap\Sigma^\prime))$ is erasable if it is isomorphic to one of these: $dc^\sharp(n)$, $aa(1+1+1)$, 
$n$-comb for $n\geq4$, $az(3,3)$, $ae_6(6)$, $ae_6(5)$, $ae_7(7)$, $ae_7(6)$, $dy(4,4)$, $dz(4,4)$, $dz(4,3)$, $dz(4,2)$, $dz(4,1)$, $de_6(7)$, $de_7(8)$, $de_8(8)$.

The strongly $\Delta$-connected component $(S^p\cap\supp(\Sigma^\prime),\Sigma^\prime,\mathbf A(S\cap\Sigma^\prime))$ is quasi-erasable if it is isomorphic to one of these: $dc^\ast(n)$, $ds(n)$, $ds^\ast(4)$.

We examine the remaining strongly $\Delta$-connected components case by case: \begin{itemize}
\item[-] $d(n)$ for $n\geq3$. If it is not isolated, it is neither erasable nor quasi-erasable.
\item[-] $aa(1,1)$. Analogously, if it is not isolated, it is neither erasable nor quasi-erasable. The strongly $\Delta$-connected component $aa(1,1)$ will be denoted also by $d(2)$, since it can be considered as limit case of $d(n)$ for $n>2$.
\item[-] $aa(n)$ for $n\geq2$. It is erasable only if it lies at one end of a connected component (in the usual sense) of the Dynkin diagram.
\item[-] $aa(1)$. If it has no projective elements then the functionals of its two colours take a negative value respectively on two distinct spherical roots. Hence, it does not lie at the end of a connected component of the Dynkin diagram or the spherical system admits the following erasable localizations: $aa(1)+ac^\ast(n)$, $aa(1)+ac^\flat(n)$, $aa(1)+ac^{\flat\flat}(4)$ (see Table~\ref{T3}).
\item[-] $ac^\ast(n)$. In case its support contains a ``bifurcation'' of the Dynkin diagram, i.e.\ a simple root nonorthogonal to three simple roots, and there are no projective elements, we have necessarily one of the following two erasable localizations: $ac^\ast(n)+aa(1)$, $ac^\ast(3)+2\mathrm{-comb}$ (see Table~\ref{T3}). In case it does not lie on a bifurcation, it is erasable if at the end of a connected component of the Dynkin diagram or if $n$ is odd, it is quasi-erasable if $n$ is even.
\item[-] $3$-comb. It has a projective element or its support is made of three mutually orthogonal simple roots joined by a bridge. In the latter case it is part of one of the following erasable localizations: $3\mathrm{-comb}+2\mathrm{-comb}$, $3\mathrm{-comb}+aa(n)$ (see Table~\ref{T3}).
\item[-] $2$-comb. Analogously, it has a projective element or its support is made of two orthogonal simple roots joined by
 a bridge.
\item[-] $ax(1,1,1)$. If there are no projective elements then either it lies at the end of a connected component of the Dynkin diagram and it is quasi-erasable, either it is part of $ax(1+p+1,1)$ that is quasi-erasable or it is part of $af(5)$ that is isolated.
\item[-] $ay(p,p)$. There are different colours that may be a bridge. If there is no bridge then it is erasable. Otherwise it is part of one of the following erasable localizations: $ay(p+q+p)$, $dy^\ast(3+q+3)$, $ay(2,2)+2\mathrm{-comb}$ (see Table~\ref{T3}).
\item[-] $ay(p,p-1)$. If there is no bridge then it is quasi-erasable, or erasable if at the end of a Dynkin diagram. Otherwise it is part of one of the following localizations: $ay(p+q+(p-1))$ erasable, $dy^\ast(3+q+2)$ isolated, $ay(2,1)+2\mathrm{-comb}$ isolated (see Table~\ref{T3}).
\item[-] $az(3,2)$. If there is no bridge then it is erasable. Otherwise it is part of one of the following localizations: $dy(3+q+2)$ isolated, $dz(3+q+2)$ erasable (see Table~\ref{T3}).
\item[-] $az(3,1)$. There is only one colour that may be a bridge and in this case we can have only $dz(3+q+1)$ that is isolated, if this colour is not a bridge then $az(3,1)$ is quasi-erasable or even erasable if it does not lie on a bifurcation.
\item[-] $dy(p,p-1)$. There is only one colour that may be a bridge and in this case we can have only $dy(p+q+(p-1))$ that is isolated, if this colour is not a bridge the strongly $\Delta$-connected component itself is erasable.
\item[-] $dy(p,p-2)$. There is only one colour that may be a bridge and in this case we can have only $dy^\prime(p+q+(p-2))$ that is isolated, if this colour is not a bridge the strongly $\Delta$-connected component itself is quasi-erasable.
\end{itemize}

We are ready to complete the list of primitive spherical systems: we have to ``put together'' the components obtained above in all possible ways such that the obtained system turns out to be indecomposable, without projective colours. Thanks to Lemma~\ref{L2a} and subsequent considerations, the procedure is quite straightforward and leads to the list of Table~\ref{T3}.

\subsubsection*{Example}

Let us start with $az(3,3)$: first of all, it is a primitive system itself. We go on adding other components, taking into account the fact that $az(3,3)$ is erasable whenever it is a part of a larger system. Therefore we may suppose that there is no other strongly $\Delta$-connected component (more precisely, no union of such components) being quasi-erasable.

The list at the end of the preceding paragraph gives us the other possible components: $d(n)$, $aa(n)$, $2$-comb. A copy of $d(n)$ or $aa(n)$ can be added to $az(3,3)$, obtaining the systems: $az^\sim(3+q+3)$, $az(3,3)+d(n)$. An additional copy of $d(n)$ gives the systems $az^\sim(3+q+3)+d(n)$ and $az(3,3)+d(n_1)+d(n_2)$.

All other possibilities (including those containing a $2$-comb) bring some union of the added components being quasi-erasable.

\section{The proof of the primitive cases}

In this section we complete the proof of Theorem \ref{thmAD}, checking that every primitive spherical system corresponds to a wonderful variety, unique up to $G$-isomorphism. We divide the proof for each case in two parts. In the ``uniqueness'' part we suppose that there exists a wonderful subgroup $H$ such that the wonderful embedding of $G/H$ has the required spherical system, and we prove that $H$ is uniquely determined up to conjugation. This will provide a candidate $H$ and in the ``existence'' part we prove that it is wonderful and that it actually corresponds to our spherical system.

For all cases where $G$ is of type $\mathsf A$ this has already been proven in \cite{L01}
, thus we suppose $G$ be semisimple adjoint of mixed type $\mathsf A$~$\mathsf D$, not of type $\mathsf A$. Wherever it is convenient, we will work on a covering of $G$ rather than on $G$ itself.

The next two technical lemmas will be widely used through all proofs.

\begin{lemma}[\cite{L01}] \label{L4} Let $H$ be a wonderful subgroup of $G$ and let $(S^p,\Sigma,\mathbf A)$ be its spherical system. Then the following identities hold: \begin{enumerate} \renewcommand{\theenumi}{\roman{enumi}} \item\label{l} $\dim G/P_{S^p} + \rank \Sigma = \dim G/H$, \item\label{ll} $\rank \Xi(H)= \rank \Delta - \rank \Sigma$. \end{enumerate}\end{lemma}

\begin{lemma}\label{L5} \label{L4a} Let $(S^p,\Sigma,\mathbf A)$ be a spherical system of type $\mathsf A$~$\mathsf D$. \begin{enumerate} \renewcommand{\theenumi}{\roman{enumi}} \item In general, $\rank\Sigma\leq \rank S$ and $\dim G/P_{S^p}=\dim P_{S^p}^u\leq \dim B^u$. \item If $\Sigma$ contains a spherical root of type $d_m$ or $a_1\times a_1$ then $\rank \Sigma < \rank S$. \item \label{am} If $\Sigma$ contains a spherical root of type $a_m$, for any $m\geq1$, then $\rank \Sigma < \rank \Delta$ or $(S^p,\Sigma,\mathbf A)$ can be localized to $ax(1,1,1)$ or $ds^\ast(4)$. \item \label{aone} Every strongly $\Delta$-connected component with only spherical roots of type $a_1$ different from $ax(1,1,1)$ and $ds^\ast(4)$ is such that $\rank \Delta-\rank\Sigma=1$. \end{enumerate}\end{lemma}

\begin{proof} All the statements are quite obvious; in order to prove \textit{\ref{am}} in case of $m=1$ and \textit{\ref{aone}} it is sufficient to check them for all cuspidal strongly $\Delta$-connected spherical systems with only spherical roots of type $a_1$. \end{proof}

\subsection{Reductive cases}\label{reductivecases}

Let $H$ be a spherical subgroup of $G$. Then the homogeneous space $G/H$ is affine, or equivalently the subgroup $H$ is reductive, if and only if there exists a linear combination \[\xi=\sum_{\gamma\in\Sigma_{G/H}}c_\gamma\,\gamma\] with $c_\gamma\geq0$ such that $\langle\rho(\delta),\xi\rangle>0$ for all $\delta\in\Delta_{G/H}$. The primitive systems which satisfy this property are: $do(p+q)$, $dc(n)$ for $n$ even, $dc^\prime(n)$, $dc(n)$ for $n$ odd, $do(p+2)$, $dd(p,p)$, $do(n)$, $ds^\ast(4)$ and $aa(1)+ds^\ast(4)$.

On the other hand, reductive spherical subgroups are classified (\cite{Kr79}, \cite{B87}); we report in Table \ref{T4} the list of the indecomposable connected reductive spherical subgroups.

\begin{small}
\begin{longtable}{|p{3.3cm}p{3.8cm}p{4cm}|}
\caption{Indecomposable connected reductive spherical subgroups of semisimple groups.}\label{reductive spherical subgroups} 
\label{T4} \\
\hline
$H$ & $G$ & \\
\hline
\hline
\endfirsthead
\caption{Indecomposable connected reductive spherical subgroups of semisimple groups (continuation).} \\
\hline
$H$ & $G$ & \\
\hline
\hline
\endhead
\hline
\endfoot
$SO(n)$ & $SL(n)$ & $n\geq2$ \\
$S(GL(n)\times GL(m))$ & $SL(n+m)$ & $n\geq m\geq1$ \\
$SL(n)\times SL(m)$ & $SL(n+m)$ & $n>m\geq1$ \\
$Sp(2n)$ & $SL(2n)$ & $n\geq2$ \\
$Sp(2n)$ & $SL(2n+1)$ & $n\geq1$ \\
$GL(1)\cdot Sp(2n)$ & $SL(2n+1)$ & $n\geq1$ \\
\hline
$GL(n)$ & $SO(2n)$ & $n\geq2$ \\
$SL(n)$ & $SO(2n)$ & $n\geq3$ odd \\
$GL(n)$ & $SO(2n+1)$ & $n\geq2$ \\
$SO(n)\times SO(m)$ & $SO(n+m)$ & $n\geq m\geq1$ and $n+m\geq 3$ \\
$Spin(7)$ & $SO(8)$ & \\
$Spin(7)$ & $SO(9)$ & \\
$G_2$ & $SO(7)$ & \\
$G_2$ & $SO(8)$ & \\
$SL(2)\cdot Sp(4)$ & $SO(8)$ & \\
$SO(2)\times Spin(7)$ & $SO(10)$ & \\
\hline
$GL(n)$ & $Sp(2n)$ & $n\geq1$ \\
$Sp(2n)\times Sp(2m)$ & $Sp(2n+2m)$ & $n\geq m\geq1$ \\
$GL(1)\times Sp(2n-2)$ & $Sp(2n)$ & $n\geq1$ \\
\hline
$A_2$ & $G_2$ & \\
$A_1\times A_1$ & $G_2$ & \\
\hline
$B_4$ & $F_4$ & \\
$C_3\times A_1$ & $F_4$ & \\
\hline
$C_4$ & $E_6$ & \\
$F_4$ & $E_6$ & \\
$D_5$ & $E_6$ & \\
$GL(1)\cdot D_5$ & $E_6$ & \\
$A_5\times A_1$ & $E_6$ & \\
$GL(1)\cdot E_6$ & $E_7$ & \\
$A_7$ & $E_7$ & \\
$D_6\times A_1$ & $E_7$ & \\
$D_8$ & $E_8$ & \\
$A_1\times E_7$ & $E_8$ & \\
\hline
$H$ & $H\times H$ & $H$ simple group embedded diagonally in $H\times H$ \\
\hline
$SL(2)$ & $SL(2)\times SL(2)\times SL(2)$ & embedded diagonally \\
\hline
$SL(2)\times Sp(2n)$ & $SL(2)\times Sp(2n+2)$ & $(u,v)\mapsto (u,u\oplus v)$ \\
\hline
$Sp(4)\times Sp(2n)$ & $Sp(4)\times Sp(2n+4)$ & $(u,v)\mapsto (u,u\oplus v)$ \\
\hline
$SO(n)$ & $SO(n)\times SO(n+1)$ & $u\mapsto (u,i(u))$ \\
\hline
$Spin(7)$ & $Spin(7)\times SO(8)$ & $u\mapsto (u,i(u))$ \\
\hline
\mbox{$SL(2)\times Sp(2m)\times$} \mbox{$\quad \times Sp(2n)$} & $Sp(2m+2)\times Sp(2n+2)$ & $(t,u,v)\mapsto (t\oplus u,t\oplus v)$ \\
\hline
\mbox{$SL(2)\times SL(2)\times$} \mbox{$\quad \times Sp(2m)\times Sp(2n)$} & \mbox{$Sp(4)\times Sp(2m+2)\times$} \mbox{$\quad \times Sp(2n+2)$} & \mbox{$(t,u,v,w)\mapsto$} \mbox{$\quad \mapsto(t\oplus u,t\oplus v,u\oplus w)$} \\
\hline
\mbox{$SL(2)\times Sp(2l)\times$} \mbox{$\quad \times Sp(2m)\times Sp(2n)$} & \mbox{$Sp(2l+2)\times Sp(2m+2)\times$} \mbox{$\quad \times Sp(2n+2)$} & \mbox{$(t,u,v,w)\mapsto$} \mbox{$\quad\mapsto (t\oplus u,t\oplus v,t\oplus w)$} \\
\hline
$SL(n)\times GL(1)$ & $SL(n)\times SL(n+1)$ & $(u,\lambda)\mapsto (u,\lambda u\oplus \lambda^{-n})$ \\
\hline
\mbox{$SL(2)\times SL(m)\times$} \mbox{$\quad \times Sp(2n)$} & $SL(m+2)\times Sp(2n+2)$ & $(u,v,w)\mapsto (u\oplus v,u\oplus w)$ \\
\hline
\mbox{$SL(2)\times SL(m)\times$} \mbox{$\quad \times Sp(2n)\times GL(1)$} & $SL(m+2)\times Sp(2n+2)$ & \mbox{$(u,v,w,\lambda)\mapsto$} \mbox{$\quad\mapsto (\lambda^{m}u\oplus \lambda^{-2}v,u\oplus w)$} \\
\hline
\end{longtable}
\end{small}

In particular, for $G$ of type $\mathsf D$ we have the following indecomposable connected reductive spherical subgroups.

\begin{itemize}
\item[-] $GL(n)\subset SO(2n)$, $n\geq2$, $\dim G/H=n(n-1)$, $\rank \Xi(H)=1$.
\item[-] $SL(n)\subset SO(2n)$, $n\geq3$ odd, $[N(H):H]$ infinite.
\item[-] $SO(n)\times SO(m)\subset SO(n+m)$, $n\geq m\geq1$, $n+m\geq4$ even, $\dim G/H =nm$, $\rank\Xi(H)$ equal to $0$ if $m,n\neq2$, $1$ if $m=2$ or $n=2$, $2$ if $m=n=2$.
\item[-] $Spin(7)\subset SO(8)$, $\dim G/H=7$, $\rank\Xi(H)=0$.
\item[-] $G_2\subset SO(8)$, $\dim G/H=14$, $\rank\Xi(H)=0$.
\item[-] $SL(2)\cdot Sp(4)\subset SO(8)$, $\dim G/H=15$, $\rank\Xi(H)=0$.
\item[-] $SO(2)\times Spin(7)\subset SO(10)$, $\dim G/H=23$, $\rank\Xi(H)=1$.
\item[-] $SO(2n)\subset SO(2n)\times SO(2n)$ embedded diagonally, $n\geq2$, $\dim G/H=n(2n-1)$, $\rank\Xi(H)=0$.
\end{itemize}

This table provides a candidate $H$ for each of the above primitive systems, or at least a candidate for the connected component $H^\circ$ of $H$ containing the identity. In the latter case we have to guess $H$, knowing that $H^\circ\subset H\subset N_G(H)=N_G(H^\circ)$, where the last equality is known from the theory of spherical varieties.

In general, the normalizer $N_G(H)$ of a wonderful subgroup $H$ acts on the colours of $G/H$ and the set of the colours of $G/N_G(H)$ is a quotient of the set of the colours of $G/H$ such that the colours permuted by $N_G(H)$ are identified. The normalizer $N_G(H)$ can permute only the colours having the same associated functional, and we have that the only case where two colours can have the same functional is the case of two colours $\delta^+_{\alpha}$ and $\delta^-_{\alpha}$ moved by the same simple root $\alpha$ being also a spherical root.

This fact can be easily deduced from the list of the rank two wonderful varieties, and in type $\mathsf A$~$\mathsf D$ it holds that the spherical system associated $N_G(H)$, if it is different from the spherical system associated to $H$, can be obtained from the latter by substituting some spherical roots $\alpha\in S\cap\Sigma$ with the spherical roots $2\alpha$ (see also Proposition~\ref{rigidity}).

Now, once we have a candidate $H$, ad-hoc considerations are needed in order to prove that $H$ satisfy our requirements, with the help of lemmas \ref{L4} and \ref{L5} and the stringent restrictions imposed by the axioms of spherical systems.


$\mathbf{do(p+q)}$ ($p\geq1$, $q\geq3$). To see that $H$ must be reductive here we can take $\xi=c_1\,2\alpha_1+\ldots+c_p\, 2\alpha_p+c_{p+1}\,\gamma_{p+1}$, where $\Sigma=\{2\alpha_1,\ldots,2\alpha_p,\gamma_{p+1}\}$, with $c_i=f(i)$ for any concave function of the interval $[1,p+1]$ with $2f(1)>f(2)$ and $f(p)<f(p+1)$, for example $c_i=i+1-(1/(p+2-i))$. For the next cases we will omit this easy check that $H$ is reductive.

{\em Uniqueness.} We have here $G=SO(2p+2q)$. The spherical system is such that $\dim G/P_{S^p}=\dim P_{S^p}^u =(p+1)(2q+p-2)$, $\rank \Sigma=p+1$ and $\rank \Delta=p+1$. Hence by Lemma~\ref{L4} we have $\dim G/H=(p+1)(2q+p-1)$ and $\rank\Xi(H)=0$.
Looking at Table~\ref{T4} we see that these conditions force the connected component $H^\circ$ of the identity in $H$ to be conjugated to $SO(p+1)\times SO(2q+p-1)$, or $GL(p+q)$ if $p+q=(q-1)^2$. Let us prove that the latter choice is impossible. If $p+q$ is even then $[N_G(GL(p+q):GL(p+q)]=2$. Now $H$ cannot be isomorphic to $GL(p+q)$ since $\rank \Xi(GL(p+q))=1$. The only other chance would be $H=N_G(GL(p+q)$ up to conjugation, but the existence part of the proof for $dc^\prime(n)$ (page \pageref{dcprimen}) shows that $N_G(GL(p+q))$ is wonderful and the associated spherical system is of type $dc^\prime(n)$. If $p+q$ is odd $N_G(GL(p+q))=GL(p+q)$ and $\rank \Xi(GL(p+q))=1$. Therefore, $H^\circ$ is conjugated to $SO(p+1)\times SO(2q+p-1)$. The normalizer is $S(O(p+1)\times O(2q+p-1))$ and $[N(H^\circ):H^\circ]=2$, hence $H$ can only be $H^\circ$ or $N(H^\circ)$. Now, we will see in the existence part of the proof that $N(H^\circ)$ is wonderful and the associated spherical system is actually $do(p+q)$. By the results on spherically closed subgroups in \cite{L01}
the lattice $\Xi_{G/H^\circ}$ contains properly $\Xi_{G/H}$, hence $H^\circ$ and $N(H^\circ)$ cannot have the same spherical system. This forces $H=N(H^\circ)$.

{\em Existence.} Let us put $H=S(O(p+1)\times O(2q+p-1))$. Since $\rank \Xi(H)=0$, we have $\rank \Sigma=\rank \Delta$, hence the associated spherical system is cuspidal and, by Lemma~\ref{L4a}, the set $\Sigma$ does not contain roots of type $a_m$ for $m\geq1$. Indeed, the case of a component of type $ds^\ast(4)$ can be excluded directly by noticing that there is a spherical root of type $a_3$ that can be nonorthogonal only to a spherical root of type $a_m$ (see Table~\ref{T2}). Moreover, $\dim P^u_{S^p}+\rank \Sigma=\dim G/H=(p+q)(p+q-1)-q(q-2)$. In the case of $S^p=\emptyset$, we would have $\rank \Sigma=p+q-q(q-2)$ since $P^u_{S^p}=B^u$ and no roots of type $d_m$, but this is impossible. Hence, we have that $S^p\neq\emptyset$ and $\Sigma$ contains at least one root of type $d_m$. If this root does not lie at the end of the diagram like in $do(p+q)$, we get $dc^\prime(p+q)$, since a spherical root of type $d_3$ on a Dynkin diagram of type $\mathsf A$ can be nonorthogonal only to a root of the same type. The case $dc^\prime(p+q)$ can be excluded since it is not invariant under the exchange of the simple roots $\alpha_{p+q-1}$ and $\alpha_{p+q}$, while if consider the corresponding outer involutive automorphism $\sigma$ we have $\sigma(H)=H$. Hence, there is only one root of type $d_m$, it lies at the end of the diagram and the unique possible spherical system is $do(p+q)$.

$\mathbf{do(p+2)}$. Suppose at first $n=p+2>4$. {\em Uniqueness.} This case is similar to $do(p+q)$, here we can take $p=n-2$ and $q=2$. The subgroup $H$ must be reductive. We have $\dim G/H=(n-1)(n+1)$ and $\rank\Xi(H)=0$. If $G=SO(2n)$, by Table~\ref{T4} we get that $H^\circ\cong SO(n-1)\times SO(n+1)$, and analogously to $do(p+q)$ we get that $H=N_G(H^\circ)=S(O(n-1)\times O(n+1))$.

{\em Existence.} $H$ is wonderful. The identities $\dim P^u_{S^p}+\rank \Sigma=(n-1)(n+1)$ and $\rank \Delta-\rank \Sigma=0$ force $\rank \Sigma\geq n-1$. In general $\rank\Sigma\leq n$. If $\rank\Sigma= n$ then $\rank S^p=1$, but no root of type $a_m$ can belong to $\Sigma$, hence $\rank \Sigma= n-1$ and $S^p=\emptyset$. We obtain the spherical system $do(p+2)$.

Let $n=4$. Notice that in this case we have three different spherical systems of rank $3$.
\[\begin{picture}(16050,3600)\put(300,2100){\usebox{\segm}}\multiput(0,1200)(1800,0){2}{\usebox{\aprime}}\put(2100,600){\usebox{\plusaoneaone}}\put(6750,2100){\usebox{\segm}}\multiput(8250,1200)(1200,1200){2}{\usebox{\aprime}}\put(8550,900){\usebox{\bifurc}}\put(12750,2100){\usebox{\segm}}\multiput(14250,1200)(1200,-1200){2}{\usebox{\aprime}}\put(14550,900){\usebox{\bifurc}}\put(6750,2100){\circle{600}}\put(12750,2100){\circle{600}}\put(9750,900){\circle{600}}\put(15750,3300){\circle{600}}\put(6750,1800){\line(0,-1){900}}\put(6750,900){\line(1,0){2700}}\put(12750,2400){\line(0,1){900}}\put(12750,3300){\line(1,0){2700}}\end{picture}\]
We have $\dim G/H=15$ and $\rank \Xi(H)=0$. Our $G$ is $SO(8)$, and by Table~\ref{T4} we get that $H^\circ$ can be isomorphic only to $SO(3)\times SO(5)$ or $SL(2)\cdot Sp(4)$. In the first case we have that $N_{SO(8)}(SO(3)\times SO(5))=(SO(3)\times SO(5))\cdot C_{SO(8)}$ as above. In the second case we have that $SL(2)\cdot Sp(4)\cong Spin(3)\cdot Spin(5) \subset Spin(7)\subset SO(8)$ is equal to its normalizer. Let us consider the outer involutive automorphism $\sigma$ of $SO(8)$ that exchanges the simple roots $\alpha_3$ and $\alpha_4$: we have that $\sigma(S(O(3)\times O(5)))=S(O(3)\times O(5))$ while $\sigma(SL(2)\cdot Sp(4))$ is isomorphic to $SL(2)\cdot Sp(4)$ but they are not conjugated. Therefore, we have three possible subgroups up to conjugation: $S(O(3)\times O(5))$, $SL(2)\cdot Sp(4)$ and $\sigma(SL(2)\cdot Sp(4))$.

They are spherical and hence wonderful, and their spherical systems are all of type $do(p+2)$ (one proceeds like in the case $n\geq4$). The first subgroup has the first spherical system of the above figure, since it is the unique one fixed by $\sigma$. The other subgroups are analogously in correspondence with the remaining two cases.

$\mathbf{do(n)}$. {\em Uniqueness.} We have $\dim G/H=n^2$ and $\rank \Xi(H)=0$. If $G=SO(2n)$, by Table~\ref{T4} we get $H^\circ\cong SO(n)\times SO(n)$, and $[N_G(H^\circ),H^\circ]=4$. The subgroup $N_G(H^\circ)$ is wonderful and, thanks to the existence part, the spherical system associated to it is $do(n)$; we conclude that $H$ must be $N_G(H^\circ)$ as in case $do(p+q)$.

{\em Existence.} From the identities $\dim P^u_{S^p}+\rank \Sigma=n^2$ and $\rank\Delta-\rank\Sigma=0$ we get $S^p=\emptyset$, $\rank\Sigma=n$ and $\rank \Delta=n$. Therefore, we get the spherical system $do(n)$.

$\mathbf{dc(n)}$ ($n\geq4$ even). {\em Uniqueness.} By Lemma~\ref{L4} we have $\dim G/H=(n(n-1)-n/2)+n/2=n(n-1)$ and $\rank \Xi(H)=1$. Our $G$ is $SO(2n)$, from Table~\ref{T4} we see that if $n>4$ then $H^\circ\cong GL(n)$ and $H=H^\circ$, since $\rank \Xi(N_G(GL(n)))=0$.

{\em Existence.} The spherical subgroup $H$ has index $2$ in its normalizer $N_G(H)$. Let us consider the functions $f_1,f_2\in\mathbb C[SO(2n)]$ defined by the minors $(n+1\ \ldots\ 2n\,|\,1\ \ldots\ n)$ and $(n+1\ \ldots\ 2n\,|\,n+1\ \ldots\ 2n)$, respectively. Since they are $(B\times H)$-proper and they are exchanged by the right action of $N_G(H)$, there are two colours of $G/H$ exchanged by $N(H)$. Therefore $H$ is spherically closed and hence wonderful (\cite{L01},\cite{Kn96}). We have $\dim P^u_{S^p}+\rank \Sigma=n(n-1)$ and $\rank \Delta-\rank \Sigma=1$. Hence, $S^p$ is non-empty and there exists at least a spherical root of type $d_m$, since by the above identities any root of type $a_m$ with $m\geq3$ cannot be in $\Sigma$. It remains only the spherical system $dc(n)$.

$\mathbf{dc^\prime(n)}$ ($n\geq4$ even)\label{dcprimen}. {\em Uniqueness.} Analogously, $\dim G/H=n(n-1)$ and $\rank\Xi(H)=0$. Here $G=SO(2n)$; by Table~\ref{T4} and, by the existence part of the proof for the case $do(p+q)$, we get $H\cong N_G(GL(n))$.

{\em Existence.} The subgroup $H$ is wonderful, since it is equal to its normalizer. Moreover, we already know the spherical system of $H^\circ$ and hence we get the case $dc^\prime(n)$.

$\mathbf{dc(n)}$ ($n\geq5$ odd). {\em Uniqueness.} We have $\dim G/H=n(n-1)$ and $\rank\Xi(H)=1$. Here $G=SO(2n)$; by Table~\ref{T4} we get $H^\circ\cong GL(n)$. Moreover, if $n$ is odd then $GL(n)=N_G(GL(n))$ , hence $H=H^\circ$.

{\em Existence.} $H$ is wonderful. We have $\dim P^u_{S^p}+\rank \Sigma=n(n-1)$ and $\rank \Delta-\rank \Sigma=1$. Hence, $S^p\neq\emptyset$ and there exists at least a spherical root of type $a_m$ with $m\geq3$, since we cannot obtain any spherical system satisfying the above identities with only roots of type $d_m$. There are no other choices than the spherical system $dc(n)$.

$\mathbf{dd(p,p)}$ ($p\geq4$). {\em Uniqueness.} We have $\dim G/H=p(2p-1)$. We have $G=SO(2p)\times SO(2p)$, by Table~\ref{T4} we get that $H^\circ$ is equal to $SO(2p)$ embedded diagonally in $G$. The normalizer $N_G(H^\circ)$ of $H^\circ$ is $C_G\cdot H^\circ$, therefore $H$ is $N_G(SO(2p))$.

{\em Existence.} From $\rank\Delta=\rank\Sigma$ follows that the spherical system of $H$ is cuspidal and without any root of type $a_m$. Since $H$ is not a direct product, there exists a root of type $a_1\times a_1$ whose support is not included in one single connected components of the Dynkin diagram. The spherical system is forced to be $dd(p,p)$.

$\mathbf{ds^\ast(4)}$. {\em Uniqueness.} We have $\dim G/H=14$ and $\rank\Xi(H)=0$. Our $G$ is $SO(8)$, by Table~\ref{T4} we get $H^\circ\cong G_2$. Moreover, we have $N_{SO(8)}(G_2)=G_2\cdot C_{SO(8)}$.

{\em Existence.} By the identities $\dim P^u_{S^p}+\rank \Sigma=\dim B^u +2$ and $\rank\Delta=\rank\Sigma$ we get $\rank S^p=1$ and hence at least one root of type $a_3$. Therefore, as seen in Lemma~\ref{L4a}, we get the spherical system $ds^\ast(4)$.

$\mathbf{aa(1)+ds^\ast(4)}$. {\em Uniqueness.} We have $\dim G/H=15$ and $\rank\Xi(H)=1$. If $G=SO(10)$, by Table~\ref{T4} we get $H^\circ\cong SO(2)\times Spin(7)$. Moreover, we have $N_{SO(10)}(H^\circ)=H^\circ$.

{\em Existence.} By the identities $\dim P^u_{S^p}+\rank \Sigma=\dim B^u +3$ and $\rank\Delta=\rank\Sigma+1$ we get that $\rank S^p=1$ and there are roots of type $a_m$. In particular, the only possibility is to take three roots of type $a_3$ and one root of type $a_1$. Hence, we must attach to $ds^\ast(4)$ a root of type $a_1$. There are two possibilities: in one case $H$ would be included, by Proposition~\ref{morphisms}, in the parabolic induction of $G_2$, but this is false. We get the spherical system $aa(1)+ds^\ast(4)$.



\subsection{Nonreductive cases}

As above, for each case we prove uniqueness and existence of the associated wonderful subgroup $H$. For the uniqueness part we use the following strategy, with some marginal changes. It is the typical argument used by Luna in \cite{L01}.

Let $(S^p,\Sigma,\mathbf A)$ be a primitive spherical system for $G$ and let us suppose that there exists a subgroup $H$ of $G$ such that the wonderful embedding of $G/H$ exists and has $(S^p,\Sigma,\mathbf A)$ as its spherical system.

{\em First step.} We find a distinguished subset $\Delta^\prime$ of $\Delta$ with the $\ast$-property such that $(S^p,\Sigma,\mathbf A)/\Delta^\prime$ is a parabolic induction of a reductive case and there exists a minimal parabolic distinguished subset $\Delta^{\prime\prime}$ such that it contains $\Delta^\prime$.

The quotient $(S^p,\Sigma,\mathbf A)/\Delta^\prime$ corresponds to a unique wonderful variety $X_1$ thanks to the previous paragraph, and we have $G$-morphisms: \[ X\stackrel{\Phi_{\Delta^\prime}}{\longrightarrow} X_1\stackrel{\Psi}{\longrightarrow}G/P\] where $\Psi\circ\Phi_{\Delta^\prime}=\Phi_{\Delta^{\prime\prime}}$ and $P$ is the parabolic subgroup containing $B^-$ and associated to the subset of simple roots $S^p/\Delta^{\prime\prime}$.

Let $H_1\subset G$ be a subgroup such that $X_1$ is the wonderful embedding of $G/H_1$: we may suppose that $H\subset H_1\subset P$.

$P$ is a minimal parabolic subgroup containing $H$ thanks to the properties of $\Delta^{\prime\prime}$. We have also the inclusion $H^u\subset P^u$ between the unipotent radicals: indeed, the morphism: \[H\hookrightarrow P\longrightarrow P/P^u\] yields the inclusion: \[H/(H\cap P^u)\hookrightarrow P/P^u\] $P/P^u$ is reductive and its subgroup $H/(H\cap P^u)$ is not included in any proper parabolic subgroup, hence it is reductive and this implies $H^u\subset P^u$. Therefore we may suppose to have two Levi decompositions $H=L\cdot H^u$ and $P=L(P)\cdot P^u$, where $L\subset L(P)$ and $H^u\subset P^u$; this also implies that $L$ is spherical in $L(P)$ and not contained in any proper parabolic subgroup. Moreover, we may suppose to have a Levi decomposition $H_1=L_1\cdot H_1^u$, where $L\subset L_1$.

The commutator subgroup $(L,L)$\index{$(H,H)$, commutator subgroup of $H$} is a reductive spherical subgroup of $(L_1,L_1)$. Let $Q$ be the minimal parabolic subgroup containing $B^-$ such that $Q^r\subset H_1\subset Q$: in such a case we have $Q^u=H_1^u$. In many cases $Q$ is equal to $P$.

{\em Second step.} We show that $(L,L)$ is forced to be equal to $(L_1,L_1)$, so that we also have $H^u\subset H_1^u$. This is proven checking all possible candidates for $(L,L)$ provided by the list of reductive spherical subgroups in Table~\ref{T4}. For each such candidate, we consider $H^\prime = P^r\cdot(L,L)$, which will be always wonderful. The inclusion $H\subset H^\prime$ would induce a $G$-morphism $\Phi$ with connected fibers from $X$ onto the wonderful embedding of $G/H^\prime$, whose spherical system is known thanks to the previous paragraph. But now the associated distinguished subset of colours $\Delta_\Phi$ will not exist, unless we had chosen $(L,L)=(L_1,L_1)$.

Therefore we are in the following situation: \[ \begin{array}{ccccccc} H   &=&    C   &\cdot& (L,L)     &\cdot& H^u \\ & &  \cap  & &      \|        &     & \cap \\ H_1 &=&    C_1 &\cdot& (L_1,L_1) &\cdot& H_1^u\\ \end{array}\] where $C$ and $C_1$ are the connected centers of $H$ and $H_1$, respectively. We have that the Lie algebra of $H^u$ is an $(L,L)$-submodule of the Lie algebra of $H^u_1$. Moreover, since we know $C_1$ and $\dim C$, we can compute the codimension of $H^u$ in $H^u_1$ using Lemma~\ref{L4}.

{\em Third step.} We analyze the decomposition of the Lie algebra of $H^u_1$ into simple $L_1$-modules: ad-hoc considerations about the spherical system arise conditions on $C$ and the Lie algebra of $H^u$ in such a way that $H$ is uniquely determined up to conjugation.

For the existence part of the proof, namely once we have a candidate $H$, we prove it is spherical and then we use the same kind of indirect arguments of those in the previous section. To prove that $H$ is spherical we use the following:

\begin{lemma}[\cite{B87}]\label{L6} Let $Q$ and $L(Q)$ be the above defined subgroups of $G$, let $\mathfrak n_1$ and $\mathfrak n$ denote the Lie algebras of $H^u_1$ and $H^u$, respectively. If there exists a Borel subgroup $B^\prime$ of $L(Q)$ such that $B^\prime\,L$ is open in $L(Q)$ and $B^\prime\cap L$ has an open orbit in $\mathfrak n_1/\mathfrak n$, then $H$ is spherical.\end{lemma}

Moreover, notice that the proof of the existence for a primitive spherical system implies the existence for everyone of its localizations. Therefore, for any localized case it is necessary to provide only the uniqueness part of the proof.

\subsubsection*{Strongly $\Delta$-connected cases}

In the following strongly $\Delta$-connected cases we get always that $P=Q$.


$\mathbf{dc^\ast(n)}$. We have $S^p=\emptyset$, $\Sigma=\{\alpha_1+\alpha_2,\ \ldots,\ \alpha_{n-3}+\alpha_{n-2},\ \alpha_{n-2}+\alpha_{n-1},\ \alpha_{n-2}+\alpha_n\}$ and $\mathbf A=\emptyset$. The set of colours is $\Delta=\{\delta_{\alpha_1},\ldots,\delta_{\alpha_n}\}$.

{\em Uniqueness.} Suppose $n$ odd. The set $
\Delta^{\prime\prime}=\{\delta_{\alpha_1},\delta_{\alpha_2},\ldots,\delta_{\alpha_{n-2}}\}$
is the minimum parabolic subset of colours, and it is associated
to the parabolic subgroup $P_{S^p/\Delta^{\prime\prime}}$ where
$S^p/\Delta^{\prime\prime}=S\setminus\{\alpha_{n-1},\alpha_n\}$.
The Luna diagram of the quotient spherical system
$(S^p,\Sigma,\mathbf A)/\Delta^{\prime\prime}$ is the following.
\[\begin{picture}(10500,3000)\multiput(0,1500)(1800,0){2}{\usebox{\segm}}\put(7200,1500){\usebox{\segm}}\put(3600,1500){\usebox{\susp}}\put(9000,300){\usebox{\bifurc}}\multiput(10200,300)(0,2400){2}{\circle{600}}\end{picture}\]
Consider the distinguished set of colours $\Delta^\prime=\{\delta_{\alpha_1},\ldots,\delta_{\alpha_{2i+1}},\ldots,\delta_{\alpha_{n-2}}\}\subset\Delta^{\prime\prime}$: it is such that $S^p/\Delta^\prime=\{\alpha_1,\ldots,\alpha_{2i+1},\ldots,\alpha_{n-2}\}$ and $\Sigma(\Delta^\prime)=\{\alpha_{n-2}+\alpha_{n-1},\alpha_{n-2}+\alpha_n\}$. We have \[\Xi/\Delta^\prime=\{\xi\colon \langle\rho(\delta),\xi\rangle=0\ \forall\delta\in\Delta^\prime\ \mathrm{and}\ \langle\gamma^\ast,\xi\rangle=0\ \forall\gamma\in\Sigma(\Delta^\prime)\},\] so the spherical roots of the quotient are the indecomposable elements of the semigroup \[\{\sum_{i=1}^{n-3}c_i(\alpha_i+\alpha_{i+1})\colon c_i\geq0\ \forall i\ \mathrm{and}\ c_i-c_{i+1}=0\ \forall i\ \mathrm{odd}\},\] i.\ e.\  $\Sigma/\Delta^\prime=\{\alpha_1+2\alpha_2+\alpha_3,\ldots,\alpha_{2i-1}+2\alpha_{2i}+\alpha_{2i+1},\ldots,\alpha_{n-4}+2\alpha_{n-3}+\alpha_{n-2}\}$. The quotient spherical system $(S^p,\Sigma,\mathbf A)/\Delta^\prime$ is a parabolic induction of $ac(n-2)$.
\[\begin{picture}(15900,3000)\multiput(0,1200)(3600,0){2}{\usebox{\dthree}}\put(7200,1500){\usebox{\susp}}\put(10800,1200){\usebox{\dthree}}\put(14400,300){\usebox{\bifurc}}\multiput(15600,300)(0,2400){2}{\circle{600}}\end{picture}\]
Our $G$ is $SO(2n)$. Let $Q$ denote the parabolic subgroup $P^-_{S^p/\Delta^{\prime\prime}}$, we have $Q=Q^u\,L(Q)$, where $L(Q)\cong GL(n-1)\times GL(1)$. Let $H_1$ denote the wonderful subgroup associated to the spherical system $(S^p,\Sigma,\mathbf A)/\Delta^\prime$. It is equal to $H^u_1\,L_1$, where $H^u_1=Q^u$, $(L_1,L_1)\cong Sp(n-1)$ and $C_1$ is the connected center of $L(Q)$, where $\dim C_1=2=\rank (\Delta\setminus \Delta^\prime) - \rank \Sigma/\Delta^\prime$. Indeed, the wonderful subgroup of $SL(n-1)$ associated to the spherical system $ac(n-2)$ is isomorphic to $Sp(n-1)\cdot C_{SL(n-1)}$. So $(L_1,L_1)$ is the subgroup of the matrices
\begin{scriptsize}
\[\left(\begin{array}{c|c|c|c}
A & 0 & 0 & 0 \\
\hline
0 & 1 & 0 & 0 \\
\hline
0 & 0 & 1 & 0 \\
\hline
0 & 0 & 0 & J\,^{t}\!A^{-1}J
\end{array}\right)\]
\end{scriptsize}where $A\in Sp(n-1)$, and $C_1$ is the subgroup of the matrices
\begin{scriptsize}
\[\left(\begin{array}{c|c|c|c}
c_1\,I & 0 & 0 & 0 \\
\hline
0 & c_2 & 0 & 0 \\
\hline
0 & 0 & c_2^{-1} & 0 \\
\hline
0 & 0 & 0 & c_1^{-1}\,I
\end{array}\right)\]
\end{scriptsize}where $c_1,c_2\in GL(1)$. The Lie algebra $\mathfrak n_1\subset\mathfrak {so}(2n)$ of $H^u_1$ is the subalgebra of the matrices skew symmetric with respect to the skew diagonal
\begin{scriptsize}
\[\left(\begin{array}{c|c|c|c}
0 & 0 & 0 & 0 \\
\hline
M_1 & 0 & 0 & 0 \\
\hline
M_2 & 0 & 0 & 0 \\
\hline
M_3 & -J\,^t\!M_2 & -J\,^t\!M_1 & 0
\end{array}\right)\]
\end{scriptsize}where $M_1$ and $M_2$ are row $(n-1)$-vectors and $M_3$ is a square $(n-1)$-matrix skew symmetric with respect to the skew diagonal.\\ Let $H=H^u\,L$ be a wonderful subgroup with spherical system $dc^\ast(n)$. We can suppose $H\subset H_1$ and $L\subset L_1$. By Lemma~\ref{L4} we have that $\dim H_1-\dim H=n$ and $\dim C=1$. The subgroup $(L,L)$ of $(L_1,L_1)$ is spherical; it is also reductive because it is contained in no proper parabolic subgroup (thanks to the minimality of $\Delta^{\prime\prime}$). By Table~\ref{T4}, $(L,L)$ can be isomorphic only to $Sp(n-1)$ and hence it is equal to $(L_1,L_1)$. In particular, $\dim H^u_1 - \dim H^u = n-1$.\\ Call $\mathfrak m_i$ the subspace of $\mathfrak n_1$ given by the condition $M_j=0$ for all $j\neq i$. As a $(L,L)$-module, $\mathfrak n_1$ has three isotypic components: the first component is $\mathfrak m_1+\mathfrak m_2$ where $\mathfrak m_1$ and $\mathfrak m_2$ are simple $L_1$-submodules of dimension $n-1$ (two copies of the dual of the standard representation $\mathbb C^{n-1}$), the second component is a simple $L_1$-submodule of dimension $((n-1)(n-2)/2)-1$ and the third component is a simple $L_1$-submodule of dimension $1$. The sum of these two last isotypic components is $\mathfrak m_3\cong(\bigwedge^2\mathbb C^{n-1})^*$. The Lie algebra $\mathfrak n \subset \mathfrak n _1$ of $H^u$ is an $L$-submodule of codimension $n-1$. Being a Lie subalgebra, it must be the sum of $\mathfrak m_3$ and a submodule $\mathfrak m\subset\mathfrak m_1+\mathfrak m_2$ of codimension $n-1$. The subgroup $C_1$ acts as follows, for all $c\in C_1$ and for all $M_i\in\mathfrak m_i$: $c.M_1=c_2\,c_1^{-1}\,M_1$, $c.M_2=c_2^{-1}\,c_1^{-1}\,M_2$. Moreover, $H$ must be wonderful, hence of finite index in its normalizer. Therefore, $C=\{c_2\,c_1^{-1}=c_2^{-1}\,c_1^{-1}\}^\circ=\{c_2=1\}$ and $\mathfrak m$ is diagonal in $\mathfrak m_1+\mathfrak m_2$, namely defined by an equation $\mu_1 M_1+\mu_2 M_2=0$ with $\mu_i\neq 0$. This determines $C$ and $H^u$ up to conjugation. Since $H$ contains the center $C_G$ of the group $G$, the Levi subgroup $L$ is isomorphic to $C_G\cdot C\cdot Sp(n-1)$. Moreover $H=N_G(H)$, hence it is uniquely determined up to conjugation.

Now suppose $n$ even. The subset
$\Delta^{\prime\prime}=\{\delta_{\alpha_2},\delta_{\alpha_3},\ldots,\delta_{\alpha_{n-2}}\}$
is a minimum parabolic subset and
$S^p/\Delta^{\prime\prime}=S\setminus\{\alpha_1,\alpha_{n-1},\alpha_n\}$.
The set
$\Delta^\prime=\{\delta_{\alpha_2},\ldots,\delta_{\alpha_{2i}},\ldots,\delta_{\alpha_{n-2}}\}$
is distinguished, with
$S^p/\Delta^\prime=\{\alpha_2,\ldots,\alpha_{2i},\ldots,\alpha_{n-2}\}$
and
$\Sigma(\Delta^\prime)=\{\alpha_1+\alpha_2,\alpha_{n-2}+\alpha_{n-1},\alpha_{n-2}+\alpha_n\}$,
so that $\Sigma/\Delta^\prime=\{\alpha_2+2\alpha_3+\alpha_4,\ldots,$
$\alpha_{2i}+2\alpha_{2i+1}+\alpha_{2i+2},\ldots,$
$\alpha_{n-4}+2\alpha_{n-3}+\alpha_{n-2}\}$. If $n>4$ the quotient
spherical system is a parabolic induction of $ac(n-3)$, if $n=4$
it has rank zero.
\[\begin{picture}(18000,3000)\put(2100,0){\begin{picture}(15900,3000)\multiput(0,1200)(3600,0){2}{\usebox{\dthree}}\put(7200,1500){\usebox{\susp}}\put(10800,1200){\usebox{\dthree}}\put(14400,300){\usebox{\bifurc}}\multiput(15600,300)(0,2400){2}{\circle{600}}\end{picture}}\put(300,1500){\usebox{\segm}}\put(300,1500){\circle{600}}\end{picture}\]
Again, $G=SO(2n)$. Set $Q=P^-_{S^p/\Delta^{\prime\prime}}$, we have $L(Q)\cong GL(1)\times GL(n-2)\times GL(1)$. Here $H_1$ is such that $(L_1,L_1)\cong Sp(n-2)$ and $C_1$ is the connected center of $L(Q)$, where $\dim C_1=3$. Let $I$ and $J$ denote the diagonal and skew diagonal unitary $(n-2)\times (n-2)$ matrices; we have that $(L_1,L_1)$ is the subgroup of the matrices
\begin{scriptsize}
\[\left(\begin{array}{c|c|c|c|c|c}
1 & 0 & 0 & 0 & 0 & 0 \\
\hline
0 & A & 0 & 0 & 0 & 0 \\
\hline
0 & 0 & 1 & 0 & 0 & 0 \\
\hline
0 & 0 & 0 & 1 & 0 & 0 \\
\hline
0 & 0 & 0 & 0 & J\,^{t}\!A^{-1}J & 0 \\
\hline
0 & 0 & 0 & 0 & 0 & 1
\end{array}\right)\]
\end{scriptsize}where $A\in Sp(n-2)$. The connected center $C_1$ is the subgroup of the matrices
\begin{scriptsize}
\[\left(\begin{array}{c|c|c|c|c|c}
c_1 & 0 & 0 & 0 & 0 & 0 \\
\hline
0 & c_2\,I & 0 & 0 & 0 & 0 \\
\hline
0 & 0 & c_3 & 0 & 0 & 0 \\
\hline
0 & 0 & 0 & c_3^{-1} & 0 & 0 \\
\hline
0 & 0 & 0 & 0 & c_2^{-1}\,I & 0 \\
\hline
0 & 0 & 0 & 0 & 0 & c_1^{-1}
\end{array}\right)\]
\end{scriptsize}where $c_1,c_2,c_3\in GL(1)$. The Lie algebra $\mathfrak n_1\subset\mathfrak {so}(2n)$ of $H^u_1$ is the subalgebra of the matrices skew symmetric with respect to the skew diagonal
\begin{scriptsize}
\[\left(\begin{array}{c|c|c|c|c|c}
0 & 0 & 0 & 0 & 0 & 0 \\
\hline
M_1 & 0 & 0 & 0 & 0 & 0 \\
\hline
M_2 & M_3 & 0 & 0 & 0 & 0 \\
\hline
M_4 & M_5 & 0 & 0 & 0 & 0 \\
\hline
M_6 & M_7 & -J\,^t\!M_5 & -J\,^t\!M_3 & 0 & 0 \\
\hline
0 & -^t\!M_6J & -M_4 & -M_2 & -^t\!M_1J & 0
\end{array}\right)\]
\end{scriptsize}where $M_7$ is a square $(n-2)$-matrix skew symmetric with respect to the skew diagonal.\\ Let $H$ be a wonderful subgroup with the given spherical system, $H=H^u\,L$, $H\subset H_1$ and $L\subset L_1$. Analogously to the odd case we have that $\dim H_1-\dim H=\rank \Sigma - \rank \Sigma/\Delta^\prime+\rank S^p/\Delta^\prime=n$. Moreover, $\dim C=1$. The subgroup $(L,L)$ of $(L_1,L_1)$ can be isomorphic only to $Sp(n-2)$ and hence it is equal to $(L_1,L_1)$. In particular, $\dim H^u_1 - \dim H^u = n-2$.\\ The $(L,L)$-module $\mathfrak n_1$ decomposes into three isotypic components: the first one is $\mathfrak m_1+\mathfrak m_3+\mathfrak m_5+\mathfrak m_6$ where these modules are all simple $L_1$-modules of dimension $n-2$ (the standard representation $\mathbb C^{n-2}$ or its dual), the second component is a simple $L_1$-submodule of dimension $((n-2)(n-3)/2)-1$ and the third component is a direct sum of three simple $L_1$-submodules of dimension $1$ (the sum of these two last isotypic components is $\mathfrak m_2+\mathfrak m_4 +\mathfrak m_7$). Notice that if $n=4$ the second component is zero. The Lie subalgebra $\mathfrak n \subset \mathfrak n _1$ is an $L$-submodule of codimension $n-2$. Since it is a Lie subalgebra, it is the sum of $\mathfrak m_2+\mathfrak m_4+\mathfrak m_6+\mathfrak m_7$ and a submodule $\mathfrak m\subset\mathfrak m_1+\mathfrak m_3+\mathfrak m_5$ of codimension $n-2$. The subgroup $C_1$ acts as follows, for all $c\in C_1$ and for all $M_i\in\mathfrak m_i$: $c.M_1=c_2\,c_1^{-1}\,M_1$, $c.M_3=c_3\,c_2^{-1}\,M_3$, $c.M_5=c_3^{-1}\,c_2^{-1}\,M_5$. Therefore, $C=\{c_2\,c_1^{-1}=c_3\,c_2^{-1}=c_3^{-1}\,c_2^{-1}\}^\circ=\{c_3=1,c_1=c_2^2\}$ and $\mathfrak m$ is diagonal in $\mathfrak m_1+\mathfrak m_3+\mathfrak m_5$, namely defined by an equation $\mu_1 M_1+\mu_3 M_3+\mu_5 M_5=0$ with $\mu_i\neq 0$, hence unique up to conjugation. Therefore $L\cong C_G\cdot C\cdot Sp(n-2)$, and Analogously to the odd case $H$ is determined up to conjugation.

{\em Existence} We can suppose $n$ odd, the even case follows by
localization. The subgroup $H$ is spherical. Indeed, consider:
\[\begin{scriptsize} g=
\left(\begin{array}{c|c}A & 0 \\ \hline 0 &
J\,^t\!A^{-1}J\end{array}\right),
\end{scriptsize}\] where
$A\in GL(n)$ is the following:
\[A=\begin{scriptsize}\left(\begin{array}{c|c|c}I_{\frac{n-1}{2}} &
0 & 0 \\ \hline 0 & 1 & 0 \\ \hline J_{\frac{n-1}{2}}& 0 &
I_{\frac{n-1}{2}} \end{array}\right),\end{scriptsize}\] then
$B\,gHg^{-1}$ is open in $G$, since $\mathfrak b+g\mathfrak h
g^{-1}=\mathfrak{so}(2n)$. Moreover, $H$ is wonderful, because $H=N_G(H)$. We have $\dim P^u_{S^p}+\rank
\Sigma=\dim G/H=(n+1)(n-1)$ and $\rank\Delta-\rank\Sigma=\dim
C=1$. Since $\rank\Sigma\leq n$ we have $\rank S^p\leq1$, but the case
of $\rank S^p=1$ gives $\rank\Sigma=n$ and hence $S^p=\emptyset$,
i.e.\ a contradiction. Therefore, $S^p=\emptyset$ and
$\rank\Sigma=n-1$. The spherical system must be cuspidal, hence
$\Sigma$ contains a root of type $a_2$ and it remains only the
case of $dc^\ast(n)$.

$\mathbf{dy(p,p)}$ ($p\geq4$). {\em Uniqueness.} We have $\Sigma=S$. The set $\Delta^{\prime\prime}=\Delta\setminus\{\delta_{\alpha_{p-1}}^-$, $\delta_{\alpha_p}^-$, $\delta_{\alpha^\prime_{p-1}}^-$, $\delta_{\alpha^\prime_p}^-\}$ is the minimum parabolic subset of $\Delta$, associated to $P_{S^p/\Delta^{\prime\prime}}$, where $S^p/\Delta^{\prime\prime}=\{\alpha_1,\ldots,\alpha_{p-2},\ \alpha^\prime_1,\ldots,\alpha^\prime_{p-2}\}$. We take $\Delta^\prime=\{\delta_{\alpha_1}^+,\ \delta_{\alpha_2}^+,\ \ldots,\ \delta_{\alpha_{p-2}}^+,\ \delta_{\alpha_{p-1}}^+\}$ $\subset\Delta^{\prime\prime}$. Now $S^p/\Delta^\prime=\emptyset$ and $\Sigma(\Delta^\prime)=\{\alpha_{p-1},\ \alpha_p,\ \alpha^\prime_{p-1},\ \alpha^\prime_p\}$, and so we have that $\Sigma/\Delta^\prime=\{\alpha_1+\alpha^\prime_{p-2},\ \ldots,\ \alpha_i+\alpha^\prime_{p-i-1},\ \ldots,\ \alpha_{p-2}+\alpha^\prime_1\}$. The quotient is a parabolic induction of $aa(p-2,p-2)$:
\[\begin{picture}(15300,4200)\multiput(0,0)(8100,0){2}{\put(300,2700){\usebox{\segm}}\put(2100,2700){\usebox{\susp}}\put(5700,1500){\usebox{\bifurc}}\multiput(300,2700)(1800,0){2}{\circle{600}}\put(5700,2700){\circle{600}}\multiput(6900,1500)(0,2400){2}{\circle{600}}}\multiput(300,2400)(13500,0){2}{\line(0,-1){2400}}\put(300,0){\line(1,0){13500}}\put(2100,2400){\line(0,-1){2100}}\put(2100,300){\line(1,0){9900}}\multiput(12000,300)(0,300){4}{\line(0,1){150}}\put(10200,2400){\line(0,-1){1800}}\put(10200,600){\line(-1,0){6300}}\multiput(3900,600)(0,300){4}{\line(0,1){150}}\multiput(5700,2400)(2700,0){2}{\line(0,-1){1500}}\put(5700,900){\line(1,0){2700}}\end{picture}\]
Our $G$ is $SO(2p)\times SO(2p)$. Consider $Q=P^-_{S^p/\Delta^{\prime\prime}}$: we have $Q=Q^u\,L(Q)$ where $L(Q)=GL(p-1)\times GL(1)\times GL(p-1)\times GL(1)$. Also, $H_1=H^u_1\,L_1$, where $H^u_1=Q^u$, $(L_1,L_1)=SL(p-1)$ embedded diagonally in $GL(p-1)\times 1\times GL(p-1)\times 1\subset L(Q)$ and $C_1$ is the connected center of $L(Q)$, and $\dim C_1=4$. The elements of $(L_1,L_1)$ have this form:
\begin{scriptsize}
\[\big(\left(\begin{array}{c|c|c|c}
A & 0 & 0 & 0 \\
\hline
0 & 1 & 0 & 0 \\
\hline
0 & 0 & 1 & 0 \\
\hline
0 & 0 & 0 & \Big(J\,^{t}\!A^{-1}J\Big)
\end{array}\right),\left(\begin{array}{c|c|c|c}
A & 0 & 0 & 0 \\
\hline
0 & 1 & 0 & 0 \\
\hline
0 & 0 & 1 & 0 \\
\hline
0 & 0 & 0 & \Big(J\,^{t}\!A^{-1}J\Big)
\end{array}\right)\big)\]
\end{scriptsize}where $A\in SL(p-1)$, while $C_1$ is the subgroup of the couples
\begin{scriptsize}
\[\big(\left(\begin{array}{c|c|c|c}
c_1\,I & 0 & 0 & 0 \\
\hline
0 & c_2 & 0 & 0 \\
\hline
0 & 0 & c_2^{-1} & 0 \\
\hline
0 & 0 & 0 & c_1^{-1}\,I
\end{array}\right),\left(\begin{array}{c|c|c|c}
c_3\,I & 0 & 0 & 0 \\
\hline
0 & c_4 & 0 & 0 \\
\hline
0 & 0 & c_4^{-1} & 0 \\
\hline
0 & 0 & 0 & c_3^{-1}\,I
\end{array}\right)\big)\]
\end{scriptsize}where $c_1,c_2,c_3,c_4\in GL(1)$. The Lie algebra $\mathfrak n_1$ of $H^u_1$ is the subalgebra of the couples of skew symmetric matrices with respect to the skew diagonal
\begin{scriptsize}
\[\big(\left(\begin{array}{c|c|c|c}
0 & 0 & 0 & 0 \\
\hline
M_1 & 0 & 0 & 0 \\
\hline
M_2 & 0 & 0 & 0 \\
\hline
M_3 & -J\,^t\!M_2 & -J\,^t\!M_1 & 0
\end{array}\right),\left(\begin{array}{c|c|c|c}
0 & 0 & 0 & 0 \\
\hline
M_4 & 0 & 0 & 0 \\
\hline
M_5 & 0 & 0 & 0 \\
\hline
M_6 & -J\,^t\!M_5 & -J\,^t\!M_4 & 0
\end{array}\right)\big)\]
\end{scriptsize}where $M_1,M_2,M_4,M_5$ are row $(p-1)$-vectors and $M_3,M_6$ are square $(p-1)$-matrices skew symmetric with respect to the skew diagonal.\\ As usual we take $H=H^u\,L$ with $H\subset H_1$ and $L\subset L_1$. We have that $\dim H_1-\dim H=p+2$. Moreover, $\dim C=1$, hence $\dim (L_1,L_1)-\dim (L,L)\leq p-1$. Table~\ref{T4} shows that $(L,L)$ can be isomorphic only to $SL(p-1)$, $SO(p-1)$ or $Sp(p-1)$ (if $p-1$ is even), but from the above inequality we get that $(L,L)=(L_1,L_1)\cong SL(p-1)$. In particular, $\dim H^u_1 - \dim H^u = p-1$.\\ The $(L,L)$-module $\mathfrak n_1$ decomposes into two isotypic components: the first one is $\mathfrak m_1+\mathfrak m_2+\mathfrak m_4+\mathfrak m_5$ where each of these is a simple $L_1$-submodule of dimension $p-1$ and of highest weight $\lambda_2$ (the dual of the standard representation $\mathbb C^{p-1}$), the second component is $\mathfrak m_3+\mathfrak m_6$ where these are simple $L_1$-submodules of dimension $(p-1)(p-2)/2$ and of highest weight $\lambda_1$ (the dual of $\bigwedge^2\mathbb C^{p-1}$). The Lie algebra $\mathfrak n \subset \mathfrak n _1$ of $H^u$ is an $L$-submodule of codimension $p-1$. Since it is a Lie subalgebra, it is the sum of $\mathfrak m_3+\mathfrak m_6$ and a submodule $\mathfrak m\subset\mathfrak m_1+\mathfrak m_2+\mathfrak m_4+\mathfrak m_5$ of codimension $p-1$. The subgroup $C_1$ acts as follows, for all $c\in C_1$ and for all $M_i\in\mathfrak m_i$: $c.M_1=c_2\,c_1^{-1}\,M_1$, $c.M_2=c_2^{-1}\,c_1^{-1}\,M_2$, $c.M_4=c_4\,c_3^{-1}\,M_4$, $c.M_5=c_4^{-1}\,c_3^{-1}\,M_5$. Therefore, $C=\{c_2\,c_1^{-1}=c_2^{-1}\,c_1^{-1}=c_4\,c_3^{-1}=c_4^{-1}\,c_3^{-1}\}^\circ=\{c_1=c_3,\ c_2=c_4=1\}$ and $\mathfrak m$ is diagonal in $\mathfrak m_1+\mathfrak m_2+\mathfrak m_4+\mathfrak m_5$, namely defined by an equation $\mu_1 M_1+\mu_2 M_2+\mu_4 M_4+\mu_5 M_5=0$ with $\mu_i\neq 0$, hence unique up to conjugation. Since $H$ contains the center $C_G$ of the group $G$, we have $L=C_G\cdot GL(p-1)$, where $GL(p+1)$ is embedded diagonally in $L(Q)$. Moreover, $H=N_G(H)$ so it is uniquely determined up to conjugation.

{\em Existence.} $L^\circ$ is isomorphic to $GL(p-1)$. Now $B^\prime=B_{GL(p-1)}\times GL(1)\times B^-_{GL(p-1)} \times GL(1)$ is a Borel subgroup of $L(Q)$ that satisfies the hypothesis of the Lemma~\ref{L6}. Indeed, $B^\prime\cap L^\circ=T_{GL(p-1)}$ has an open orbit in $\mathbb C^{p-1}$. Hence, $H$ is wonderful, since it is spherical and equal to its normalizer. By Lemma~\ref{L4}, $\dim G/H=2p^2$ and $\dim C=1$. Since $\dim B^u=2p^2-2p$, we have $S^p=\emptyset$ and $\rank \Sigma=2p=\rank S$. The spherical system is hence cuspidal and $\Sigma$ does not contain roots of type $a_m$ for $m\geq2$, nor roots of type $a_1\times a_1$ or $d_m$. A root of type $a_1^\prime$ would cause the spherical system to be decomposable, since it cannot be nonorthogonal to any root of type $a_1$, but this is not the case, since $H$ is not a direct product. Therefore, $\Sigma=S$ and by Lemma~\ref{L4a}(\textit{\ref{aone}}) the spherical system is strongly $\Delta$-connected. Since $H$ is not solvable, by the Proposition~\ref{morphisms} the spherical system cannot be a comb, hence we have $dy(p,p)$ or $dz(4,4)$ if $p=4$. The latter is excluded by the corresponding uniqueness proof.

$\mathbf{dy(p,(p-1))}$ ($p\geq4$). {\em Uniqueness.} The set $\Delta^{\prime\prime}=\Delta\setminus\{\delta_{\alpha_{p-1}}^-,\delta_{\alpha^\prime_{p-1}}^-,\delta_{\alpha^\prime_p}^-\}$ is the minimum parabolic subset of $\Delta$ and $S^p/\Delta^{\prime\prime}=\{\alpha_1,\ldots,\alpha_{p-2},\ \alpha^\prime_1,\ldots,\alpha^\prime_{p-2}\}$. Set $\Delta^\prime=\{\delta_{\alpha_1}^+,\ \delta_{\alpha_2}^+,\ \ldots,\ \delta_{\alpha_{p-2}}^+,\ \delta_{\alpha_{p-1}}^+\}\subset\Delta^{\prime\prime}$. We have $S^p/\Delta^\prime=\emptyset$ and $\Sigma/\Delta^\prime=\{\alpha_1+\alpha^\prime_{p-2},\ \ldots,\ \alpha_i+\alpha^\prime_{p-i-1},\ \ldots,\ \alpha_{p-2}+\alpha^\prime_1\}$, the quotient spherical system is a parabolic induction of $aa(p-2,p-2)$. Let us pose $G=SL(p)\times SO(2p)$. We have $L(Q)=GL(p-1)\times GL(1)\times GL(p-1)\times GL(1)$. Here $(L_1,L_1)=SL(p-1)$ embedded diagonally in $GL(p-1)\times 1\times GL(p-1)\times 1\subset L(Q)$ and $\dim C_1=3$. We have that $\dim H_1-\dim H=p+1$ and $\dim C=1$, hence $\dim (L_1,L_1)-\dim (L,L)\leq p-1$. As for $dy(p,p)$, the subgroup $(L,L)$ can be isomorphic only to $SL(p-1)$ and hence it is equal to $(L_1,L_1)$. In particular, $\dim H^u_1 - \dim H^u = p-1$.\\ The $(L,L)$-module $\mathfrak n_1$ decomposes into two isotypic components: the first one is $\mathfrak m_1+\mathfrak m_2+\mathfrak m_3$ (three simple $L_1$-submodules of dimension $p-1$), the second component $\mathfrak m_4$ is a simple $L_1$-submodule of dimension $(p-1)(p-2)/2$. The Lie algebra $\mathfrak n \subset \mathfrak n _1$ of $H^u$ is then $\mathfrak m_4+\mathfrak m$ where $\mathfrak m\subset\mathfrak m_1+\mathfrak m_2+\mathfrak m_3$ has codimension $p-1$. As usual, $C$ is determined and $\mathfrak m$ is diagonal in $\mathfrak m_1+\mathfrak m_2+\mathfrak m_3$. And $L$ is determined too, since $H\supset C_G$. We have $H=N_G(H)$, hence determined up to conjugation.


$\mathbf{dy(p,(p-2))}$ ($p\geq4$). {\em Uniqueness.} The set $\Delta^{\prime\prime}=\Delta\setminus\{\delta_{\alpha^\prime_{p-1}}^-,\delta_{\alpha^\prime_p}^-\}$ is the minimum parabolic subset of $\Delta$ and $S^p/\Delta^{\prime\prime}=\{\alpha_1,\ldots,\alpha_{p-2},\ \alpha^\prime_1,\ldots,\alpha^\prime_{p-2}\}$. We take $\Delta^\prime=\{\delta_{\alpha_1}^+,\ \delta_{\alpha_2}^+,\ \ldots,\ \delta_{\alpha_{p-2}}^+,\ \delta_{\alpha^\prime_1}^+\}\subset\Delta^{\prime\prime}$, so $S^p/\Delta^\prime=\emptyset$ and $\Sigma/\Delta^\prime=\{\alpha_1+\alpha^\prime_{p-2},\ \ldots,\ \alpha_i+\alpha^\prime_{p-i-1},\ \ldots,\ \alpha_{p-2}+\alpha^\prime_1\}$. The quotient spherical system is a parabolic induction of $aa(p-2,p-2)$. Let us pose $G=SL(p-1)\times SO(2p)$. We have $L(Q)=SL(p-1)\times GL(p-1)\times GL(1)$. Here $(L_1,L_1)=SL(p-1)$ embedded diagonally in $SL(p-1)\times GL(p-1)\times 1\subset L(Q)$ and $\dim C_1=2$. We have that $\dim H_1-\dim H=p$ and $\dim C=1$, hence $\dim (L_1,L_1)-\dim (L,L)\leq p-1$; so $(L,L)=(L_1,L_1)$ and $\dim H^u_1 - \dim H^u = p-1$. The $(L,L)$-module $\mathfrak n_1$ decomposes into two isotypic components: the first one is $\mathfrak m_1+\mathfrak m_2$ where these are simple $L_1$-submodules of dimension $p-1$, the second one $\mathfrak m_3$ is a simple $L_1$-submodule of dimension $(p-1)(p-2)/2$. The Lie subalgebra $\mathfrak n$ is $\mathfrak m_3+\mathfrak m$ where $\mathfrak m\subset\mathfrak m_1+\mathfrak m_2$ has codimension $p-1$. Hence, $C$ is determined and $\mathfrak m$ is diagonal in $\mathfrak m_1+\mathfrak m_2$. Finally, $H=N_G(H)$ and it is uniquely determined up to conjugation.


$\mathbf{dy(7)}$. {\em Uniqueness.} Take $\Delta^{\prime\prime}=\{\delta_{\alpha_1}^+,\delta_{\alpha_1}^-,\delta_{\alpha_2}^+,\delta_{\alpha_2}^-,\delta_{\alpha_3}^-\}$: we have $S^p/\Delta^{\prime\prime}=\{\alpha_1,\alpha_2,\alpha_5,\alpha_6\}$. If $\Delta^\prime=\{\delta_{\alpha_1}^-,\delta_{\alpha_2}^-,\delta_{\alpha_3}^-\}$, then $S^p/\Delta^\prime=\emptyset$ and $\Sigma/\Delta^\prime=\{\alpha_1+\alpha_5,\alpha_2+\alpha_6\}$. The quotient spherical system is a parabolic induction of $aa(2,2)$. Our $G$ is $SO(14)$, so $L(Q)=GL(3)\times GL(1)\times GL(3)$ and $(L_1,L_1)=SL(3)$ embedded skew diagonally in $GL(3)\times 1\times GL(3)\subset L(Q)$, with $\dim C_1=3$. We have $\dim H_1-\dim H=5$ and $\dim C=1$, hence $\dim (L_1,L_1)-\dim (L,L)\leq 3$. $(L,L)$ can be isomorphic only to $SL(3)$ and hence it is equal to $(L_1,L_1)$. In particular, $\dim H^u_1 - \dim H^u = 3$. The Lie subalgebra $\mathfrak n \subset \mathfrak n _1$ is the sum of a codimension $3$ submodule $\mathfrak m$ diagonal in $\mathfrak m_1+\mathfrak m_2+\mathfrak m_3$ and the rest of the simple $L_1$-submodules of $\mathfrak n_1$; $C$ is determined, $H=N_G(H)$ and it is uniquely determined up to conjugation.

{\em Existence.} By Lemma~\ref{L6}, the subgroup $H$ is spherical and wonderful. By Lemma~\ref{L4} we have $S^p=\emptyset$ and $\rank \Sigma=7$. The spherical system is cuspidal with only roots of type $a_1$. Therefore, $\Sigma=S$, the spherical system is strongly $\Delta$-connected and it is not a comb, so it is $de_6(7)$, $dy(7)$, $de_7(7)$ or $de_8(7)$. We can exclude the other cases by the corresponding uniqueness proofs.

$\mathbf{dz(4,4)}$. {\em Uniqueness.} Take $\Delta^{\prime\prime}=\{\delta_{\alpha_1}^+,\delta_{\alpha_2}^+,\delta_{\alpha_2}^-\}$, $\Delta^\prime=\{\delta_{\alpha_1}^+,\delta_{\alpha_2}^+\}$, so $S^p/\Delta^{\prime\prime}=\{\alpha_2,\alpha^\prime_2\}$ and $S^p/\Delta^\prime=\emptyset$, $\Sigma/\Delta^\prime=\{\alpha_2+\alpha^\prime_2\}$. The quotient by $\Delta^\prime$ is a parabolic induction of $aa(1,1)$. $G$ is $SO(8)\times SO(8)$, so $L(Q)=GL(1)\times GL(2)\times GL(1)\times GL(1)\times GL(2)\times GL(1)$. $(L_1,L_1)=SL(2)$ embedded diagonally in $1\times GL(2)\times 1\times 1\times GL(2)\times 1\subset L(Q)$ and $\dim C_1=6$. Here $\dim H_1-\dim H=7$ and $\dim C=1$, hence $\dim (L_1,L_1)-\dim (L,L)\leq 2$. Thus $(L,L)$ can be isomorphic only to $SL(2)$ and hence it is equal to $(L_1,L_1)$. In particular, $\dim H^u_1 - \dim H^u = 2$. Here $\mathfrak n_1$ has two isotypic components: the first one is the sum of eight simple $L_1$-submodules of dimension $2$ and the second component is composed by six summands of dimension $1$. $\mathfrak n \subset \mathfrak n _1$ is the sum of $\mathfrak m$ diagonal in $\mathfrak m_1+\mathfrak m_2+\mathfrak m_3+\mathfrak m_4+\mathfrak m_5+\mathfrak m_6$ (Table~\ref{T3}) having codimension $2$, and the rest of the simple $L_1$-submodules. $C$ is analogously determined. $H=N_G(H)$ and uniquely determined up to conjugation.

{\em Existence.} By Lemma~\ref{L6}, $H$ is spherical and wonderful. By Lemma~\ref{L4}, $\dim P^u_{S^p}+\rank \Sigma=32$ and $\rank \Delta-\rank \Sigma=1$, thus $S^p$ is empty and $\rank \Sigma=8$. The spherical system is hence cuspidal, without roots of type $a_m$ for $m\geq2$, nor $a_1\times a_1$ nor $d_m$. Being irreducible, there are also no roots of type $a_1^\prime$. So $\Sigma=S$ and the spherical system is strongly $\Delta$-connected. It cannot be a comb, thus we have $dy(4,4)$ or $dz(4,4)$. The former is excluded by the corresponding uniqueness proof.

$\mathbf{dz(4,3)}$.  {\em Uniqueness.} Take $\Delta^{\prime\prime}=\{\delta_{\alpha_1}^+,\delta_{\alpha_2}^+,\delta_{\alpha_2}^-\}$ and $\Delta^\prime=\{\delta_{\alpha_1}^+,\delta_{\alpha_2}^+\}$, so $S^p/\Delta^\prime=\emptyset$ and $\Sigma/\Delta^\prime=\{\alpha_2+\alpha^\prime_2\}$. The quotient by $\Delta^\prime$ is a parabolic induction of $aa(1,1)$. $G$ is $SL(4)\times SO(8)$ and $(L_1,L_1)=SL(2)$ embedded diagonally in $L(Q)$. We have $\dim C_1=5$, $\dim H_1-\dim H=6$ and $\dim C=1$, hence $\dim (L_1,L_1)-\dim (L,L)\leq 2$. As usual, $(L,L)=(L_1,L_1)$ so $\dim H^u_1 - \dim H^u = 2$. The Lie subalgebra $\mathfrak n \subset \mathfrak n _1$ is the sum of a codimension $2$ submodule $\mathfrak m$ diagonal in the direct sum of five simple $L_1$-submodules of dimension $2$, plus the rest of the simple $L_1$-submodules. $C$ is determined, $H=N_G(H)$ and it is uniquely determined up to conjugation.


$\mathbf{dz(4,2)}$. {\em Uniqueness.} Take $\Delta^{\prime\prime}=\{\delta_{\alpha_1}^+,\delta_{\alpha_2}^+,\delta_{\alpha_2}^-\}$, $\Delta^\prime=\{\delta_{\alpha_1}^+,\delta_{\alpha_2}^+\}$. The quotient by $\Delta^\prime$ is a parabolic induction of $aa(1,1)$. Our $G$ is $SL(3)\times SO(8)$ and $(L_1,L_1)\cong SL(2)$; also $\dim C_1=4$ and $\dim (L_1,L_1)-\dim (L,L)\leq 2$. Again, $(L,L)=(L_1,L_1)$ hence $\dim H^u_1 - \dim H^u = 2$. The Lie subalgebra $\mathfrak n \subset \mathfrak n _1$ is the sum of a codimension $2$ submodule $\mathfrak m$ diagonal in the direct sum of four simple $L_1$-submodules of dimension $2$ and the rest of the simple $L_1$-submodules. $C$ is determined, and $H=N_G(H)$ so it is uniquely determined up to conjugation.


$\mathbf{dz(4,1)}$. {\em Uniqueness.} By renaming the simple roots as usually, $\Delta^{\prime\prime}=\{\delta_{\alpha_1}^+$, $\delta_{\alpha_1}^-$, $\delta_{\alpha^\prime_2}^+\}$ and $\Delta^\prime=\{\delta_{\alpha_1}^+,\delta_{\alpha^\prime_2}^+\}$. The Lie subalgebra $\mathfrak n \subset \mathfrak n _1$ is the sum of a codimension $2$ submodule $\mathfrak m$ diagonal in the direct sum of three simple $L_1$-submodules of dimension $2$, plus the rest of the simple $L_1$-submodules. $C$ is determined and $H=N_G(H)$, so it is uniquely determined up to conjugation.


$\mathbf{de_6(7)}$. {\em Uniqueness.} Take $\Delta^{\prime\prime}=\{\delta_{\alpha_1}^+,\delta_{\alpha_2}^+,\delta_{\alpha_2}^-\}$ and $\Delta^\prime=\{\delta_{\alpha_1}^+,\delta_{\alpha_2}^+\}$, so $S^p/\Delta^{\prime\prime}=\{\alpha_2,\alpha_5\}$, $S^p/\Delta^\prime=\emptyset$ and $\Sigma/\Delta^\prime=\{\alpha_2+\alpha_5\}$. The quotient by $\Delta^\prime$ is a parabolic induction of $aa(1,1)$. Let $G$ be $SO(14)$. We have $L(Q)=GL(1)\times GL(2)\times GL(1)\times GL(2)\times GL(1)$, with $(L_1,L_1)=SL(2)$ embedded diagonally in $1\times GL(2)\times 1\times GL(2)\times 1\subset L(Q)$ and $\dim C_1=5$. We have $\dim H_1-\dim H=6$ and $\dim C=1$, hence $\dim (L_1,L_1)-\dim (L,L)\leq 2$. $(L,L)$ can be isomorphic only to $SL(2)$ and hence it is equal to $(L_1,L_1)$. In particular, $\dim H^u_1 - \dim H^u = 2$. The Lie subalgebra $\mathfrak n \subset \mathfrak n _1$ is the sum of a codimension $2$ submodule $\mathfrak m$ diagonal in $\mathfrak m_1+\mathfrak m_2+\mathfrak m_3+\mathfrak m_4+\mathfrak m_5$, plus the rest of the simple $L_1$-submodules. $C$ is determined, and $H=N_G(H)$ so it is uniquely determined up to conjugation.

{\em Existence.} By Lemmas~\ref{L6} and ~\ref{L4}, $H$ is spherical and wonderful, with $S^p=\emptyset$ and $\rank \Sigma=7$. The spherical system is cuspidal and with only roots of type $a_1$. Therefore, $\Sigma=S$ and the spherical system is strongly $\Delta$-connected. It cannot be a comb, and we have $de_6(7)$, $dy(7)$, $de_7(7)$ or $de_8(7)$. We can exclude the other cases by the corresponding uniqueness proofs.

$\mathbf{de_6(6)}$. {\em Uniqueness.} By renaming the simple roots, $\Delta^{\prime\prime}=\{\delta_{\alpha_1}^+,\delta_{\alpha_1}^-,\delta_{\alpha_2}^+\}$, $\Delta^\prime=\{\delta_{\alpha_1}^+,\delta_{\alpha_2}^+\}$. The quotient by $\Delta^\prime$ is a parabolic induction of $aa(1,1)$. Here $G=SO(12)$, $(L_1,L_1)\cong SL(2)$; also $\dim C_1=4$, $\dim H_1-\dim H=5$ and $\dim C=1$. Since $(L,L)=(L_1,L_1)$, $\dim H^u_1 - \dim H^u = 2$. Here $\mathfrak n \subset \mathfrak n _1$ and $C$ are determined as before, $H=N_G(H)$ and it is uniquely determined up to conjugation.


$\mathbf{de_7(8)}$. {\em Uniqueness.} Take $\Delta^{\prime\prime}=\{\delta_{\alpha_1}^+,\delta_{\alpha_2}^+,\delta_{\alpha_2}^-,\delta_{\alpha_3}^+,\delta_{\alpha_3}^-\}$ and $\Delta^\prime=\{\delta_{\alpha_1}^+,\delta_{\alpha_2}^+,\delta_{\alpha_3}^+\}$, so that $S^p/\Delta^{\prime\prime}=\{\alpha_2,\alpha_3,\alpha_5,\alpha_6\}$ and the quotient by $\Delta^\prime$ is a parabolic induction of $aa(2,2)$. Our $G$ is $SO(16)$, so $L(Q)=GL(1)\times GL(3)\times GL(3)\times GL(1)$ and $(L_1,L_1)=SL(3)$ embedded skew diagonally in $1\times GL(3)\times GL(3)\times 1\subset L(Q)$. We have $\dim C_1=4$, $\dim H_1-\dim H=6$ and $\dim C=1$; since $\dim (L_1,L_1)-\dim (L,L)\leq 3$ then $(L,L)=(L_1,L_1)$, therefore $\dim H^u_1 - \dim H^u = 3$. $\mathfrak n \subset \mathfrak n _1$ is the sum of a codimension $3$ submodule $\mathfrak m\subset\mathfrak m_1+\mathfrak m_2+\mathfrak m_3+\mathfrak m_4$ and the rest of the simple $L_1$-submodules. Also, $\mathfrak m_2=\mathfrak m_2^\prime+\mathfrak m_2^{\prime\prime}$: now $\mathfrak m_2^\prime\cong\mathfrak m_1\cong\mathfrak m_3\cong\mathfrak m_4$ are all isomorphic to the dual of $\bigwedge^2\mathbb C^3$, while $\mathfrak m_2^{\prime\prime}$ is isomorphic to the dual of $\Sym^2 \mathbb C^3$. So $\mathfrak m$ is the sum of $\mathfrak m_2^{\prime\prime}$ and a submodule diagonal in $\mathfrak m_1+\mathfrak m_2^\prime+\mathfrak m_3+\mathfrak m_4$. $C$ is determined and $H=N_G(H)$, so it is uniquely determined up to conjugation.

{\em Existence.} By Lemmas~\ref{L6} and ~\ref{L4}, $H$ is spherical and wonderful with $\Sigma=S$, the spherical system is strongly $\Delta$-connected and it is not a comb. We have $de_7(8)$ or $de_8(8)$. We can exclude the latter by the corresponding uniqueness proof.

$\mathbf{de_7(7)}$. {\em Uniqueness.} Take $\Delta^{\prime\prime}=\{\delta_{\alpha_1}^+,\delta_{\alpha_1}^-,\delta_{\alpha_2}^+,\delta_{\alpha_2}^-,\delta_{\alpha_4}^+\}$, $\Delta^\prime=\{\delta_{\alpha_1}^+,\delta_{\alpha_2}^+,\delta_{\alpha_4}^+\}$. The quotient by $\Delta^\prime$ is a parabolic induction of $aa(2,2)$. Let $G$ be $SO(14)$: $(L_1,L_1)\cong SL(3)$ and $\dim C_1=3$. We have $\dim H_1-\dim H=3$ and $\dim C=1$, $(L,L)=(L_1,L_1)$ and $\dim H^u_1 - \dim H^u = 3$. The Lie subalgebra $\mathfrak n \subset \mathfrak n _1$ and $C$ are determined as usual; $H=N_G(H)$ and it is uniquely determined up to conjugation.


$\mathbf{de_8(8)}$. {\em Uniqueness.} Take $\Delta^{\prime\prime}=\{\delta_{\alpha_1}^-,\delta_{\alpha_2}^+,\delta_{\alpha_2}^-,\delta_{\alpha_3}^+,\delta_{\alpha_3}^-,\delta_{\alpha_4}^+,\delta_{\alpha_4}^-\}$ and $\Delta^\prime=\{\delta_{\alpha_1}^-$, $\delta_{\alpha_2}^+$, $\delta_{\alpha_3}^+$, $\delta_{\alpha_4}^-\}$, so $S^p/\Delta^{\prime\prime}=\{\alpha_2,\alpha_3,\alpha_4,\alpha_6,\alpha_7,\alpha_8\}$ and the quotient by $\Delta^\prime$ is a parabolic induction of $aa(3,3)$. Here $G=SO(16)$, $L(Q)\cong GL(1)\times GL(4)\times SO(6)$, $(L_1,L_1)\cong SL(4)$ and $\dim C_1=2$. We have $\dim H_1-\dim H=5$ and $\dim C=1$, hence $\dim (L_1,L_1)-\dim (L,L)\leq 4$. $(L,L)=(L_1,L_1)$ thus $\dim H^u_1 - \dim H^u = 4$. The Lie subalgebra $\mathfrak n \subset \mathfrak n _1$ is the sum of a codimension $4$ submodule $\mathfrak m\subset\mathfrak m_1+\mathfrak m_2$ and the rest of the simple $L_1$-submodules. Now $\mathfrak m_2\cong\big(\bigwedge^2\mathbb C^4\big)\otimes\big(\mathbb C^4\big)=\mathfrak m_2^\prime+\mathfrak m_2^{\prime\prime}$ where $\mathfrak m_2^\prime\cong\mathfrak m_1\cong\mathbb C^4$ (the standard representation), while $\mathfrak m_2^{\prime\prime}$ is simple with $\dim \mathfrak m_2^{\prime\prime}=20$. So $\mathfrak m$ is the sum of $\mathfrak m_2^{\prime\prime}$ and a submodule diagonal in $\mathfrak m_1+\mathfrak m_2^\prime$. $C$ is determined, $H=N_G(H)$ and it is uniquely determined up to conjugation.

{\em Existence.} By Lemmas~\ref{L6} and ~\ref{L4}, $H$ is spherical and wonderful with $\Sigma=S$, the spherical system is strongly $\Delta$-connected and it is not a comb. We have $de_7(8)$ or $de_8(8)$. We can exclude the former by the corresponding uniqueness proof.

$\mathbf{de_8(7)}$. {\em Uniqueness.} We take $\Delta^{\prime\prime}=\{\delta_{\alpha_1}^+,\delta_{\alpha_1}^-,\delta_{\alpha_2}^+,\delta_{\alpha_2}^-,\delta_{\alpha_3}^+,\delta_{\alpha_3}^-,\delta_{\alpha_7}^+\}$, $\Delta^\prime=\{\delta_{\alpha_1}^+,\delta_{\alpha_2}^+,\delta_{\alpha_3}^-,\delta_{\alpha_7}^+\}$. The quotient by $\Delta^\prime$ is a parabolic induction of $aa(3,3)$. Let $G$ be $SO(14)$: $(L_1,L_1)\cong SL(4)$, $\dim C_1=1$, $\dim H_1-\dim H=4$ and $\dim C=1$. $(L,L)=(L_1,L_1)$ thus $\dim H^u_1 - \dim H^u = 4$. Here $\mathfrak n_1=\mathfrak m_1+\mathfrak m_2+\mathfrak m_3$ where these are simple $L_1$ modules, of dimensions $4$, $6$, $20$; $\mathfrak n$ is obtained from $\mathfrak n_1$ by removing $\mathfrak m_1$ and $C=C_1$. Therefore $H$ is uniquely determined up to conjugation.



\subsubsection*{Remaining cases}

For the following cases with more than one strongly $\Delta$-connected component we use the same arguments as above.

As in type $\mathsf A$, the following lemma is useful for the existence parts. We maintain the above notations.

\begin{lemma} [\cite{L01}] \label{L7} If $\Sigma^\prime\subset\Sigma/\Delta^\prime$ gives an isolated strongly $\Delta$-connected component of the spherical system $(S^p,\Sigma,\mathbf A)/\Delta^\prime$ such that its localization in $\supp(\Sigma^\prime)$ is isomorphic to $aa(p+q+p)$ or $aa(q)$, then: \begin{itemize}\item[-]if $q=1$, the corresponding spherical root $\alpha$ of type $a_1$ belongs to $\Sigma$ and $\Delta(\alpha)\cap\Delta(\beta)=\emptyset$ for all $\beta\in S$,\item[-] if $q>1$, the localization of the spherical system $(S^p,\Sigma,\mathbf A)/\Delta^\prime$ in $\supp(\Sigma/\Delta^\prime)$ is associated to $L_1$ with $(L_1,L_1)\cong K^\prime\times K^{\prime\prime}$, where $K^\prime\cong SL(p+1)\times SL(p+q)$, and if the factor of $K^\prime$ isomorphic to $SL(p+q)$ acts trivially on $\mathfrak n_1/\mathfrak n$, the corresponding spherical root $\gamma$ of type $a_q$ belongs to $\Sigma$.\end{itemize}\end{lemma}



$\mathbf{aa(1)+ds(n)}$, $n\geq4$. 
{\em Uniqueness}. Choose $\Delta^{\prime\prime}=\{\delta_{\alpha_1}^+,\delta_{\alpha_1}^-,\delta_{\alpha_2}\}$ and $\Delta^\prime=\{\delta_{\alpha_1}^+,\delta_{\alpha_2}\}$, so that $S^p/\Delta^{\prime\prime}=S\setminus\{\alpha_n,\alpha_{n+1}\}$ and $\Sigma/\Delta^\prime=\{\alpha_1+\ldots+\alpha_n\}$. The quotient by $\Delta^\prime$ is a parabolic induction of $aa(n)$. Let $G$ be $SO(2n+2)$. In this case we have $P\subsetneq Q$, since $Q=N(H^u_1)$ is the parabolic subgroup containing $B^-$ associated to $S^p/\Delta^{\prime\prime}\cup \{\alpha_n\}$, while $P=P_{-S^p/\Delta^{\prime\prime}}$. We have $L(P)\cong GL(n)\times GL(1)$, and $(L_1,L_1)\cong SL(n)$ with $\dim C_1=2$. Also, we have $\dim H_1- \dim H=n$ and $\dim C=2$. By Table~\ref{T4} the subgroup $(L,L)$ can be isomorphic only to $SL(n)$, $SO(n)$ or $Sp(n)$ (if $n$ is even). If $(L,L)\cong SO(n)$, respectively $Sp(n)$, we have that $H$ is included in a parabolic induction of $SO(n)$, respectively $Sp(n)$. In both cases we get a contradiction, since there are no distinguished subsets of $\Delta$ giving rise to the corresponding quotients. Therefore $(L,L)=(L_1,L_1)$ and $\dim H^u_1-\dim H^u=n$.
The Lie algebra $\mathfrak n_1$ decomposes into two simple $L_1$-submodules $\mathfrak m_1+\mathfrak m_2$ of dimensions $n$ and $n(n-1)/2$. The Lie subalgebra $\mathfrak n$ is equal to $\mathfrak m_2$ and $C=C_1$.

{\em Existence}. The subgroup $H$ is spherical and equal to $N_G(H)$, hence wonderful. We have $\dim P^u_{S^p}+\rank \Sigma=\dim B^u -(n-2)(n-3)/2 + 3$ and $\rank \Delta-\rank\Sigma=2$. Hence, $S^p\neq\emptyset$ and it can be proven that there are no root of type $d_m$ and that there is at least one root of type $a_m$. Now it is possible to see that $S^p$ is equal to $\{\alpha_3,\ldots,\alpha_{n-1}\}$ or $\{\alpha_2,\ldots,\alpha_{n-2}\}$ with at least a root of type $a_{n-1}$. One can easily prove that there cannot be a projective colour, by looking at the projective fibration it would be associated to, hence we have $aa(1)+ds(n)$ or $aa(n-1)+(2\mathrm{-comb})$. The latter can be excluded since it admits a morphism into a rank zero spherical system whose set of simple roots of type $p$ is equal to $\{\alpha_1,\ldots,\alpha_{n-2}\}\subsetneq S^p/\Delta^{\prime\prime}$.

$\mathbf{ac^\ast(n)+aa(1)}$. {\em Uniqueness}. Let $n$ be even. Take $\Delta^{\prime\prime}=\Delta\setminus\{\delta_{\alpha_1},\delta_{\alpha_{n+1}}^-\}$ and $\Delta^\prime=\{\delta_{\alpha_2},\ldots,\delta_{\alpha_{2i}},\ldots,\delta_{\alpha_n}\}$, so that $S^p/\Delta^{\prime\prime}=S\setminus\{\alpha_1,\alpha_{n+1}\}$ and $\Sigma/\Delta^\prime=\{\alpha_2+2\alpha_3+\alpha_4,\ldots,\alpha_{2i}+2\alpha_{2i+1}+\alpha_{2i+2},\ldots,\alpha_{n-2}+2\alpha_{n-1}+\alpha_n,\alpha_{n+1}\}$. The quotient by $\Delta^\prime$ is a parabolic induction of $dc(n)$. Let $G$ be $SO(2n+2)$, we have that $N(H^u_1)$ is the parabolic subgroup containing $B^-$ associated to $S\setminus\{\alpha_1\}=S^p/\Delta^{\prime\prime}\cup\{\alpha_{n+1}\}$. Here $(L_1,L_1)\cong SL(n)$ and $\dim C_1=2$. We have $\dim H_1-\dim H=n$ and $\dim C=2$, hence $(L,L)=(L_1,L_1)$. The Lie algebra $\mathfrak n_1$ decomposes into two conjugated simple $L_1$-submodules of dimension $n$. The Lie subalgebra $\mathfrak n$ is equal to one of them, and $H=N_G(H)$ so it is uniquely determined up to conjugation.

Let $n$ be odd and greater than $3$. Take $\Delta^{\prime\prime}=\Delta\setminus\{\delta_{\alpha_n},\delta_{\alpha_{n+1}}^+\}$ and $\Delta^\prime=\{\delta_{\alpha_2},\ \ldots$, $\delta_{\alpha_{n-1}}$, $\delta_{\alpha_{n+1}}^-\}$, so that $S^p/\Delta^{\prime\prime}=S\setminus\{\alpha_n,\alpha_{n+1}\}$ and $\Sigma/\Delta^\prime=\{\alpha_1+\ldots+\alpha_n\}$. The quotient spherical system is a parabolic induction of $aa(n)$. Let $G$ be $SO(2n+2)$, we have that $N(H^u_1)$ is the parabolic subgroup containing $B^-$ associated to $S\setminus\{\alpha_n\}=S^p/\Delta^{\prime\prime}\cup\{\alpha_{n+1}\}$; moreover $(L_1,L_1)\cong SL(n)$ and $\dim C_1=2$. We have $\dim H_1-\dim H=n(n-1)/2$ and $\dim C=2$, hence $(L,L)$ is equal to $(L_1,L_1)$. The Lie algebra $\mathfrak n_1$ decomposes into two simple $L_1$-submodules $\mathfrak m_1+\mathfrak m_2$ of dimension $n$ and $n(n-1)/2$, respectively. The Lie subalgebra $\mathfrak n$ is equal to $\mathfrak m_1$. The subgroup $H$ is equal to its normalizer and uniquely determined up to conjugation.

{\em Existence}. The subgroup $H$ is wonderful, we have $\dim P^u_{S^p}+\rank \Sigma=(n-1)(n+1)$ and $\rank\Delta-\rank\Sigma=2$. Hence, $S^p=\emptyset$, $\rank\Sigma=n-1$ and $\rank\Delta=n+1$. We have two cases with two strongly $\Delta$-connected components of type $ac^\ast(n)$ and $aa(1)$. One of such cases has a projective element and it can be excluded like in the previous case.

Let $n$ be equal to $3$. We have three different spherical systems of type $ac^\ast(3)+aa(1)$.
The proof is analogous to the case of $n>3$. In this case we have $n=n(n-1)/2$ and we have three subgroups up to conjugation: $H^\prime$ with $\mathfrak n=\mathfrak m_1$, $H^{\prime\prime}=\sigma(H^\prime)$ and $H^{\prime\prime\prime}$ with $\mathfrak n=\mathfrak m_2$ ($\sigma(H^{\prime\prime\prime})=H^{\prime\prime\prime}$). We can conclude analogously to the case of $do(p+2)$.

$\mathbf{2\!\!-\!\!comb+ac^\ast(3)}$. {\em Uniqueness}. Here $\Delta^{\prime\prime}=\{\delta_{\alpha_1}^+,\delta_{\alpha_3},\delta_{\alpha_4}^-\}$ and $\Delta^\prime=\{\delta_{\alpha_1}^+,\delta_{\alpha_3}\}$, and the quotient by $\Delta^\prime$ is a parabolic induction of $aa(3)$. Our $G$ is $SO(10)$, and $(L_1,L_1)\cong SL(3)$ with $\dim C_1=3$. We have $\dim H_1-\dim H=4$ and $\dim C=2$, hence $(L,L)=(L_1,L_1)$ and $\dim H^u_1-\dim H^u=3$. The $(L,L)$-module $\mathfrak n_1$ decomposes into three isotypic components, one is the sum of two $L_1$-submodules $\mathfrak m_1+\mathfrak m_2$ of dimension $3$ and of highest weight $\lambda_1$. The Lie algebra $\mathfrak n \subset \mathfrak n_1$ is the sum of a codimension $3$ submodule diagonal in $\mathfrak m_1+\mathfrak m_2$ and the other two isotypic components. Therefore $C$ is determined, $H$ is equal to its normalizer and uniquely determined up to conjugation.

{\em Existence}. The subgroup $H$ is wonderful. We have $\dim P^u_{S^p}+\rank \Sigma=\dim B^u + 4$ and $\rank \Delta-\rank\Sigma=2$. Hence, $S^p=\emptyset$ and there are two strongly $\Delta$-connected components. We know that $G/H\to G/H_1$ extends to a dominant $G$-morphism with connected fibers between the two wonderful embeddings: now $2\mathrm{-comb}+ac^\ast(3)$ is the unique spherical system satisfying the above requirements and admitting $(S^p,\Sigma,\mathbf A)/\Delta^\prime$ as quotient.

$\mathbf{ac^\ast(n_1)+d(n_2)}$, $n_2\geq2$. For $n_2=2$ we have $ac^\ast(n_1)+aa(1,1)$.
{\em Uniqueness}. Let $n_1$ be odd and $n_2$ be equal to $2$. Take $\Delta^{\prime\prime}=\Delta\setminus\{\delta_{\alpha_1},\delta_{\alpha_{n_1}}\}$ and $\Delta^\prime=\{\delta_{\alpha_2},\ldots,\delta_{\alpha_{2i}},\ldots,\delta_{\alpha_{n_1-1}}\}$; the quotient by $\Delta^\prime$ is a parabolic induction of two factors respectively of type $ac(n_1-2)$ and $aa(1,1)$. We have $G=SO(2n_1+4)$ and $L(Q)\cong GL(1)\times GL(n_1-1)\times SO(4)$. The universal cover of $(L_1,L_1)$ is isomorphic to $Sp(n_1-1)\times SL(2)$. The connected center $C_1$ has dimension $2$. We have $\dim H_1-\dim H=n$, $\dim C=1$ and $\dim (L_1,L_1)-\dim (L,L)\leq n-1$. Also, $(L,L)=(L_1,L_1)$ and $\dim H_1^u-\dim H^u=n_1-1$. The Lie subalgebra $\mathfrak n\subset \mathfrak n_1$ is the sum of a codimension $n_1-1$ submodule $\mathfrak m\subset \mathfrak m_1+\mathfrak m_2$ and the rest of the simple $L_1$-submodules. The submodule $\mathfrak m_2$ decomposes into two simple $L_1$-submodules $\mathfrak m_2=\mathfrak m_2^\prime+\mathfrak m_2^{\prime\prime}$, with $\mathfrak m_2^{\prime\prime}\cong\mathfrak m_1$ as $(L,L)$-modules. The submodule $\mathfrak m$ is the sum of $\mathfrak m_2^\prime$ and a submodule diagonal in $\mathfrak m_1+\mathfrak m_2^{\prime\prime}$.

Let $n_1$ be odd and $n_2$ be greater than $2$. The proof is completely analogous; the quotient by $\Delta^\prime$ is a parabolic induction of $ac(n_1-2)$ and $d(n_2)$. We have $G=SO(2n_1+2n_2)$ and $L(Q)\cong GL(1)\times GL(n_1-1)\times SO(2n_2)$. The subgroup $(L_1,L_1)$ is isomorphic to $Sp(n_1-1)\times SO(2n_2-1)$. As above, $\dim C_1=2$, $\dim H_1-\dim H=n$, $\dim C=1$, $\dim (L_1,L_1)-\dim (L,L)\leq n-1$ and $(L,L)=(L_1,L_1)$ with $\dim H_1^u-\dim H^u=n_1-1$. The Lie subalgebra $\mathfrak n\subset \mathfrak n_1$ is the sum of $\mathfrak m\subset \mathfrak m_1+\mathfrak m_2$ and the rest of the simple $L_1$-submodules. The submodule $\mathfrak m_2$ decomposes into $\mathfrak m_2=\mathfrak m_2^\prime+\mathfrak m_2^{\prime\prime}$, with $\mathfrak m_2^{\prime\prime}\cong\mathfrak m_1$ as $(L,L)$-modules. The submodule $\mathfrak m$ is the sum of $\mathfrak m_2^\prime$ and a submodule diagonal in $\mathfrak m_1+\mathfrak m_2^{\prime\prime}$.

Let $n_1$ be even. Let us treat together the cases $n_2=2$ and $n_2>2$. Take $\Delta^{\prime\prime}=\Delta\setminus\{\delta_{\alpha_{n_1}}\}$ and $\Delta^\prime=\{\delta_{\alpha_1},\ldots,\delta_{\alpha_{2i+1}},\ldots,\delta_{\alpha_{n_1-1}}\}$; the quotient by $\Delta^\prime$ is a parabolic induction of two factors respectively of type $ac(n_1-1)$ and $d(n_2)$. Here $(L_1,L_1)\cong Sp(n_1)\times SO(2n_2-1)$ (if $n_2=2$ its universal cover is $Sp(n_1)\times SL(2)$). Again $(L,L)=(L_1,L_1)$, and here $\dim H_1^u-\dim H^u=n_1$. The Lie subalgebra $\mathfrak n\subset \mathfrak n_1$ is the sum of a codimension $n_1$ submodule $\mathfrak m^\prime\subset \mathfrak m$ and the complementary $L_1$-submodule of $\mathfrak m$. The submodule $\mathfrak m$ decomposes into two simple $L_1$-submodules $\mathfrak m=\mathfrak m^\prime+\mathfrak m^{\prime\prime}$, with $\mathfrak m^{\prime\prime}$ of dimension $n_1$.

{\em Existence}. We have $\dim P^u_{S^p}+\rank \Sigma=\dim B^u - (n_2-1)(n_2-2) + n_1$ and $\rank \Delta-\rank\Sigma=1$. If $n_2=2$, $S^p=\emptyset$, $\rank\Sigma=n_1$ and $\rank\Delta=n_1+1$. The spherical system is cuspidal and there is a spherical root of type $aa(1,1)$ on the bifurcation. If $n_2>2$, analogously there is a spherical root of type $d(n_2)$ on the bifurcation and $S^p$ is determined.

$\mathbf{ay(p+q+p)+d(n)}$, $n\geq2$. {\em Uniqueness}. Here $\Delta^{\prime\prime}=\Delta\setminus\{\delta_{\alpha_{p+q}},\delta_{\alpha_{2p+q}}^-\}$, and $\Delta^\prime=\{\delta_{\alpha_1}^+,\ldots,\delta_{\alpha_p}^+\}$ is such that the quotient is a parabolic induction of $aa((p-1)+(q+1)+(p-1))$ and $d(n)$. Let $G$ be $SO(4p+2q+2n)$. In this case $N(H^u_1)$ is the parabolic subgroup containing $B^-$ associated to $S\setminus\{\alpha_{2p+q}\}=S^p/\Delta^{\prime\prime}\cup\{\alpha_{p+q}\}$. In this case $(L_1,L_1)\cong SL(p+q)\times SL(p)\times SO(2n-1)$ (if $n=2$ its universal cover is $SL(p+q)\times SL(p)\times SL(2)$) and $\dim C_1=2$. We have $\dim H_1-\dim H=p+q$ and $\dim C=2$. By Table~\ref{T4} the subgroup $(L,L)$ can be isomorphic only to $SL(p+q)\times SL(p)\times SO(2n-1)$, $SL(p+q)\times SO(p)\times SO(2n-1)$,  $SL(p+q)\times SO(p)\times SO(2n-1)$ (for $p$ even), $SL(p+q)\times SL(2)$ (for $p=n=2$), $SL(p+q)\times SL(p)\times SO(r)\times SO(2n-1-r)$, $SL(p+q)\times SL(p)\times Spin(7)$ (for $n=5$), $SL(p+q)\times SL(p)\times G_2$ (for $n=4$) and some intersections of two of them. Any choice with $(L,L)\neq (L_1,L_1)$ produces a $H$ contained in a parabolic induction, and none of these cases corresponds to any quotient of our spherical system. Therefore, $(L,L)=(L_1,L_1)$ and $\dim H^u_1-\dim H^u=p+q$. The Lie subalgebra $\mathfrak n\subset \mathfrak n_1$ is the sum of a codimension $p+q$ submodule $\mathfrak m^\prime\subset \mathfrak m$ and the complementary $L_1$-submodule of $\mathfrak m$. The submodule $\mathfrak m$ decomposes into two simple $L_1$-submodules $\mathfrak m=\mathfrak m^\prime+\mathfrak m^{\prime\prime}$, with $\mathfrak m^{\prime\prime}$ of dimension $p+q$.

{\em Existence}. We have $\dim P^u_{S^p}+\rank\Sigma=\dim B^u -(q-1)(q-2)/2 -(n-1)(n-2) +2p+2$ and $\rank\Delta-\rank\Sigma=2$. We know that $G/H\to G/H_1$ extends to a dominant $G$-morphism with connected fibers between the corresponding wonderful embeddings, hence $S^p\subset S^p/\Delta^\prime$. If $q>2$, it can be shown that the inclusion is proper, there is exactly one root of type $d_n$ on the bifurcation and one root of type $a_q$. The rest is a strongly $\Delta$-connected component of roots of type $a_1$. We have only three cases: $ay(p+q+p)+d(n)$, if $p=3$ $az^\sim(3+q+3)+d(n)$ or if $p=2$ $ay^\ast(2+q+2)+d(n)$. The proof can be concluded thanks to the corresponding uniqueness parts. If $q\leq2$, we know analogously that there are three strongly $\Delta$-connected components and the proof ends by noticing that $ay(p+q+p)+d(n)$ is the unique spherical system admitting $(S^p,\Sigma,\mathbf A)/\Delta^{\prime\prime}$ as quotient.

$\mathbf{ay(p,p)+d(n_1)+d(n_2)}$, $n_1,n_2\geq2$. {\em Uniqueness}. Here $\Delta^{\prime\prime}=\Delta\setminus\{\delta_{\alpha_p}^-,\delta_{\alpha^\prime_p}^-\}$, and $\Delta^\prime=\{\delta_{\alpha_1}^+,\ldots,\delta_{\alpha_p}^+\}$ is such that the quotient is a parabolic induction of $aa(p-1,p-1)$, $d(n_1)$ and $d(n_2)$. Let $G$ be $SO(2p+2n_1)\times SO(2p+2n_2)$; we have $(L_1,L_1)\cong SL(p)\times SO(2n_1-1)\times SO(2n_2-1)$ and $\dim C_1=2$. We have $\dim H_1-\dim H=p+1$, $\dim C=1$ and $\dim (L_1,L_1)-\dim (L,L)\leq p$. As usual $(L,L)=(L_1,L_1)$ and $\dim H_1^u-\dim H^u=p$. The Lie subalgebra $\mathfrak n\subset \mathfrak n_1$ is the sum of $\mathfrak m\subset \mathfrak m_1+\mathfrak m_2$ and the rest of the simple $L_1$-submodules. The submodules $\mathfrak m_1$ and $\mathfrak m_2$ decompose into $\mathfrak m_1=\mathfrak m_1^\prime+\mathfrak m_1^{\prime\prime}$ and $\mathfrak m_2=\mathfrak m_2^\prime+\mathfrak m_2^{\prime\prime}$, with $\mathfrak m_1^{\prime\prime}\cong\mathfrak m_2^{\prime\prime}$ as $(L,L)$-modules. The submodule $\mathfrak m$ is the sum of $\mathfrak m_1^\prime+\mathfrak m_2^\prime$ and a submodule diagonal in $\mathfrak m_1^{\prime\prime}+\mathfrak m_2^{\prime\prime}$.

{\em Existence}. We have $\dim P^u_{S^p}+\rank\Sigma=\dim B^u -(n_1-1)(n_1-2)-(n_2-1)(n_2-2) +2p+2$ and $\rank\Delta-\rank\Sigma=1$. If $n_1=n_2=2$ the unique possibility is $\rank\Sigma=2p+2$ and $S^p=\emptyset$. Moreover, we can deduce that the spherical system consists of a strongly $\Delta$-connected component with roots of type $a_1$ on a Dynkin diagram of type $A_p\times A_p$ plus two roots of type $a_1\times a_1$ on the bifurcations. We can easily exclude the case of $\mathrm {comb}+aa(1,1)+aa(1,1)$. The are only two other possibilities: $ay(p,p)+aa(1,1)+aa(1,1)$ and $az(3,3)+aa(1,1)+aa(1,1)$ if $p=3$. The latter can be excluded by uniqueness. If $n_1>2$ or $n_2>2$, we can deduce that $S^p\neq\emptyset$ and hence that $\Sigma$ contains at least one root of type $d_m$. Moreover, no root of type $a_m$ with $m\geq2$ or $a^\prime_1$ can occur. We can have two roots of type $d_m$, or one root of type $d_m$ and one root of type $a_1\times a_1$, in both cases both roots on the bifurcations. All the other spherical roots are of type $a_1$ lying in a strongly $\Delta$-connected component. 

$\mathbf{ay(p,p)+d(n)}$, $n\geq2$. {\em Uniqueness}. Take $\Delta^{\prime\prime}=\Delta\setminus\{\delta_{\alpha_p}^-,\delta_{\alpha^\prime_p}^-\}$ and $\Delta^\prime=\{\delta_{\alpha_1}^+,\ldots,\delta_{\alpha_p}^+\}$. The quotient by $\Delta^\prime$ is a parabolic induction of $aa(p-1,p-1)$ and $d(n)$. Let $G$ be $SL(p+1)\times SO(2p+2n)$. The subgroup $(L_1,L_1)$ is isomorphic to $SL(p)\times SO(2n-1)$ and $\dim C_1=2$. We have $\dim H_1-\dim H=p+1$, $\dim C=1$. The subgroup $(L,L)$ is equal to $(L_1,L_1)$ and $\dim H_1^u-\dim H^u=p$. The Lie subalgebra $\mathfrak n\subset \mathfrak n_1$ is the sum of $\mathfrak m\subset \mathfrak m_1+\mathfrak m_2$ and the rest of the simple $L_1$-submodules. The submodule $\mathfrak m_2$ decompose into $\mathfrak m_2=\mathfrak m_2^\prime+\mathfrak m_2^{\prime\prime}$, with $\mathfrak m_1\cong\mathfrak m_2^{\prime\prime}$ as $(L,L)$-modules. The submodule $\mathfrak m$ is the sum of $\mathfrak m_2^\prime$ and a submodule diagonal in $\mathfrak m_1+\mathfrak m_2^{\prime\prime}$.

$\mathbf{ay(p,(p-1))+d(n)}$, $n\geq2$. {\em Uniqueness}. Analogously, set $\Delta^{\prime\prime}=\Delta\setminus\{\delta_{\alpha^\prime_p}^-\}$. The distinguished subset $\Delta^\prime=\{\delta_{\alpha^\prime_1}^+,\ldots,\delta_{\alpha^\prime_p}^+\}$ is such that the quotient spherical system is a parabolic induction of two factors respectively of type $aa(p-1,p-1)$ and $d(n)$. Let $G$ be $SL(p)\times SO(2p+2n)$. The subgroup $(L_1,L_1)$ is isomorphic to $SL(p)\times SO(2n-1)$ and $\dim C_1=1$. We have $\dim H_1-\dim H=p$, $\dim C=1$. The subgroup $(L,L)$ is equal to $(L_1,L_1)$ and $\dim H_1^u-\dim H^u=p$. The Lie subalgebra $\mathfrak n\subset \mathfrak n_1$ is the sum of $\mathfrak m^\prime\subset \mathfrak m$ and the complementary $L_1$-submodule of $\mathfrak m$. The submodule $\mathfrak m$ decomposes into two simple $L_1$-submodules $\mathfrak m=\mathfrak m^\prime+\mathfrak m^{\prime\prime}$, with $\mathfrak m^{\prime\prime}$ of dimension $p$.

$\mathbf{ay(p+q+(p-1))+d(n)}$, $n\geq2$. {\em Uniqueness}. Take $\Delta^{\prime\prime}=\Delta\setminus\{\delta_{\alpha_p}$, $\delta_{\alpha_{2p+q-1}}^-\}$ and $\Delta^\prime=\{\delta_{\alpha_{p+q}}^+, \ldots$, $\delta_{\alpha_{2p+q-1}}^+\}$. The quotient by $\Delta^\prime$ is a parabolic induction of two factors respectively of type $aa((p-1)+q+(p-1))$ and $d(n)$. Let $G$ be $SO(4p+2q-2+2n)$. The subgroup $(L_1,L_1)$ is isomorphic to $SL(p)\times SL(p+q-1)\times SO(2n-1)$ and $\dim C_1=2$. We have $\dim H_1-\dim H=p$ and $\dim C=2$. The subgroup $(L,L)$ is equal to $(L_1,L_1)$. The Lie subalgebra $\mathfrak n\subset \mathfrak n_1$ is the sum of $\mathfrak m^\prime\subset \mathfrak m$ and the complementary $L_1$-submodule of $\mathfrak m$. The submodule $\mathfrak m$ decomposes into two simple $L_1$-submodules $\mathfrak m=\mathfrak m^\prime+\mathfrak m^{\prime\prime}$, with $\mathfrak m^{\prime\prime}$ of dimension $p$.

{\em Existence}. By Lemma~\ref{L7} the root $\alpha_4+\ldots+\alpha_{q+3}$ belongs to $\Sigma$. Moreover, $S^p\supset\{\alpha_5,\ldots,\alpha_{q+2}\}$ and there is a root of type $d_n$ ($a_1\times a_1$ if $n=2$) on the bifurcation. We obtain two cases: $az^\sim(3+q+2)+aa(1,1)$ if $p=3$ and $ay(p+q+(p-1))+aa(1,1)$. The former can be excluded by the corresponding uniqueness proof.

$\mathbf{dy(p+q+(p-1))}$, $p\geq3$. {\em Uniqueness}. Here $\Delta^{\prime\prime}=\Delta\setminus\{\delta_{\alpha_{p+q-1}}$, $\delta_{\alpha_{2p+q-2}}^-$, $\delta_{\alpha_{2p+q-1}}^-\}$ and $\Delta^\prime=\{\delta_{\alpha_1}^+,\ldots$, $\delta_{\alpha_{p-1}}^+$, $\delta_{\alpha_p}\}$; the quotient by $\Delta^\prime$ is a parabolic induction of $aa((p-2)+(q+1)+(p-2))$. Let $G$ be $SO(4p+2q-2)$. In this case $N(H^u_1)$ is the parabolic subgroup containing $B^-$ associated to $S\setminus\{\alpha_{2p+q-2},\alpha_{2p+q-1}\}=S^p/\Delta^{\prime\prime}\cup\{\alpha_{p+q-1}\}$. The subgroup $(L_1,L_1)$ is isomorphic to $SL(p+q-1)\times SL(p-1)$ and $\dim C_1=3$. We have $\dim H_1-\dim H=p+q$ and $\dim C=2$. The subgroup $(L,L)$ can be isomorphic only to $SL(p+q-1)\times SO(p-1)$, $SL(p+q-1)\times Sp(p-1)$ (for $p$ odd) and $SL(p+q-1)\times SL(p-1)$. If $(L,L)$ is isomorphic to one of the above proper spherical subgroups, $H$ is contained in a parabolic induction, with connected quotient. From both the cases we get a contradiction. Therefore, the subgroup $(L,L)$ is equal to $(L_1,L_1)$ and $\dim H^u_1-\dim H^u=p+q-1$. 

{\em Existence}. We have $\dim P^u_{S^p}+\rank\Sigma=\dim B^u -(q-1)(q-2)/2+2p$ and $\rank\Delta-\rank\Sigma=2$. We know that $H\subset H_1$ and that the quotient $G/H\to G/H_1$ extends to a dominant $G$-morphism with connected fibers. If $q>2$, there are no root of type $d_m$ and there is exactly one root of type $a_q$. The rest is a strongly $\Delta$-connected component of roots of type $a_1$. We have only three cases: $dy(p+q+(p-1))$ and, if $p=4$, $dy^\sim(4+q+3)$ or $dz^\sim(4+q+3)$. The proof can be concluded by uniqueness. If $q\leq2$, we know analogously that there are two strongly $\Delta$-connected components and the proof ends by noticing that $dy(p+q+(p-1))$ is the unique spherical system admitting $(S^p,\Sigma,\mathbf A)/\Delta^{\prime\prime}$ as quotient.

$\mathbf{dy(p,(p-1))+d(n)}$, $p\geq 4$ and $n\geq2$. {\em Uniqueness}. In this case $\Delta^{\prime\prime}=\Delta\setminus\{\delta_{\alpha_{p-1}}^-,\delta_{\alpha^\prime_{p-1}}^-,\delta_{\alpha^\prime_p}^-\}$, and $\Delta^\prime=\{\delta_{\alpha_1}^+,\ldots,\delta_{\alpha_{p-1}}^+\}$ is such that the quotient spherical system is a parabolic induction of $aa(p-1,p-1)$ and $d(2)$. Let $G$ be $SO(2p-2+2n)\times SO(2p)$. The subgroup $(L_1,L_1)$ is isomorphic to $SL(p-1)\times SO(2n-1)$ and $\dim C_1=3$. We have $\dim H_1-\dim H=p+1$, $\dim C=1$ and $\dim (L_1,L_1)-\dim (L,L)\leq p-1$. The subgroup $(L,L)$ is equal to $(L_1,L_1)$ and $\dim H^u_1-\dim H^u=p-1$. 

{\em Existence}. We have $\dim P^u_{S^p}+\rank\Sigma=\dim B^u +2p$ and $\rank\Delta-\rank\Sigma=1$. Analogously, we get only the cases $dz(4,3)+d(n)$ (if $p=4$) and $dy(p,(p-1))+d(n)$. The former is excluded by the corresponding uniqueness proof.

$\mathbf{dy^\prime(p+q+(p-2))}$, $p\geq4$. 
{\em Uniqueness}. Take $\Delta^{\prime\prime}=\Delta\setminus\{\delta_{\alpha_{p-1}}$, $\delta_{\alpha_{2p+q-3}}^-$, $\delta_{\alpha_{2p+q-2}}^-\}$ and $\Delta^\prime=\{\delta_{\alpha_1}^+,\ldots$, $\delta_{\alpha_{p-2}}^+$, $\delta_{\alpha_{p+q-1}}\}$. The quotient by $\Delta^\prime$ is a parabolic induction of $aa((p-2)+q+(p-2))$. Let $G$ be $SO(4p+2q-4)$. In this case $N(H^u_1)$ is the parabolic subgroup containing $B^-$ associated to $S\setminus\{\alpha_{2p+q-3},\alpha_{2p+q-2}\}=S^p/\Delta^{\prime\prime}\cup\{\alpha_{p+q-2}\}$. The subgroup $(L_1,L_1)$ is isomorphic to $SL(p-1)\times SL(p+q-2)$ and $\dim C_1=3$. We have $\dim H_1-\dim H=p$ and $\dim C=2$. The subgroup $(L,L)$ is equal to $(L_1,L_1)$ and $\dim H^u_1-\dim H^u=p-1$.

{\em Existence}. By Lemma~\ref{L7} the root $\alpha_{p-1}+\ldots+\alpha_{p+q-1}$ belongs to $\Sigma$. We get two cases: $dy^\prime(p+q+(p-2))$ and $dz^\sim(4+q+2)$ (if $p=4$). The latter is excluded by the corresponding uniqueness proof.

$\mathbf{ay^\ast(2+q+2)+d(n)}$, $n\geq2$. {\em Uniqueness}. Set $\Delta^{\prime\prime}=\{\delta_{\alpha_1}^+,\delta_{\alpha_2}^+,\delta_{\alpha_2}^-,\delta_{\alpha_3},\delta_{\alpha_{q+5}}\}$ and $\Delta^\prime=\{\delta_{\alpha_1}^+,\delta_{\alpha_2}^+\}$. The quotient by $\Delta^\prime$ is a parabolic induction of two factors respectively of type $aa(1+q+1)$ and $d(n)$. Let $G$ be $SO(2q+8+2n)$. The subgroup $(L_1,L_1)$ is isomorphic to $SL(2)\times SL(q+1)\times SO(2n-1)$ and $\dim C_1=3$. We have $\dim H_1-\dim H=3$ and $\dim C=2$. The subgroup $(L,L)$ is equal to $(L_1,L_1)$.
 
{\em Existence}. By Lemma~\ref{L7} the root $\alpha_3+\ldots+\alpha_{q+2}$ belongs to $\Sigma$. Moreover, $S^p\supseteq\{\alpha_4,\ldots,\alpha_{q+1}\}$ and there is a root of type $d_n$ ($a_1\times a_1$ if $n=2$) on the bifurcation. We obtain two cases: $ay(2+q+2)+d(n)$ and $ay^\ast(2+q+2)+d(n)$. The former can be excluded by the corresponding uniqueness proof.

$\mathbf{dy^\ast(3+q+3)}$. {\em Uniqueness}. Choose $\Delta^{\prime\prime}=\Delta\setminus\{\delta_{\alpha_1}^+,\delta_{\alpha_{q+3}},\delta_{\alpha_{q+6}}^+\}$ and $\Delta^\prime=\{\delta_{\alpha_1}^-,\delta_{\alpha_2}^-,\delta_{\alpha_3}^-\}$; the quotient by $\Delta^\prime$ is a parabolic induction of $aa(2+q+2)$. Let $G$ be $SO(2q+12)$. The subgroup $(L_1,L_1)$ is isomorphic to $SL(3)\times SL(q+2)$ and $\dim C_1=3$. We have $\dim H_1-\dim H=4$ and $\dim C=2$. The subgroup $(L,L)$ is equal to $(L_1,L_1)$. 

{\em Existence}. By Lemma~\ref{L7} the root $\alpha_4+\ldots+\alpha_{q+3}$ belongs to $\Sigma$. Moreover, $S^p=\{\alpha_5,\ldots,\alpha_{q+2}\}$. The unique possible spherical system is $dy^\ast(3+q+3)$.

$\mathbf{dy^\ast(3+q+2)}$. {\em Uniqueness}. Choose $\Delta^{\prime\prime}=\Delta\setminus\{\delta_{\alpha_{q+2}},\delta_{\alpha_{q+5}}^+\}$, and $\Delta^\prime=\{\delta_{\alpha_1}^-,\delta_{\alpha_2}^-,\delta_{\alpha_{q+4}}^-\}$ is such that the quotient spherical system is a parabolic induction of $aa(2+q+2)$. Let $G$ be $SO(2q+10)$. The subgroup $(L_1,L_1)$ is isomorphic to $SL(3)\times SL(q+2)$ and $\dim C_1=2$. We have $\dim H_1-\dim H=3$ and $\dim C=2$. The subgroup $(L,L)$ is equal to $(L_1,L_1)$.

$\mathbf{dy^\sim(4+q+3)}$. {\em Uniqueness}. Choose $\Delta^{\prime\prime}=\{\delta_{\alpha_1}^+,\delta_{\alpha_1}^-,\delta_{\alpha_2}^+,\delta_{\alpha_2}^-,\delta_{\alpha_3}^-,\delta_{\alpha_4}\}$ and $\Delta^\prime=\{\delta_{\alpha_1}^-,\delta_{\alpha_2}^-,\delta_{\alpha_3}^-\}$. The quotient by $\Delta^\prime$ is a parabolic induction of $aa(2,2)$ and $aa(q)$. Let $G$ be $SO(2q+14)$. In this case $N(H^u_1)$ is the parabolic subgroup containing $B^-$ associated to $S^p/\Delta^{\prime\prime}\cup\{\alpha_{q+3}\}$. The subgroup $(L_1,L_1)$ is isomorphic to $SL(3)\times SL(q)$ and $\dim C_1=4$. We have $\dim H_1-\dim H=5$ and $\dim C=2$. By Table~\ref{T4} the subgroup $(L,L)$ is equal to $(L_1,L_1)$ and $\dim H^u_1-\dim H^u=3$. If $q\geq2$, the Lie subalgebra $\mathfrak n \subset \mathfrak n _1$ is the sum of a codimension $3$ submodule $\mathfrak m\subset\mathfrak m_1+\mathfrak m_2+\mathfrak m_3$ and the rest of the simple $L_1$-submodules.
If $q=1$, we can substitute the colours $\delta_{\alpha_4}$ and $\delta_{\alpha_q+3}$ with $\delta_{\alpha_4}^+$ and $\delta_{\alpha_4}^-$, respectively. In the decomposition of $\mathfrak n_1$ we get two other submodules, $\mathfrak m_4$ and $\mathfrak m_5$, isomorphic to $\mathfrak m_1$, $\mathfrak m_2$, $\mathfrak m_3$ as $(L_1,L_1)$-modules. As above we need three of them where to take a codimension two diagonal submodule. If $\mathfrak n$ contains $\mathfrak m_3$ we have that $H^u$ contains the unipotent radical of a proper parabolic subgroup and it is associated to a parabolic induction, but this is not the case. The same thing occurs if $\mathfrak n$ contains  $\mathfrak m_4+\mathfrak m_1$ or $\mathfrak m_5+\mathfrak m_2$. Analogously, if $\mathfrak n$ contains $\mathfrak m_4+\mathfrak m_2$ or $\mathfrak m_1+\mathfrak m_5$, we have that $H$ is equal to a wonderful subgroup associated to the following, \[\begin{picture}(11850,5700)\multiput(300,3000)(1800,0){3}{\usebox{\segm}}\multiput(5700,3000)(1800,0){2}{\usebox{\segm}}\put(9300,1200){\usebox{\longbifurc}}\multiput(0,2100)(1800,0){3}{\usebox{\aone}}\multiput(7200,2100)(1800,0){2}{\usebox{\aone}}\multiput(10800,300)(0,3600){2}{\usebox{\aone}}\put(5400,2100){\usebox{\aone}}\put(9300,3600){\usebox{\tose}}
\put(11100,5400){\usebox{\tosw}}
\multiput(1200,3400)(1800,0){2}{\put(150,0){\usebox{\tow}}}\multiput(8400,3400)(-1800,0){2}{\put(150,0){\usebox{\tow}}}\put(4200,3400){\put(50,0){\usebox{\toe}}}
\multiput(300,3900)(9000,0){2}{\line(0,1){750}}\put(300,4650){\line(1,0){9000}}\put(2100,3900){\line(0,1){650}}\put(2100,4750){\line(0,1){650}}\put(2100,5400){\line(1,0){8700}}\multiput(300,2100)(7200,0){2}{\line(0,-1){1500}}\put(300,600){\line(1,0){10500}}\multiput(2100,2100)(7200,0){2}{\line(0,-1){900}}\put(2100,1200){\line(1,0){5300}}\put(7600,1200){\line(1,0){1700}}\put(3900,2100){\line(0,-1){700}}\put(3900,1100){\line(0,-1){400}}\put(3900,500){\line(0,-1){500}}\put(3900,0){\line(1,0){7950}}\put(11850,0){\line(0,1){4200}}\put(11850,4200){\line(-1,0){450}}\end{picture}\] 
but this is not the case; moreover, notice that in $\Delta(\alpha_4)$ there is a projective element and in this case $H$ would be included (with connected quotient) in a wonderful subgroup associated to a parabolic induction of $dy(4,3)$. Therefore, there are only two other conjugated cases, where $\mathfrak m$ is included in $\mathfrak m_1+\mathfrak m_2+\mathfrak m_3$ or $\mathfrak m_4+\mathfrak m_5+\mathfrak m_3$.  

{\em Existence}. We have $\dim P^u_{S^p}+\rank\Sigma=\dim B^u -(q-1)(q-2)/2+8$ and $\rank\Delta-\rank\Sigma=2$. By Lemma~\ref{L7} we get that the root $\alpha_4+\ldots+\alpha_{q+3}$ belongs to $\Sigma$. Moreover, $S^p=\{\alpha_5,\ldots,\alpha_{q+2}\}$ and we get only three cases: $dz^\sim(4+q+3)$, $dy(4+q+3)$ and $dy^\sim(4+q+3)$. The former two cases can be excluded by the corresponding uniqueness proofs.

$\mathbf{ay(2,2)+2\!\!-\!\!comb}$. {\em Uniqueness}. Choose $\Delta^{\prime\prime}=\{\delta_{\alpha_1}^-,\delta_{\alpha_2}^+,\delta_{\alpha_2}^-,\delta_{\alpha_3}^+,\delta_{\alpha_3}^-\}$ and $\Delta^\prime=\{\delta_{\alpha_1}^-,\delta_{\alpha_2}^-,\delta_{\alpha_3}^-\}$. The quotient by $\Delta^\prime$ is a parabolic induction of $aa(1+2+1)$. Let $G$ be $SO(12)$. The subgroup $(L_1,L_1)$ is isomorphic to $SL(2)\times SL(3)$ and $\dim C_1=3$. We have $\dim H_1-\dim H=4$, $\dim C=2$, $(L,L)=(L_1,L_1)$ and $\dim H^u_1-\dim H^u=3$.

{\em Existence}. From $\dim P^u_{S^p}+\rank\Sigma=\dim B^u+6$ and $\rank\Delta-\rank\Sigma=2$ we get $S^p=\emptyset$, $\Sigma=S$. We have two strongly $\Delta$-connected components. The spherical system $ay(2,2)+(2\mathrm{-comb})$ is the unique that admits $(S^p,\Sigma,\mathbf A)/\Delta^{\prime\prime}$ as quotient.

$\mathbf{ay(2,1)+2\!\!-\!\!comb}$. {\em Uniqueness}. We take $\Delta^{\prime\prime}=\Delta\setminus\{\delta_{\alpha_3}^+,\delta_{\alpha_5}^+\}$ and $\Delta^\prime=\{\delta_{\alpha_1}^-,\delta_{\alpha_2}^-,\delta_{\alpha_4}^-\}$; the quotient by the latter is a parabolic induction of $aa(1+2+1)$. Let $G$ be $SO(10)$. The subgroup $(L_1,L_1)$ is isomorphic to $SL(2)\times SL(3)$ and $\dim C_1=2$. We have $\dim H_1-\dim H=3$, $\dim C=2$, $(L,L)=(L_1,L_1)$ and $\dim H^u_1-\dim H^u=3$.

$\mathbf{az(3,3)+d(n_1)+d(n_2)}$, $n_1,n_2\geq2$. {\em Uniqueness}. Choose $\Delta^{\prime\prime}=\{\delta_{\alpha_1}^+$, $\delta_{\alpha_2}^+$, $\delta_{\alpha_2}^-$, $\delta_{\alpha_4}$, $\delta_{\alpha^\prime_4}\}$ and $\Delta^\prime=\{\delta_{\alpha_1}^+,\delta_{\alpha_2}^+\}$, so that the quotient by $\Delta^\prime$ is a parabolic induction of three factors of type $aa(1,1)$, $d(n_1)$ and $d(n_2)$. Let $G$ be $SO(6+2n_1)\times SO(6+2n_2)$. The subgroup $(L_1,L_1)$ is isomorphic to $SL(2)\times SO(2n_1-1)\times SO(2n_2-1)$ and $\dim C_1=4$. The subgroup $(L,L)$ is equal to $(L_1,L_1)$ and $\dim H^u_1 - \dim H^u = 2$.

{\em Existence}. If $n_1=n_2=2$, we have $\dim P^u_{S^p}+\rank \Sigma=\dim B^u + 8$ and the unique possibility is $\rank\Sigma=8$ and $S^p=\emptyset$. Moreover, we can deduce that the spherical system consists of a strongly $\Delta$-connected component with roots of type $a_1$ on a Dynkin diagram of type $A_3\times A_3$ plus two roots of type $a_1\times a_1$ on the bifurcations. There are only two other possibilities: $az(3,3)+aa(1,1)+aa(1,1)$ and $ay(3,3)+aa(1,1)+aa(1,1)$. Analogously, if $n_1>2$ or $n_2>2$ we can prove that there are two roots of type $d_{n_1}$ and $d_{n_2}$ on the bifurcations and all the other spherical roots are of type $a_1$ lying in a strongly $\Delta$-connected component.

$\mathbf{az(3,2)+d(n_1)+d(n_2)}$, $n_1,n_2\geq2$. {\em Uniqueness}. We can take $\Delta^{\prime\prime}=\{\delta_{\alpha_1}^+$, $\delta_{\alpha_1}^-$, $\delta_{\alpha_2}^+$, $\delta_{\alpha_3}$, $\delta_{\alpha^\prime_4}\}$ and $\Delta^\prime=\{\delta_{\alpha_1}^+,\delta_{\alpha_2}^+\}$. 

$\mathbf{az(3,3)+d(n)}$, $n\geq2$. {\em Uniqueness}. Here $\Delta^{\prime\prime}=\{\delta_{\alpha_1}^+,\delta_{\alpha_2}^+,\delta_{\alpha_2}^-,\delta_{\alpha^\prime_4}\}$ so that $S^p/\Delta^{\prime\prime}=\{\alpha_2,\alpha^\prime_2,\alpha^\prime_4,\alpha^\prime_5\}$, and $\Delta^\prime=\{\delta_{\alpha_1}^+,\delta_{\alpha_2}^+\}$ is such that the quotient spherical system is a parabolic induction of two factors of type $aa(1,1)$ and $d(n)$.

$\mathbf{az(3,2)+^1d(n)}$ $n\geq2$. {\em Uniqueness}. Here $\Delta^{\prime\prime}=\{\delta_{\alpha_1}^+,\delta_{\alpha_1}^-,\delta_{\alpha_2}^+,\delta_{\alpha^\prime_4}\}$ and $\Delta^\prime=\{\delta_{\alpha_1}^+,\delta_{\alpha_2}^+\}$. The quotient by $\Delta^\prime$ is a parabolic induction of two factors of type $aa(1,1)$ and $d(n)$.

$\mathbf{az(3,2)+^2d(n)}$, $n\geq2$. {\em Uniqueness}. Here $\Delta^{\prime\prime}=\{\delta_{\alpha_1}^+,\delta_{\alpha_2}^+,\delta_{\alpha_2}^-,\delta_{\alpha^\prime_3}\}$ and $\Delta^\prime=\{\delta_{\alpha_1}^+,\delta_{\alpha_2}^+\}$. The quotient by $\Delta^\prime$ is a parabolic induction of two factors of type $aa(1,1)$ and $d(n)$. 

$\mathbf{az(3,1)+d(n)}$, $n\geq2$. {\em Uniqueness}. Here $\Delta^{\prime\prime}=\{\delta_{\alpha_1}^+,\delta_{\alpha_1}^-,\delta_{\alpha^\prime_2}^+,\delta_{\alpha^\prime_4}\}$ and $\Delta^\prime=\{\delta_{\alpha_1}^+,\delta_{\alpha^\prime_2}^+\}$. The quotient by $\Delta^\prime$ is a parabolic induction of two factors of type $aa(1,1)$ and $d(n)$. 

$\mathbf{az^\sim(3+q+3)+d(n)}$, $n\geq2$. {\em Uniqueness}. Here we choose $\Delta^{\prime\prime}=\{\delta_{\alpha_1}^+$, $\delta_{\alpha_2}^+$, $\delta_{\alpha_2}^-$, $\delta_{\alpha_{q+3}}$, $\delta_{\alpha_{q+7}}\}$ and $\Delta^\prime=\{\delta_{\alpha_1}^+,\delta_{\alpha_2}^+\}$. The quotient by $\Delta^\prime$ is a parabolic induction of three factors of type $aa(1,1)$, $aa(q)$ and $d(n)$. Let $G$ be $SO(2q+12+2n)$. We have $\dim C_1=5$, $\dim H_1-\dim H=5$, $\dim C=2$ and $(L,L)=(L_1,L_1)$. In particular, $\dim H^u_1 - \dim H^u = 2$. If $q\geq2$, the Lie subalgebra $\mathfrak n \subset \mathfrak n _1$ is the sum of a codimension $2$ submodule $\mathfrak m\subset\mathfrak m_1+\mathfrak m_2+\mathfrak m_3+\mathfrak m_4$ and the rest of the simple $L_1$-submodules. If $q=1$, in the decomposition of $\mathfrak n_1$ we get two other modules isomorphic to $\mathfrak m_1$,  $\mathfrak m_2$, $\mathfrak m_3$ and $\mathfrak m_4$ as $(L_1,L_1)$-modules. As above we need four of them where to take a codimension $2$ diagonal submodule. The proof is completely analogous to that of $dy^\sim(4+q+3)$.

{\em Existence}. By Lemma~\ref{L7} we get that the root $\alpha_4+\ldots+\alpha_{q+3}$ of type $a_q$ belongs to $\Sigma$. Therefore, if $n=2$ we have $S^p=\{\alpha_5,\ldots,\alpha_{q+2}\}$ and there are two cases: $az(3,3)+aa(1,1)$ and $ay(3,3)+aa(1,1)$. The latter can be excluded by the corresponding uniqueness proof. If $n>2$ it can be deduced that $S^p$ contains properly $\alpha_5,\ldots,\alpha_{q+2}$ and there is a root of type $d_n$ on the bifurcation.

$\mathbf{az^\sim(3+q+2)+d(n)}$, $n\geq2$. {\em Uniqueness}. We can take $\Delta^{\prime\prime}=\{\delta_{\alpha_1}^+$, $\delta_{\alpha_1}^-$, $\delta_{\alpha_2}^+$, $\delta_{\alpha_3}$, $\delta_{\alpha_{q+6}}\}$ and $\Delta^\prime=\{\delta_{\alpha_1}^+,\delta_{\alpha_2}^+\}$. 

$\mathbf{dz(4,3)+d(n)}$, $n\geq2$. {\em Uniqueness}. Here $\Delta^{\prime\prime}=\{\delta_{\alpha_1}^+,\delta_{\alpha_2}^+,\delta_{\alpha_2}^-,\delta_{\alpha^\prime_4}\}$ and $\Delta^\prime=\{\delta_{\alpha_1}^+,\delta_{\alpha_2}^+\}$. The quotient by $\Delta^\prime$ is a parabolic induction of two factors of type $aa(1,1)$ and $d(n)$. Let $G$ be $SO(8)\times SO(6+2n)$. The subgroup $(L_1,L_1)$ is isomorphic to $SL(2)\times SO(2n-1)$ and $\dim C_1=5$. We have $\dim H_1-\dim H=6$ and $\dim C=1$, hence $\dim (L_1,L_1)-\dim (L,L)\leq 2$. The subgroup $(L,L)$ is equal to $(L_1,L_1)$ and $\dim H^u_1 - \dim H^u = 2$. 

{\em Existence}. If $n=2$, $\dim P^u_{S^p}+\rank\Sigma=\dim B^u+8$ and $\rank\Delta-\rank\Sigma=1$, hence $S^p=\emptyset$. We have the following cases: $dz(4,3)+aa(1,1)$ and $dy(4,3)+aa(1,1)$. The latter can be excluded by the corresponding uniqueness proof. If $n>2$, $S^p\neq\emptyset$ and there is only one root of type $d_n$ on a bifurcation.

$\mathbf{dz(4,2)+d(n)}$, $n\geq2$. {\em Uniqueness}. Here $\Delta^{\prime\prime}=\{\delta_{\alpha_1}^+,\delta_{\alpha_2}^+,\delta_{\alpha_2}^-,\delta_{\alpha^\prime_3}\}$ and $\Delta^\prime=\{\delta_{\alpha_1}^+,\delta_{\alpha_2}^+\}$. The quotient by $\Delta^\prime$ is a parabolic induction of two factors of type $aa(1,1)$ and $d(n)$. 

$\mathbf{dz^\sim(4+q+3)}$. {\em Uniqueness}. Here $\Delta^{\prime\prime}=\{\delta_{\alpha_1}^+,\delta_{\alpha_2}^+,\delta_{\alpha_2}^-,\delta_{\alpha_4}\}$ and $\Delta^\prime=\{\delta_{\alpha_1}^+,\delta_{\alpha_2}^+\}$. The quotient by $\Delta^\prime$ is a parabolic induction of two factors respectively of type $aa(1,1)$ and $aa(q)$. Let $G$ be $SO(2q+14)$. The subgroup $(L_1,L_1)$ is isomorphic to $SL(2)\times SL(q)$ and $\dim C_1=6$. We have $\dim H_1-\dim H=6$ and $\dim C=2$. The subgroup $(L,L)$ is equal to $(L_1,L_1)$ and $\dim H^u_1 - \dim H^u = 2$.

{\em Existence}. The root $\alpha_4+\ldots+\alpha_{q+3}$ belongs to $\Sigma$. Since $\dim P^u_{S^p}+\rank\Sigma=\dim B^u+(q-1)(q-2)/2+8$ and $\rank\Delta-\rank\Sigma=2$, we get $S^p=\{\alpha_5,\ldots,\alpha_{q+2}\}$. We have the following three cases: $dz^\sim(4+q+3)$, $dy(4+q+3)$ and $dy^\sim(4+q+3)$.

$\mathbf{dz^\sim(4+q+2)}$. {\em Uniqueness}. Set $\Delta^{\prime\prime}=\{\delta_{\alpha_1}^+,\delta_{\alpha_1}^-,\delta_{\alpha_2}^+,\delta_{\alpha_3}\}$. The distinguished subset $\Delta^\prime=\{\delta_{\alpha_1}^+,\delta_{\alpha_2}^+\}$ is such that the quotient spherical system is a parabolic induction of two factors respectively of type $aa(1,1)$ and $aa(q)$. 

$\mathbf{dz(3+q+2)}$. {\em Uniqueness}. Take $\Delta^{\prime\prime}=\{\delta_{\alpha_1}^+,\delta_{\alpha_2}^+,\delta_{\alpha_2}^-,\delta_{\alpha_3}\}$. The distinguished subset $\Delta^\prime=\{\delta_{\alpha_1}^+,\delta_{\alpha_2}^+\}$ is such that the quotient spherical system is a parabolic induction of $aa(1+q+1)$. Let $G$ be $SO(2q+10)$. The subgroup $(L_1,L_1)$ is isomorphic to $SL(2)\times SL(q+1)$ and $\dim C_1=4$. We have $\dim H_1-\dim H=4$ and $\dim C=2$. The subgroup $(L,L)$ is equal to $(L_1,L_1)$ and $\dim H^u_1-\dim H^u=2$. 

{\em Existence}. The root $\alpha_3+\ldots+\alpha_{q+2}$ belongs to $\Sigma$. We get two cases: $dy(3+q+2)$ and $dz(3+q+2)$. 

$\mathbf{dz(3+q+1)}$. {\em Uniqueness}. Take $\Delta^{\prime\prime}=\{\delta_{\alpha_1}^+,\delta_{\alpha_1}^-,\delta_{\alpha_2},\delta_{\alpha_{q+2}}^+\}$. The distinguished subset $\Delta^\prime=\{\delta_{\alpha_1}^+,\delta_{\alpha_{q+2}}^+\}$ is such that the quotient spherical system is a parabolic induction of $aa(1+q+1)$. 

$\mathbf{az(3,2)+2\!\!-\!\!comb}$. {\em Uniqueness}. Set $\Delta^{\prime\prime}=\{\delta_{\alpha_1}^-,\delta_{\alpha_2}^+,\delta_{\alpha_2}^-,\delta_{\alpha_4}^-\}$. The distinguished subset $\Delta^\prime=\{\delta_{\alpha_1}^-,\delta_{\alpha_2}^-,\delta_{\alpha_4}^-\}$ is such that the quotient spherical system is a parabolic induction of $aa(1,1)$. Let $G$ be $SO(14)$. The subgroup $(L_1,L_1)$ is isomorphic to $SL(2)$ and $\dim C_1=5$. We have $\dim H_1-\dim H=6$, $\dim C=2$, $(L,L)=(L_1,L_1)$ and $\dim H^u_1-\dim H^u=3$.

{\em Existence}. From $\dim P^u_{S^p}+\rank\Sigma=\dim B^u+7$ and $\rank\Delta-\rank\Sigma=2$ we get $S^p=\emptyset$, $\Sigma=S$. We have two strongly $\Delta$-connected components. The proof can be concluded by noticing that $az(3,2)+(2\mathrm{-comb})$ is the unique spherical system admitting $(S^p,\Sigma,\mathbf A)/\Delta^{\prime\prime}$ as quotient.

$\mathbf{ae_6(6)+d(n)}$, $n\geq2$. {\em Uniqueness}. Take $\Delta^{\prime\prime}=\{\delta_{\alpha_1}^+,\delta_{\alpha_2}^+,\delta_{\alpha_2}^-,\delta_{\alpha_7}\}$. The distinguished subset $\Delta^\prime=\{\delta_{\alpha_1}^+,\delta_{\alpha_2}^+\}$ is such that the quotient spherical system is a parabolic induction of two factors of type $aa(1,1)$ and $d(n)$. Let $G$ be $SO(12+2n)$. The subgroup $(L_1,L_1)$ is isomorphic to $SL(2)\times SO(2n-1)$ and $\dim C_1=4$. We have $\dim H_1-\dim H=5$ and $\dim C=1$, hence $\dim (L_1,L_1)-\dim (L,L)\leq 2$. The subgroup $(L,L)$ is equal to $(L_1,L_1)$ and $\dim H^u_1 - \dim H^u = 2$. 

{\em Existence}. If $n=2$, we have $\dim P^u_{S^p}+\rank\Sigma=\dim B^u+7$ and $\rank\Delta-\rank\Sigma=1$, hence $S^p=\emptyset$. We have the following cases: $ae_6(6)+aa(1,1)$ and $ae_7(6)+aa(1,1)$. If $n>2$, $S^p\neq\emptyset$ and there is a root of type $d_n$ on the bifurcation. 

$\mathbf{ae_6(5)+d(n)}$, $n\geq2$. {\em Uniqueness}. Take $\Delta^{\prime\prime}=\{\delta_{\alpha_1}^+,\delta_{\alpha_1}^-,\delta_{\alpha_2}^+,\delta_{\alpha_6}\}$. The distinguished subset $\Delta^\prime=\{\delta_{\alpha_1}^+,\delta_{\alpha_2}^+\}$ is such that the quotient spherical system is a parabolic induction of two factors of type $aa(1,1)$ and $d(n)$. 

$\mathbf{ae_7(7)+d(n)}$, $n\geq2$. {\em Uniqueness}. Take $\Delta^{\prime\prime}=\{\delta_{\alpha_1}^+,\delta_{\alpha_2}^+,\delta_{\alpha_2}^-,\delta_{\alpha_3}^+,\delta_{\alpha_3}^-,\delta_{\alpha_8}\}$. The distinguished subset $\Delta^\prime=\{\delta_{\alpha_1}^+,\delta_{\alpha_2}^+,\delta_{\alpha_3}^+\}$ is such that the quotient spherical system is a parabolic induction of two factors respectively of type $aa(2,2)$ and $d(n)$. Let $G$ be $SO(14+2n)$. The subgroup $(L_1,L_1)$ is isomorphic to $SL(3)\times SO(2n-1)$ and $\dim C_1=3$. We have $\dim H_1-\dim H=5$ and $\dim C=1$, hence $\dim (L_1,L_1)-\dim (L,L)\leq 3$. The subgroup $(L,L)$ is equal to $(L_1,L_1)$ and $\dim H^u_1 - \dim H^u = 3$. 

{\em Existence}. If $n=2$, $\dim P^u_{S^p}+\rank\Sigma=\dim B^u+8$ and $\rank\Delta-\rank\Sigma=1$, hence $S^p=\emptyset$. It remains only the case of $ae_7(7)+aa(1,1)$. If $n>2$, $S^p\neq\emptyset$ and there is a root of type $d_n$ on the bifurcation.

$\mathbf{ae_7(6)+d(n)}$, $n\geq2$. {\em Uniqueness}. Take $\Delta^{\prime\prime}=\{\delta_{\alpha_1}^+,\delta_{\alpha_1}^-,\delta_{\alpha_2}^+,\delta_{\alpha_2}^-,\delta_{\alpha_4}^+,\delta_{\alpha_7}\}$. The distinguished subset $\Delta^\prime=\{\delta_{\alpha_1}^+,\delta_{\alpha_2}^+,\delta_{\alpha_4}^+\}$ is such that the quotient spherical system is a parabolic induction of two factors respectively of type $aa(2,2)$ and $d(n)$. 

$\mathbf{3\!\!-\!\!comb+aa(q)}$. {\em Uniqueness}. Set $\Delta^{\prime\prime}=\{\delta_{\alpha_1}^+,\delta_{\alpha_1}^-,\delta_{\alpha_2}\}$. The distinguished subset $\Delta^\prime=\{\delta_{\alpha_1}^+,\delta_{\alpha_2}\}$ is such that the quotient spherical system is a parabolic induction of $aa((q+1))$. Let $G$ be $SO(2q+6)$. The subgroup $(L_1,L_1)$ is isomorphic to $SL(q+1)$ and $\dim C_1=3$. We have $\dim H_1-\dim H=q+2$ and $\dim C=2$. The subgroup $(L,L)$ is equal to $(L_1,L_1)$.

{\em Existence}. We have $\dim P^u_{S^p}+\rank\Sigma=\dim B^u -(q-1)(q-2)/2+4$ and $\rank\Delta-\rank\Sigma=2$. We know that $H\subset H_1$ and that the quotient $G/H\to G/H_1$ extends to a dominant $G$-morphism with connected fibers. If $q>2$, there are no root of type $d_m$ and there is exactly one root of type $a_q$. The rest is a strongly $\Delta$-connected component of roots of type $a_1$. If $q\leq2$, we know analogously that there are two strongly $\Delta$-connected components.

$\mathbf{3\!\!-\!\!comb+2\!\!-\!\!comb}$. {\em Uniqueness}. Set $\Delta^{\prime\prime}=\{\delta_{\alpha_1}^-,\delta_{\alpha_2}^+,\delta_{\alpha_2}^-\}$. The distinguished subset $\Delta^\prime=\{\delta_{\alpha_1}^-,\delta_{\alpha_2}^+\}$ is such that the quotient spherical system is a parabolic induction of $aa(2)$. Let $G$ be $SO(10)$. The subgroup $(L_1,L_1)$ is isomorphic to $SL(2)$ and $\dim C_1=4$. We have $\dim H_1-\dim H=4$ and $\dim C=2$. The subgroup $(L,L)$ is equal to $(L_1,L_1)$. 

{\em Existence}. We have that $S^p=\emptyset$ and there are two strongly $\Delta$-connected components. The spherical system $(3\mathrm{-comb})+(2\mathrm{-comb})$ is the unique that admits $(S^p,\Sigma,\mathbf A)/\Delta^{\prime\prime}$ as quotient.

$\mathbf{ax(1+p+1)+d(n)}$, $n\geq2$. {\em Uniqueness}. Set $\Delta^{\prime\prime}=\Delta\setminus\{\delta_{\alpha_{p+1}},\delta_{\alpha_{p+2}}^-\}$. The distinguished subset $\Delta^\prime=\{\delta_{\alpha_1}^+,\delta_{\alpha_2}\}$ is such that the quotient spherical system is a parabolic induction of two factors respectively $aa((p+1))$ and $d(n)$. Let $G$ be $SO(2p+4+2n)$. The subgroup $(L_1,L_1)$ is isomorphic to $SL(p+1)\times SO(2n-1)$ and $\dim C_1=2$. We have $\dim H_1-\dim H=p+1$ and $\dim C=2$. The subgroup $(L,L)$ is equal to $(L_1,L_1)$. 

{\em Existence}. We have $\dim P^u_{S^p}+\rank\Sigma=\dim B^u -(p-1)(p-2)/2 -(n-1)(n-2) +4$ and $\rank\Delta-\rank\Sigma=2$. We know that $H\subset H_1$ and that the quotient $G/H\to G/H_1$ extends to a dominant $G$-morphism with connected fibers. If $q>2$, there is exactly one root of type $d_n$ on the bifurcation and one root of type $a_q$. The rest is a strongly $\Delta$-connected component of roots of type $a_1$. If $q\leq2$, we know analogously that there are three strongly $\Delta$-connected components.

$\mathbf{af(4)+d(n)}$, $n\geq2$. {\em Uniqueness}. Set $\Delta^{\prime\prime}=\{\delta_{\alpha_1}^+,\delta_{\alpha_2}^+,\delta_{\alpha_2}^-,\delta_{\alpha_5}\}$. The distinguished subset $\Delta^\prime=\{\delta_{\alpha_1}^+,\delta_{\alpha_2}^+\}$ is such that the quotient spherical system is a parabolic induction of two factors respectively $aa(2)$ and $d(n)$. Let $G$ be $SO(8+2n)$. The subgroup $(L_1,L_1)$ is isomorphic to $SL(2)\times SO(2n-1)$ and $\dim C_1=3$. We have $\dim H_1-\dim H=3$, $\dim C=2$ and $\dim (L_1,L_1)-\dim (L,L)\leq 2$. The subgroup $(L,L)$ is equal to $(L_1,L_1)$. 

{\em Existence}. We have $\dim P^u_{S^p}+\rank\Sigma=\dim B^u -(n-1)(n-2) +5$ and $\rank\Delta-\rank\Sigma=2$. We know that $H\subset H_1$ and that the quotient $G/H\to G/H_1$ extends to a dominant $G$-morphism with connected fibers. There is exactly one root of type $d_n$ on the bifurcation and two strongly $\Delta$-connected components of roots of type $a_1$. 


\section{Normality of the Demazure embedding}

In this last section we show that the proof of D.~Luna (\cite{L02}) of the conjecture of M.~Brion (\cite{B90}) on the normality of the Demazure embedding holds also in type $\mathsf A$~$\mathsf D$.

Let us report some notations from \cite{L02}. Let $\mathfrak g$ be a semisimple complex Lie algebra, and let $\mathfrak h$ be a Lie subalgebra of $\mathfrak g$ equal to its normalizer. Moreover, $G$ denotes the connected component of the identity in $\mathrm{Aut}(\mathfrak g)$ and $H$ the stabilizer of $\mathfrak h$ in $G$. If $m=\dim\mathfrak h$, $\mathrm{Gr}_m(\mathfrak g)$ denotes the grassmannian of $m$-dimensional subspaces of $\mathfrak g$. The Demazure embedding of $G/H$ is the orbit closure $\overline{G.\mathfrak h}$ in $\mathrm{Gr}_m(\mathfrak g)$.

\subsubsection*{Rigid wonderful varieties}

A wonderful $G$-variety is said to be {\em rigid} if it is an embedding of $G/H$ where $H=N_G(H)$. In case of an adjoint semisimple group $G$ of type $\mathsf A$~$\mathsf D$, the rigidity of a wonderful variety can be characterized, as for type $\mathsf A$, using the classification:

\begin{proposition}\label{rigidity} Let $G$ be an adjoint semisimple group of type $\mathsf A$~$\mathsf D$, a wonderful $G$-variety is rigid if and only if the set $\mathbf A_X$ of its spherical system $(S^p_X,\Sigma_X,\mathbf A_X)$ does not contain two elements giving the same functional on $\Xi_X$.\end{proposition}

\begin{proof} If $N_G(H)$ and $H$ are different then they are the stabilizers of two points in the open $G$-orbits of two different wonderful $G$-varieties $X^\prime$ and $X$ associated to different spherical systems.

The normalizer $N_G(H)$ acts on the set of colours of $X$, and this action is non-trivial because otherwise the spherical system of $X^\prime$ would be equal to the one of $X$. This comes from the fact that $S^p_{X^\prime}=S^p_X$ and $\mathbf A_{X^\prime}=\mathbf A_X$, in type $\mathsf A$ and $\mathsf D$, implies $\Sigma_{X^\prime}=\Sigma_X$.

But now the action of $N_G(H)$ can only exchange couples of colours $X$ moved by the same set of simple (spherical) roots and giving the same functional on $\Xi_X$.

Vice versa, the spherical systems with two elements of $\mathbf A$ giving the same functional can be easily treated. They are parabolic inductions of spherical systems with at least one irreducible component equal to one of the following three spherical systems: $aa(p+1+p)$, $D1(i)$, $dc(n)$ for $n$ even. Each one of them is associated to a nonrigid wonderful variety, embedding of $G/H$ that is respectively equal to $SL(2p+2)/S(GL(p+1)\times GL(p+1))$, $SO(2n)/(GL(1)\times SO(2n-2))$, $SO(2n)/GL(n)$. Therefore, every spherical system with two elements of $\mathbf A$ with the same functional is associated to a nonrigid wonderful variety.\end{proof}

\subsubsection*{Normality}

\begin{theorem}Let $\mathfrak g$ be of type $\mathsf A$~$\mathsf D$, for every spherical Lie subalgebra $\mathfrak h$ equal to its normalizer, the Demazure embedding $\overline{G.\mathfrak h}$ is wonderful.\end{theorem}

It is known that under these hypotheses, for $\mathfrak g$ of any type, the normalization of the Demazure embedding is wonderful (\cite{B90}, \cite{Kn96}). Therefore, for us it is sufficient to prove that the Demazure embedding itself is normal.

As in \cite{L02} the matter can be reformulated as follows. Let $X$ be a wonderful $G$-variety, embedding of $G/H$, $m=\dim H$. The morphism $x\mapsto \mathfrak g_x$, where $x\in G/H$ and $\mathfrak g_x$ is the Lie algebra of the stabilizer of $x$, extends to a $G$-morphism $\delta_X\colon X\rightarrow \mathrm{Gr}_m(\mathfrak g)$ (the Demazure morphism). If $X$ is a rigid wonderful variety, the morphism $\delta_X\colon X\rightarrow\delta_X(X)$ is birational and finite (\cite{B90}). If $\delta_X(X)$ is normal, then $\delta_X\colon X\rightarrow\delta_X(X)$ is an isomorphism.

Let $z\in X$ denote the unique point fixed by $B^-$, and let $T_z(X)_\gamma$ be the $\gamma$-eigenspace inside the tangent space of $X$ at $z$, for $\gamma$ any spherical root. In \cite{L02} it is shown that the normality of $\delta_X(X)$ can be proven checking that the differential of $\delta_X$, restricted to $T_z(X)_\gamma$, is injective for all $\gamma$.

In \cite{L02} the proof goes on as follows. Let $G$ be an adjoint semisimple group of type $\mathsf A$. A wonderful $G$-variety is said to be {\em critical} if it is rigid, it is not a parabolic induction (i.e.\ $\supp(\Sigma)\cup S^p=S$) and it is either of rank $1$ or of rank $2$ with at least one simple spherical root (i.e.\ $\Sigma\cap S\neq\emptyset$).
\begin{enumerate}\renewcommand{\theenumi}{\roman{enumi}}
\item\label{a2} Let $X$ be a rigid wonderful $G$-variety. For every spherical root $\gamma$ there exists a subset $S^\prime$ of simple roots such that $\supp(\gamma)\cup S^p\subset S^\prime$ and the localization of $X$ in $S^\prime$ is critical.
\item\label{a3} If the theorem holds for the critical varieties then it holds in general.
\item\label{a4} The theorem holds for the critical varieties.
\end{enumerate}

The proof of \ref{a3} in \cite{L02} remains valid for the case $\mathsf A$~$\mathsf D$. We need only to generalize the proofs of \ref{a2} and \ref{a4}.

\textit{Proof of \ref{a2}.} If $\gamma\notin S$ then it is sufficient to consider the localization in $\supp(\gamma)$. We have five cases, corresponding to rank one spherical systems, four of which occurring in type $\mathsf A$ (spherical roots of type $a_m$ $m\geq2$, $a^\prime_1$, $a_1\times a_1$, $d_3$) plus one corresponding to a spherical root of type $d_m$ $m\geq4$. Each of these spherical systems is associated to a rigid wonderful variety.

If $\gamma=\alpha\in\Sigma\cap S$ then there exists a spherical root $\gamma^\prime\neq\alpha$ such that $\langle\rho(\delta_\alpha^+),\gamma^\prime\rangle\neq\langle\rho(\delta_\alpha^-),\gamma^\prime\rangle$. Therefore, it is sufficient to consider the localization in $\supp\{\alpha,\gamma^\prime\}$. From the list of rank two prime wonderful varieties (\cite{Wa96}), we have five cases, four of which occurring in type $\mathsf A$ plus one corresponding to $\gamma^\prime$ of type $d_m$ $m\geq3$, $\langle\rho(\delta_\alpha^+),\gamma^\prime\rangle=0$ and $\langle\rho(\delta_\alpha^-),\gamma^\prime\rangle=-2$, namely the rank two spherical system $D4(ii)$. Each of these spherical systems is associated to a rigid wonderful variety.

\textit{Proof of \ref{a4}.} The rank one cases are symmetric varieties, hence the theorem follows from \cite{DP83}. It remains to deal with the above critical cases of rank $2$, $\Sigma=\{\alpha,\gamma^\prime\}$, and prove the injectivity of the differential of $\delta_X$ only on $T_z(X)_\alpha$. In \cite{L02} this is reduced to verify that $s_\alpha\cdot(\gamma^\prime-\langle\rho(\delta),\gamma^\prime\rangle\alpha)\neq(\gamma^\prime-\langle\rho(\delta),\gamma^\prime\rangle\alpha)$, for $\delta\in\mathbf A(\alpha)$.

Hence, we have only to check the last of the five above rank two spherical systems: $\Sigma=\{\alpha=\alpha_1,\gamma^\prime=2\alpha_2+\ldots+2\alpha_{n-2}+\alpha_{n-1}+\alpha_n\}$. Here we have $\gamma^\prime-\langle\rho(\delta_\alpha^+),\gamma^\prime\rangle\alpha=\gamma^\prime$ and $\gamma^\prime-\langle\rho(\delta_\alpha^-),\gamma^\prime\rangle\alpha=\gamma^\prime+2\alpha$, while $s_\alpha\cdot\gamma^\prime=\gamma^\prime-\langle\alpha^\vee,\gamma^\prime\rangle\alpha=\gamma^\prime+2\alpha$ and $s_\alpha\cdot(\gamma^\prime+2\alpha)=\gamma^\prime$.

\end{document}